\documentclass[11pt,a4paper]{article}
\usepackage{amsfonts,amsmath,amssymb,accents,amsthm,bbm}
\usepackage{verbatim,bm}
\usepackage{hyperref}
\usepackage{tikz}
\usetikzlibrary{calc,arrows,positioning,fit,petri,external}
\usepackage[normalem]{ulem}%for \sout,\xout

\newcommand{\inputtikz}[1]{\includegraphics{#1}}
\newcommand{\dis}{\displaystyle}
\newcommand{\txt}{\textstyle}

\newcommand{\noi}{\noindent}
\newcommand{\halmos}{\rule{1ex}{1.4ex}}
\newcommand{\QED}{\nopagebreak{\hspace*{\fill}$\halmos$\medskip}}

\newcommand{\med}{\medskip}
\newcommand{\quand}{\quad\mbox{and}\quad}

\newtheoremstyle{mythm}
{}
{}
{\itshape}
{}
{\bfseries}
{}
{.5em}
{#1 #2 \thmnote{(#3)}}

\theoremstyle{mythm}

\newtheorem{theorem}{Theorem}
\newtheorem{proposition}[theorem]{Proposition}
\newtheorem{lemma}[theorem]{Lemma}
\newtheorem{conjecture}[theorem]{Conjecture}
\newtheorem{remark}[theorem]{Remark}
\newtheorem{defi}[theorem]{Definition}
\newcommand{\bt}{\begin{theorem}}
	\newcommand{\et}{\end{theorem}}
\newcommand{\bl}{\begin{lemma}}
	\newcommand{\el}{\end{lemma}}
\newcommand{\bp}{\begin{proposition}}
	\newcommand{\ep}{\end{proposition}}
\newcommand{\br}{\begin{remark}}
	\newcommand{\er}{\end{remark}}
\newcommand{\bcon}{\begin{conjecture}}
	\newcommand{\econ}{\end{conjecture}}

\newenvironment{Proof}[1][]{\noi\textbf{Proof #1}}{\QED}
\newcommand{\bpro}{\begin{Proof}}
	\newcommand{\epro}{\end{Proof}}

\newcommand{\be}{\begin{equation}}
\newcommand{\ee}{\end{equation}}
\newcommand{\ba}{\begin{array}}
	\newcommand{\ea}{\end{array}}
\newcommand{\bac}{\begin{array}{r@{\,}c@{\,}l}}
	\newcommand{\bc}{\be\begin{array}{r@{\,}c@{\,}l}}
		\newcommand{\ec}{\end{array}\ee}

	\expandafter\mathchardef\expandafter\varphi\number\expandafter\phi\expandafter\relax
	\expandafter\mathchardef\expandafter\phi\number\varphi
	
	\newcommand{\al}{\alpha}
	\newcommand{\bet}{\beta}
	\newcommand{\ga}{\gamma}
	\newcommand{\de}{\delta}
	\newcommand{\De}{\Delta}
	\newcommand{\eps}{\varepsilon}
	\newcommand{\la}{\lambda}
	\newcommand{\La}{\Lambda}
	\newcommand{\sig}{\sigma}
	\newcommand{\tet}{\theta}
	\newcommand{\phh}{\varphi}

	\newcommand{\Phq}{\Phi}
	\newcommand{\psib}{{\bm{\psi}}}
	
	\newcommand{\Ai}{{\cal A}}
	\newcommand{\Di}{{\cal D}}
	\newcommand{\Ei}{{\cal E}}
	\newcommand{\Fi}{{\cal F}}
	
	\newcommand{\Hi}{{\cal H}}
	\newcommand{\Ki}{{\cal K}}
	\newcommand{\Ti}{{\cal T}}
	
	\newcommand{\Vi}{{\cal V}}
	\newcommand{\Wi}{{\cal W}}
	
	\newcommand{\N}{{\mathbb N}}
	\renewcommand{\P}{{\mathbb P}}
	\newcommand{\R}{{\mathbb R}}
	\newcommand{\Z}{{\mathbb Z}}
	\newcommand{\abf}{\mathbf{a}}
	
	\newcommand{\E}{{\mathbb E}}
	\newcommand{\li}{\langle}
	\newcommand{\re}{\rangle}
	
	\newcommand{\desd}{\ensuremath{\Leftrightarrow}}
	
	\newcommand{\sub}{\subset}
	\newcommand{\beh}{\backslash}

	\newcommand{\ti}{\tilde}
	\newcommand{\ov}{\overline}
	\newcommand{\un}{\underline}
	
	\newcommand{\ffrac}[2]{{\textstyle\frac{{#1}}{{#2}}}}
	\newcommand{\ha}{\ffrac{1}{2}}
	\newcommand{\expo}{\mbox{\large\it e}}
	\newcommand{\ex}[1]{\expo^{\,\textstyle{#1}}}
	\newcommand{\nab}[1]{\nabla_{\!#1}}
	\newcommand{\dia}{\#}
	\newcommand{\dgg}{\dagger}
	
	\newcommand{\rbf}{\mathbf{r}}
	
\newcommand{\lis}{1\leq s\leq\sig}
\newcommand{\Dec}{\mathbf{D}}

\newcommand{\Aa}{{\rm A}}
\newcommand{\Ll}{{\rm L}}

\begin{document}

\makeatletter\@addtoreset{equation}{section}
\makeatother\def\theequation{\thesection.\arabic{equation}} 

\renewcommand{\labelenumi}{{\rm (\roman{enumi})}}
\renewcommand{\theenumi}{\roman{enumi}}

\title{Peierls bounds from random Toom contours}
\author{Jan M.~Swart\footnote{The Czech Academy of Sciences, Institute of Information Theory and Automation, Pod vod\'arenskou v\v{e}\v{z}\'i 4, 18200 Praha 8, Czech republic. swart@utia.cas.cz},
R\'eka Szab\'o\footnote{Bernoulli Institute, University of Groningen, Nijenborgh 9, 9747 AG Groningen, The Netherlands. r.szabo@rug.nl},
and Cristina Toninelli\footnote{CEREMADE, CNRS, Universit\'e Paris-Dauphine, PSL University, Place du Mar\'echal de Lattre de Tassigny, 75016 Paris, France and DMA, Ecole normale sup\'erieure, PSL University,  45 rue d’Ulm, F-75230 Cedex 5 Paris, France. toninelli@ceremade.dauphine.fr}}

\date{\today}

\maketitle

\begin{abstract}\noi
For deterministic monotone cellular automata on the $d$-dimensional integer lattice, Toom (1980) has given necessary and sufficient conditions for the all-one fixed point to be stable against small random perturbations. We are interested in the open problem of extending Toom's result to monotone cellular automata with intrinsic randomness, where the unperturbed evolution is random with i.i.d.\ update rules attached to the space-time points. For some applications it is also desirable to consider a more general graph structure, so we assume that the underlying lattice is an arbitrary countable group. Toom's proof of stability is based on a Peierls argument. In previous work, we demonstrated that this Peierls argument can also be used to prove stability for cellular automata with intrinsic randomness, but in this case estimating the Peierls sum becomes much harder than in the deterministic case. In the present paper, we develop a method based on random contours to estimate the Peierls sum and apply it to prove new stability results for monotone cellular automata with intrisic randomness. We also demonstrate the limitations of the method by constructing an example where the Peierls sum is infinite for arbitrary small perturbations even though stability is believed to hold.
\end{abstract}
\vspace{.5cm}

\noi
{\it Acknowledgments.} The first author is supported by grant 20-08468S of the Czech Science Foundation (GA\v{C}R). The second and third authors are supported by ERC Starting Grant 680275 ``MALIG''.

\newpage

{\setlength{\parskip}{-2pt}\tableofcontents}

\newpage

\section{Introduction and main results}\label{S:intro}

\subsection{Monotone cellular automata}\label{S:auto}

Let $\La$ be a countable group with product denoted by $(i,j)\mapsto ij$ and unit element $0$. We will especially be interested in the case that $\La=\Z^{d+1}$ equipped with the additive group structure $ij:=i+j$ $(i,j\in\Z^{d+1})$, but for some purposes it is useful to be more general. Denote by $\{0,1\}^{\La}$ the set of configurations $x=((x(i))_{i\in\La}$ of zeros and ones on $\La$. By definition, a map $\phh:\{0,1\}^\La\to\{0,1\}$ is \emph{local} if $\phh$ depends only on finitely many coordinates, i.e., there exists a finite set $\De\sub\La$ and a function $\phh':\{0,1\}^\De\to\{0,1\}$ such that $\phh\big((x(i))_{i\in\La}\big)=\phh'\big((x(i))_{i\in\De}\big)$ for each $x\in\{0,1\}^\La$. We let $\De(\phh)$ denote the smallest set with this property, which may be empty: in this case $\phh$ is constantly zero or one. We let
\be\label{phh01}
\phh^0(x):=0\quand\phh^1(x):=1\qquad(x\in\{0,1\}^\La)
\ee
denote these constant functions. A local map $\phh$ is \emph{monotone} if $x\leq y$ (coordinatewise) implies $\phh(x)\leq\phh(y)$. We will use the word \emph{height function} for a group homomorphism from $\La$ to $(\Z,+)$, i.e., a function $h:\La\to\Z$ such that 
\be\label{hfnc}
h(0)=0 \quad\mbox{and}\quad h(ij)=h(i)+h(j)\quad(i,j\in\La).
\ee
We will in particular be interested in the case that
\be\label{Zd}
\La=\Z^{d+1}\quand h(i)=h(i_1,\ldots,i_{d+1})=-i_{d+1}\qquad(i\in\Z^{d+1}).
\ee
We define a \emph{monotone cellular automaton} to be a collection of monotone local maps $\Phq=(\Phq_i)_{i\in\La}$ that satisfy
\be\label{phhDe}
\De(\Phq_i)\sub\{j\in\La:h(j)=h(i)+1\}\qquad(i\in\La)
\ee
for some height function $h$. A \emph{trajectory} of a cellular automaton is a configuration $x\in \{0,1\}^{\La}$ such that
\be
x(i)=\Phq_i(x)\qquad(i\in\La).
\ee
Note that for each trajectory $x$ of $\Phq$ and for each $t\in \Z$ the values of $\{x(i): i\in \Lambda, h(i)=t\}$ are completely determined by the values of $\{x(i): i\in\Lambda, h(i)=t+1\}$. Therefore, we can think of the height function as negative time and $\{x(i): h(i)=-t\}$ as the state of the cellular automata at time $t$. We cite the following elementary lemma from \cite[Lemma~14]{SST24}.

\bl[Minimal and maximal trajectories]
Let\label{L:maxtraj} $\Phq$ be a monotone cellular automaton. Then there exist trajectories $\un x$ and $\ov x$ that are uniquely characterised by the property that each trajectory $x$ of $\Phq$ satisfies $\un x\leq x\leq\ov x$ (pointwise).
\el

\subsection{Stability}\label{S:stable}

We will be interested in the maximal trajectory of random monotone cellular automata on an arbitrary countable group $\Lambda$ with height function $h$. Let $\{\phh_0,\ldots,\phh_m\}$ be a collection of local monotone maps $\phh_k:\{0,1\}^\La\to\{0,1\}$ of which $\phh_0=\phh^0$ is the map that is constantly zero while $\phh_k$ is not constant and satisfies
\be\label{phiDe}
\De(\phh_k)\sub\{i\in\La:h(i)=1\}\quad(1\leq k\leq m).
\ee
Furhermore, let $\rbf=\big(\rbf(1),\ldots,\rbf(m)\big)$ be a probability distribution on $\{1,\ldots,m\}$ and for $p\in[0,1]$ let $\big(\mu(i)\big)_{i\in\La}$ be i.i.d.\ random variables with values in $\{0,\ldots,m\}$ such that
\be\label{Phip}
\P\big[\mu(i)=k\big]=\left\{\ba{ll}
\dis p\quad&\mbox{if }k=0,\\[5pt]
\dis (1-p)\rbf(k)\quad&\mbox{if }1\leq k\leq m,
\ea\right. \qquad(i\in\La).
\ee
We then define a \emph{random monotone cellular automaton} $\Phi^p=(\Phi^p_i)_{i\in\La}$ by
\be\label{PhH}
\Phi^p_i(x):=\phh_{\mu(i)}\big((x(ij))_{j\in\La}\big)\qquad\big(i\in\La,\ x\in\{0,1\}^\La\big).
\ee
Note that (\ref{hfnc}) and (\ref{phiDe}) imply that the functions $(\Phi^p_i)_{i\in\La}$ satisfy (\ref{phhDe}).

We see from \eqref{Phip} that $p$ is the probability with which the zero-map is applied in a given space-time point. We will mainly be interested in the case that $p$ is small but positive. We think of $\Phi^p$ as a small perturbation of $\Phi^0$. In the special case that $m=1$, we say that $\Phi^0$ is a \emph{deterministic} monotone cellular automaton. If $m\geq 2$ and $\rbf(k)<1$ for all $k$, then we say that $\Phi^0$ has \emph{intrinsic randomness}. Let $\ov X^p$ denote the maximal trajectory of $\Phi^p$. It follows from our definitions that $\big(\ov X^p(i)\big)_{i\in\La}$ is equal in law to $\big(\ov X^p(ji)\big)_{i\in\La}$ for each $j\in\La$. As a result, the \emph{density of the maximal trajectory}
\be\label{ovrho}
\ov\rho(p):=\P\big[\ov X^p(i)=1\big]\qquad(i\in\La)
\ee
does not depend on $i\in\La$. In the special setting of (\ref{Zd}), one can alternatively interpret $\ov\rho(p)$ as the density of the upper invariant law of a Markov chain with state space $\{0,1\}^{\Z^d}$, see \cite[Lemma~15]{SST24}. A simple coupling argument shows that $\ov\rho(p)$ is a nonincreasing function of $p$. The parameter
\be
p_{\rm c}:=\sup\{p\geq 0:\ov\rho(p)>0\big\}
\ee
is called the \emph{critical noise parameter}. By definition, $\Phi^0$ is \emph{stable} if
\be\label{stable}
\lim_{p\to 0}\ov\rho(p)=1,
\ee
and \emph{completely unstable} if $\ov\rho(p)=0$ for all $p>0$, or equivalently $p_{\rm c}=0$.

Let $\phh$ be a local monotone map satisfying (\ref{phiDe}). For each $D\sub\La$, let $\psib^D=(\psi^D_i)_{i\in\La}$ be the cellular automaton defined by
\be
\psi^D_i(x):=\left\{\ba{ll}
0\quad&\mbox{if }i\in D,\\[5pt]
\phh\big((x(ij))_{j\in\La}\big)\quad&\mbox{if }i\in\La\beh D
\ea\right.
\qquad\big(i\in\La,\ x\in\{0,1\}^\La\big),
\ee
and let $\ov x^D$ denote its maximal trajectory. By definition, $\phh$ is an \emph{eroder} if the set $\{i\in\La:\ov x^D(i)=0\}$ is finite for each finite $D\sub\La$. We cite the following result from \cite[Thm~5]{Too80} (see also \cite[Thm~2]{SST24}).

\begin{theorem}[Toom's stability theorem]
Assume\label{T:Toom} (\ref{Zd}) and let $\Phi^0$ be defined by a single local monotone map $\phh$ satisfying (\ref{phiDe}). Then $\Phi^0$ is stable if $\phh$ is an eroder and completely unstable otherwise.
\end{theorem}

In the setting where $\La=\Z^{d+1}$ as in (\ref{Zd}), there exists a simple criterion to check whether a local monotone map $\phh$ is an eroder, see \cite[Prop.~3]{SST24}. For deterministic monotone cellular automata on other space-time sets than $\Z^{d+1}$, or for monotone cellular automata with intrinsic randomness, our understanding of stability is much less complete.

At the heart of Toom's proof of his stablity theorem (Theorem~\ref{T:Toom}) lies an intricate Peierls argument. This Peierls argument was reformulated and extended to cellular automata with intrinsic randomness in \cite{SST24}. Without going into the details, the argument gives an estimate of the form
\be\label{Peierls}
\P\big[\ov X^p(0)=0\big]\leq\sum_{T\in\Ti_0}\P\big[T\mbox{ is present in }\Phi^p\big],
\ee
where the sum runs over a certain type of directed graphs called \emph{Toom contours rooted at} $0$, that can be present in $\Phi^p$. (A precise formulation of (\ref{Peierls}) will be given in Theorem~\ref{T:Peierls} below.)

Under the assumptions of Theorem~\ref{T:Toom}, one can show that the Peierls sum on the right-hand side of (\ref{Peierls}) tends to zero as $p\to 0$, proving stability. In more general settings, it can be hard to estimate the Peierls sum from above. A first extension of Toom's stability theorem to some monotone cellular automata with intrinsic randomness was proved in \cite[Thm~9]{SST24} (see Theorem~\ref{T:thm9} below), but this excludes many interesting cases. In the present paper, we significantly improve on this result by developing a more sophisticated method for estimating the Peierls sum.

As a warm-up, in the following three subsections (Subsections \ref{S:NEC}--\ref{S:coop}), we present three new results that we are able to prove with our new methods. These results are consequences of two more abstract and rather technical results that will be stated later as Theorems \ref{T:cycbd} and \ref{T:conbd} in Subsection~\ref{S:abstract} below. In Subsections \ref{S:edge} and \ref{S:discus} we elaborate a bit on our methods and discuss some open problems.

In Section~\ref{S:randcont} we explain our methods in detail. We first recall Toom's Peierls argument in the reformulation of \cite{SST24}, thus giving a precise meaning to the Peierls sum on the right-hand side of (\ref{Peierls}), and then explain our method for deriving upper bounds on the Peierls sum. We also demonstrate the limitations of Toom's Peierls argument by constructing an example where the Peierls sum is infinite for all $p>0$ even though stability is believed to hold. Sections \ref{S:nuproof}--\ref{S:shrink} are devoted to proofs.

\subsection{Toom's rule}\label{S:NEC}

We start by looking at a deterministic monotone cellular automaton for which stability was already proved by Toom in 1974 \cite{Too74}. We consider the cellular automaton $\Phi^p$ of the form (\ref{PhH}) with $\La:=\Z^3$, $m:=1$, and $\phh_1:=\phh$ given by
\be\label{Toomrule}
\phh(x):={\tt round}\big([x(1,0,0)+x(0,1,0)+x(0,0,1)]/3\big),
\ee
where ${\tt round}$ denotes the function that rounds off to the nearest integer. This local map satisfies (\ref{phiDe}) for the height function
\be
h(i):=i_1+i_2+i_3\qquad(i\in\Z^3).
\ee
Since $m=1$ there is no intrinsic randomness and the distribution $\rbf$ is trivial. The function $\phh$ in (\ref{Toomrule}) is known as \emph{Toom's rule}. By a simple transformation of space-time, the cellular automaton $\Phi^p$ can be transformed into a different cellular automaton $\ti\Phi^p$ that is also of the form (\ref{PhH}) with $\La:=\Z^3$ and $m:=1$, but has $\phh_1:=\phh^{\rm NEC}$ given by
\be\label{NEC}
\phh^{\rm NEC}(x):={\tt round}\big([x(1,0,-1)+x(0,1,-1)+x(0,0,-1)]/3\big).
\ee
This so-called \emph{North East Central} (NEC) map satisfies (\ref{phiDe}) for the height function
\be
h(i):=-i_3\qquad(i\in\Z^3).
\ee
It is easy to see that $\Phi^p$ and $\ti\Phi^p$ are related by the transformation of space-time
\be
i=(i_1,i_2,i_3)\mapsto(i_1,i_2,-i_1-i_2-i_3)=:i'\qquad(i\in\Z^3).
\ee
If $\Phi^p$ is defined in terms of i.i.d.\ random variables $(\mu(i))_{i\in\Z^3}$ and $\ti\Phi^p$ is defined in terms of $(\mu'(i))_{i\in\Z^3}$ given by $\mu'(i'):=\mu(i)$ $(i\in\Z^3)$, then a function $x:\Z^3\to\{0,1\}$ is a trajectory of $\Phi^p$ if and only if the function $x'$ defined by $x'(i'):=x(i)$ $(i\in\Z^3)$ is a trajectory of $\ti\Phi^p$. As a result, the density of the maximal trajectory $\ov\rho(p)$ is the same for $\Phi^p$ and $\ti\Phi^p$. The advantage of the local rule (\ref{NEC}) is that it fits the general scheme (\ref{Zd}) but the formulation of Toom's rule in (\ref{Toomrule}) exhibits more symmetry. As an application of our methods, in Subsection~\ref{S:some} we will prove the following result.

\bt[Toom's rule]
For\label{T:NEC} the deterministic monotone cellular automata defined by either of the local rules (\ref{Toomrule}) or (\ref{NEC}) one has $p_{\rm c}\geq 1/12000$.
\et

As far as we know, Toom never derived an explicit bound for $p_{\rm c}$ for this model. By a careful analysis of his method, in \cite[Prop~11]{SST24} we were able to derive the bound $p_{\rm c}\geq 3^{-21}$. Compared to this, Theorem~\ref{T:NEC} is a significant improvement even though it is still very far from the estimate $p_{\rm c}\approx 0.053$ which comes from numerical simulations.\footnote{This estimate, as well as estimates for other models mentioned later, come from simulations we have run that estimate edge speeds as a function of $p$. These sort of estimates tend to be more reliable than those that try to estimate $\ov\rho(p)$ by starting a cellular automaton in the all one initial state and running it for a long time, since the latter tend to overestimate $p_{\rm c}$ due to metastability effects.}

\subsection{A cellular automaton on the triangular lattice}\label{S:triang}

\begin{figure}[htb]
\begin{center}
\inputtikz{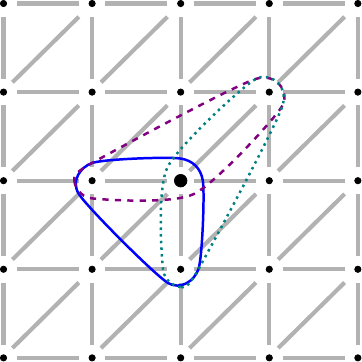}
\caption{A monotone cellular automaton with intrinsic randomness on the triangular lattice. The value of the point in the middle is replaced by the value that holds the majority in a set chosen uniformly at random from the three sets drawn in the picture.}
\label{fig:triang}
\end{center}
\end{figure}

In this subsection we consider the cellular automaton $\Phi^p$ with intrinsic randomness of the form (\ref{PhH}) with $\La:=\Z^3$, $m:=3$, $\rbf$ the uniform distribution on $\{1,2,3\}$, and the local maps $\phh_1,\phh_2,\phh_3$ given by
\bc\label{phitri}
\dis\phh^{\rm SWC}(x)&:=&\dis{\tt round}\big([x(0,0,-1)+x(0,-1,-1)+x(-1,0,-1)]/3\big),\\[5pt]
\dis\phh^{\rm SUC}(x)&:=&\dis{\tt round}\big([x(0,0,-1)+x(0,-1,-1)+x(1,1,-1)]/3\big),\\[5pt]
\dis\phh^{\rm WUC}(x)&:=&\dis{\tt round}\big([x(0,0,-1)+x(-1,0,-1)+x(1,1,-1)]/3\big).
\ec
We call $\phh^{\rm SWC}$ the \emph{South-West-Center voting rule}, $\phh^{\rm SUC}$ the \emph{South-Up-Center voting rule}, and $\phh^{\rm WUC}$ the \emph{West-Up-Center voting rule}. These maps satisfy (\ref{phiDe}) for the height function $h(i_1,i_2,i_3):=-i_3$ and hence fit in the scheme (\ref{Zd}), which means that space-time has the form $\La=\Z^2\times\Z$ where $\Z^2$ is space and $\Z$ is time.

The cellular automaton $\Phi^p$ we have just defined has a nice symmetry that becomes more apparent when we view the space $\Z^2$ as a graph in which each each point $(i_1,i_2)$ is adjacent to six other points, which are $(i_1\pm 1,i_2)$, $(i_1,i_2\pm 1)$, as well as $(i_1-1,i_2-1)$ and $(i_1+1,i_2+1)$. As illustrated in Figure~\ref{fig:triang}, $\Z^2$ equipped with this graph structure is the triangular lattice. As an application of our methods, in Subsection~\ref{S:tripro} we will prove the following result.

\bt[Majority rules on the triangular lattice]
The\label{T:triang} random monotone cellular automaton that applies the maps $\phh^{\rm SWC},\phh^{\rm WUC}$, and $\phh^{\rm SUC}$ with equal probabilities is stable. Moreover, for this cellular automaton, $p_{\rm c}\geq 7.7\cdot 10^{-13}$.
\et

Numerical simulations suggest that for this model the true value of the critical noise parameter is $p_{\rm c}\approx 0.031$. This cellular automaton does not satisfy the conditions of \cite[Thm~9]{SST24} and as far as we know its stability does not follow from any previously known results.

\subsection{Cooperative branching and the identity map}\label{S:coop}

We finally consider the cellular automaton on $\La=\Z^3$ that applies the maps
\bc\label{coopid}
\dis\phh^{\rm coop}(x)&:=&\dis x(0,0,-1)\vee[x(0,1,-1)\wedge x(1,0,-1)],\\[5pt]
\dis\phh^{\rm id}(x)&:=&\dis x(0,0,-1),
\ec
with probabilities $r$ and $1-r$, respectively. In the notation of Subsection~\ref{S:stable}, this means that we set $m:=2$, $\phh_1:=\phh^{\rm coop}$, $\phh_2:=\phh^{\rm id}$, and $\rbf$ is the probability distribution given by $\rbf(1):=r$ and $\rbf(2):=1-r$. This is again a cellular automaton of the form (\ref{Zd}) so we can view $\La$ as the product $\Z^2\times\Z$ of space $\Z^2$ and time $\Z$. The map $\phh^{\rm coop}$ is a \emph{cooperative branching map} where two occupied sites $(0,1)$ and $(1,0)$ can together produce offspring at $(0,0)$. The map $\phh^{\rm id}$ is the \emph{identity map} that just copies the state that was previously at $0$, so for $r=0$ trajectories of the cellular automaton are configurations that are constant in time. We let $\ov\rho(p,r):=\rho_\rbf(p)$ as in (\ref{ovrho}) denote the probability that the maximal trajectory has a one at the origin and write $p_{\rm c}(r)$ for the critical noise parameter in dependence on~$r$.

If we rescale time by a factor $r$ and at the same time send $r$ to zero, then our cellular automaton converges to a continuous-time interacting particle system in which the cooperative branching map is applied at times of a Poisson point process with intensity one. Gray \cite{Gra99} has proved stability for a class of monotone continuous-time interacting particle systems, which includes the one that applies $\phh^{\rm coop}$ at Poissonian times. In view of this, it is natural to conjecture that the cellular automaton defined by $\phh^{\rm coop}$ and $\phh^{\rm id}$ is stable for all $r>0$ and that $p_{\rm c}(r)\sim cr$ as $r\to 0$ for some $c>0$. Our methods allow us to confirm the first conjecture, but we stop short of proving the second one. In Subsection~\ref{S:coopro} we will prove the following result.

\bt[Cooperative branching and the identity map]
One\label{T:coop} has:
\begin{enumerate}
\item $\dis\lim_{p\to 0}\ov\rho(p,r)=1$ for each $r>0$.
\item There exists a constant $c>0$ such that $\ov\rho(p,r)\to 1$ if $p\to 0$ and $r\to 0$ in such a way that $p\leq cr^2$.
\end{enumerate}
\et

For $r=1$, stability for this cellular automaton follows from Toom's stability theorem \cite{Too80} but for $0<r<1$ we are not aware of any previously known results (including our own \cite[Thm~9]{SST24}) that would imply stability. Alternatively to our proof, one may try to apply the (highly nontrivial) methods of \cite{Gra99}. The fact that part~(ii) of Theorem~\ref{T:coop} seems to be suboptimal is an indication that our methods are not optimal. Possible reasons for this and possible ways to improve our methods will be discussed at the end of Subsection~\ref{S:infin} below.

\subsection{Edge speeds}\label{S:edge}

In this subsection, we discuss edge speeds, which are important for an intuitive understanding why certain monotone cellular automata are stable while others are not, and will also allow us to compare our results with earlier work such as our own \cite{SST24}.

To simplify the discussion, we assume for the moment that (\ref{Zd}) holds, so that space-time has the product form $\La=\Z^d\times\Z$ where $\Z^d$ represents space and $\Z$ represents time. We denote elements of $\La$ as $(i,t)$ where $i\in\Z^d$ is the spatial coordinate and $t\in\Z$ is the time coordinate. Then (\ref{Zd}) and (\ref{phiDe}) imply that for each $0\leq k\leq m$, there exists a local monotone function $\phi_k:\{0,1\}^{\Z^d}\to\{0,1\}$ such that (compare (\ref{PhH}))
\be
\Phi^p_{(i,t)}:=\phh_{\mu(i,t)}\big((x(i+j,t+u))_{(j,u)\in\La}\big)
=\phi_{\mu(i,t)}\big((x(i+j,t-1))_{j\in\Z^d}\big)
\ee
for each $(i,t)\in\La$ and $x\in\{0,1\}^\La$. The functions $\phh_k$ and $\phi_k$ are related by
\be\label{phhphi}
\phh_k(x)=\phi_k\big((x(i,-1))_{i\in\Z^d}\big)\qquad\big(x\in\{0,1\}^\La\big).
\ee
In particular, $\phi_0=\phi^0$ is the function that is constantly zero while $\phi_1,\ldots,\phi_m$ are not constant. For each local monotone function $\phi:\{0,1\}^{\Z^d}\to\{0,1\}$, we define a map $\Psi_\phi:\{0,1\}^{\Z^d}\to\{0,1\}^{\Z^d}$ by
\be\label{Psiphi}
\Psi_\phi(x)(i):=\phi\big((x(i+j))_{j\in\Z^d}\big)\qquad\big(i\in\Z^d,\ x\in\{0,1\}^{\Z^d}\big),
\ee
i.e., $\Psi_\phi$ describes one step of the time evolution of the deterministic cellular automaton defined by $\phi$. We say that the ``space'' map $\phi$ is an \emph{eroder} if the associated ``space-time'' map $\phh$ (\ref{phhphi}) is an eroder, in the sense defined at the end of Subsection~\ref{S:stable}. Equivalently, this says that for each configuration $x\in\{0,1\}^{\Z^d}$ that contains finitely many zeros, there exists a $t\in\N$ such that $\Psi_\phi^t(x)$ (the $t$ times iterated map applied to $x$) is the all-one configuration.

For any linear form $\ell:\R^d\to\R$ that is not identically zero and $r\in\R$, let $H^\ell_r\in\{0,1\}^{\Z^d}$ denote the half-space configuration defined by
\be\label{halfspace}
H^\ell_r(i):=\left\{\ba{ll}
1\quad&\mbox{if }\ell(i)\geq r,\\[5pt]
0\quad&\mbox{if }\ell(i)<r
\ea\right.
\ee
It turns out that deterministic cellular automata map half-space configurations into half-space configurations. To formulate this, we need one more definition. A \emph{one-set} of a monotone local map $\phh:\{0,1\}^{\La}\to\{0,1\}$ is a finite set $A\sub\La$ such that $\phh(1_A)=1$, where $1_A$ is the indicator function of set $A$. A \emph{minimal one-set} is a one-set that does not contain other one-sets as a proper subset. We let $\Ai(\phh)$ denote the set of all minimal one-sets of $\phh$. Note that by minimality $A\sub\De(\phh)$ for all $A\in\Ai(\phh)$. We cite the following fact from \cite[Lemma~7]{SST24}.

\bl[Edge speeds]
Let\label{L:edge} $\phi:\{0,1\}^{\Z^d}\to\{0,1\}$ be a non-constant local monotone function and let $\Psi_\phi$ be defined in (\ref{Psiphi}). Then
\be\label{edge}
\Psi^t_\phi(H^\ell_r)=H^\ell_{r-t\eps_\phi(\ell)}\quad(t\geq 0)
\quad\mbox{with}\quad
\eps_\phi(\ell):=\sup_{\Aa\in\Ai(\phi)}\inf_{i\in\Aa}\ell(i),
\ee
where $\Ai(\phi)$ denotes the set of all minimal one-sets of $\phi$.
\el

We call $\eps_\phi(\ell)$ defined in (\ref{edge}) the \emph{edge speed} of $\phi$ in the direction $\ell$. By definition, a \emph{linear polar function} of dimension $\sig\geq 2$ is a linear function $\Ll:\R^d\to\R^\sig$ such that $\sum_{s=1}^\sig\Ll_s(z)=0$ for all $z\in\R^d$. Then \cite[Lemma~12]{Pon13} (see also \cite[Lemma~10]{SST24}) says the following.

\bl[Erosion criterion]
Let\label{L:erode} $\phi:\{0,1\}^{\Z^d}\to\{0,1\}$ be a non-constant monotone function. Then $\phi$ is an eroder if and only if there exists a linear polar function $\Ll$ of dimension $\sig\geq 2$ such that $\sum_{s=1}^\sig\eps_\phi(\Ll_s)>0$.
\el

It is easy to see that the condition $\sum_{s=1}^\sig\eps_\phi(\Ll_s)>0$ implies that $\phi$ is an eroder. We first observe that $\bigvee_{s=1}^\sig H^{\Ll_s}_{r'_s}$ is the all-one configuration when $\sum_{s=1}^\sig r'_s\leq 0$. Indeed, by the defining property of a linear polar function, the latter implies that for each $i\in\Z^d$ one has $\sum_{s=1}^\sig\big(\Ll_s(i)-r'_s\big)\geq 0$ and hence $\Ll_s(i)\geq r'_s$ for at least one $s$. Now if $x$ is a configuration that contains finitely many zeros, then we can find $r_s\in\R$ such that $x\geq H^{\Ll_s}_{r_s}$ $(\lis)$. Using monotonicity, it follows that
\be\label{edgeros}
\Psi_\phi^t(x)\geq\bigvee_{s=1}^\sig H^{\Ll_s}_{r_s-t\eps_\phi(\Ll_s)}\qquad(t\geq 0),
\ee
where the right-hand side is the all-one configuration for all $t$ large enough such that $\sum_{s=1}^\sig r_s\leq t\sum_{s=1}^\sig\eps_\phi(\Ll_s)$.

Toom's Peierls argument was designed to make use of edge speeds. His stability theorem (Theorem~\ref{T:Toom}) basically says that the stability of deterministic monotone cellular automata against small random perturbations is determined by the edge speeds of the unperturbed dynamics. In the presence of intrinsic randomness, the unperturbed dynamics is random and half-space configurations are not mapped into half-space configurations. Nevertheless, we believe that for monotone cellular automata $\Phi^p$ with intrinsic randomness, it should often be possible to define a sort of ``effective edge speeds'' of the unperturbed evolution, and these should determine whether $\Phi^0$ is stable or not. It is not completely clear, however, how such ``effective edge speeds'' should be defined. They are certainly not a simple function of the edge speeds of the individual maps $\phi_1,\ldots,\phi_m$. As a first result in this direction, in \cite{SST24}, we proved the following theorem, that uses a very strong definition of edge speed.

\bt[Theorem~9 of \cite{SST24}]
Let\label{T:thm9} $\Phi^p$ be the monotone cellular automaton defined in (\ref{PhH}). Assume (\ref{Zd}). Then $\Phi^0$ is stable in the sense of (\ref{stable}) if there exists a linear polar function $\Ll$ of dimension $\sig\geq 2$ such that
\be\label{thm9}
\sum_{s=1}^\sig\inf_{1\leq k\leq m}\eps_{\phi_k}(\Ll_s)>0.
\ee
\et

The quantity $\inf_{1\leq k\leq m}\eps_{\phi_k}(\Ll_s)$ is the \emph{worst-case edge speed} in the direction $\Ll_s$. By the argument in (\ref{edgeros}), the condition (\ref{thm9}) implies that under the unperturbed evolution $\Phi^0$, any finite collection of zeros disappears after a finite \emph{deterministic} time. This is a very strong condition that is clearly not satisfied by the cellular automata of Theorems \ref{T:triang} and \ref{T:coop} (except when $r=1$ in the latter case). In our present paper, we demonstrate that the restrictive condition (\ref{thm9}) is not a fundamental drawback of Toom's Peierls argument. Indeed, Theorems \ref{T:triang} and \ref{T:coop} are based on the same Peierls argument that we used in \cite{SST24} to prove Theorem~\ref{T:thm9}. To get rid of (\ref{thm9}), however, one needs better methods to estimate the Peierls sum. These are the main contribution of the present paper.

Although edge speeds are easiest to formulate and understand if space-time has the simple structure $\La=\Z^d\times\Z$ as in (\ref{Zd}), the general idea is more widely applicable. Let $\La$ be a general countable group. On $\La$, we define a \emph{polar function} of \emph{dimension} $\sig\geq 2$ to be a function $L:\La\to\R^\sig$ such that
\be\label{polar2}
\sum_{s=1}^\sig L_s(i)=0\qquad(i\in\La).
\ee
We will say that $L$ is \emph{$\La$-linear} if $L$ is a group homomorphism from $\La$ to $(\R^\sig,+)$, i.e., if
\be\label{Lalin}
L_s(0)=0\quand L_s(ij)=L_s(i)+L_s(j)\qquad(\lis,\ i,j\in\La).
\ee
The following is \cite[Thm.~1]{Too80}, specialised to our set-up.

\bt[Theorem~1 of \cite{Too80}]
Let\label{T:ToomL} $\phh$ be a local monotone map satisfying (\ref{phiDe}) and let $\Phi^0$ be the deterministic monotone cellular automaton that applies the map $\phh$ in each lattice point. Assume that there exists a $\La$-linear polar function $L:\La\to\R^\sig$ of dimension $\sig\geq 2$ and minimal one-sets $A_1,\dots,A_\sig\in\Ai(\phh)$ such that
\be
L_s(j)>0\qquad(j\in A_s,\ \lis).
\ee
Then $\Phi^0$ is stable.
\et

Toom \cite{Too80} has shown that in the special setting of (\ref{Zd}), Theorem~\ref{T:ToomL} implies his celebrated stability theorem (Theorem~\ref{T:Toom}). In the next subsection, we likewise formulate two abstract results that as we will show imply the concrete Theorems \ref{T:NEC}, \ref{T:triang}, and \ref{T:coop}.

\subsection{Two abstract bounds}\label{S:abstract}

In this subsection we work in the general setting of Subsections \ref{S:auto} and \ref{S:stable}. We state two abstract lower bounds on $\ov\rho(p)$ defined in (\ref{ovrho}). These abstract bounds will be proved in Subsection~\ref{S:subpro} and in Subsections \ref{S:some}--\ref{S:coopro} it will be shown that they imply Theorems \ref{T:NEC}, \ref{T:triang}, and \ref{T:coop}. The bounds depend on a large number of parameters that can be freely chosen. We first need some definitions.

A \emph{subprobability distribution} on a finite set $S$ is a function $\pi:S\to\R$ such that $\pi(i)\geq 0$ for all $i\in S$ and $\sum_{i\in S}\pi(i)\leq 1$. Throughout the present subsection, we make the following assumptions:
\begin{enumerate}
\item $L$ is a $\La$-linear polar function of dimension $\sig$ and $\la\in\R$ is a real constant.
\item $\abf=\dis(\abf^\star_{s,k})^{\star\in\{\bullet,\circ\}}_{\lis,\ 1\leq k\leq m}$ is a collection of subprobability distributions on $\La$ such that $\abf^\bullet_{s,k}$ is concentrated on $A_{s,k}$ and $\abf^\circ_{s,k}$ is concentrated on $\De_k:=\bigcup_{s=1}^\sig A_{s,k}$.
\item $\hat\rbf=\dis(\hat\rbf_t)_{t\in\{\circ,1,\ldots,\sig\}}$ is a collection of probability distributions on $\{1,\ldots,m\}$.
\end{enumerate}
We let $\dis(\al^\star_{s,k})^{\star\in\{\bullet,\circ\}}_{\lis,\ 1\leq k\leq m}$ denote the smallest possible constants such that the following inequalities hold
\be\label{alpha}
\ex{-\la L_s(i)}\leq\al^\bullet_{s,k}\abf^\bullet_{s,k}(i)\quad(i\in A_{s,k})\quand\ex{-\la L_s(i)}\leq\al^\circ_{s,k}\abf^\circ_{s,k}(i)\quad(i\in\De_k),
\ee
and similarly, in the special case that $\sig=2$, we let $\dis(\ti\al^\circ_{s,k})_{1\leq s\leq 2,\ 1\leq k\leq m}$ denote the smallest possible constants such that
\be\label{varalpha}
\ex{-\la L_s(i)}\leq\ti\al^\circ_{1,k}\abf^\circ_{1,k}(i)\quad(i\in A_{2,k})\quand\ex{-\la L_s(i)}\leq\ti\al^\circ_{2,k}\abf^\circ_{2,k}(i)\quad(i\in A_{1,k}).
\ee
Also, we let $\dis(\bet_{t,k})_{t\in\{\circ,1,\ldots,\sig\},\ 1\leq k\leq m}$ denote the smallest constants such that
\be\label{beta}
\rbf(k)\leq\bet_{t,k}\hat\rbf_t(k)\qquad\big(t\in\{\circ,1,\ldots,\sig\},\ 1\leq k\leq m\big),
\ee
and we define constants $B_\bullet$ and $B_\circ$ by
\be\label{BB}
B_\star:=\sum_{j\in\La}\sup_{1\leq k\leq m}\abf^\star_{\sig,k}(j)\qquad\big(\star\in\{\bullet,\circ\}\big).
\ee
Here is our first abstract bound.

\bt[Bound based on a polar function of dimension two]
Assume\label{T:cycbd} that $\sig=2$ and $\hat\rbf_\circ=\hat\rbf_2=\rbf$.\\ Let $\dis C_\star:=B_\star\prod_{s=1}^2\ga^\star_s$ $\big(\star\in\{\bullet,\circ\}\big)$ with
\be\label{betcyc}
\ga^\bullet_s:=\sup_{1\leq k\leq m}\bet_{s,k}\al^\bullet_{s,k},\quad\ga^\circ_s:=\sup_{1\leq k\leq m}\ti\al^\circ_{s,k}\qquad(1\leq s\leq 2).
\ee
Let $\hat p>0$, assume that $p_\circ:=1-(1-\hat p)^{-\sig}C_\bullet>0$, and set $\eps:=\hat p(1-\hat p)p_\circ C_\circ^{-1}$. Then
\be\label{cycbd}
1-\ov\rho(p)\leq\frac{p}{\hat p(1-\hat p)}\qquad(p\leq\eps).
\ee
\et

We next state a result that holds for polar functions of any dimension. This case is more complicated and depends on additional parameters that can be freely chosen. We assume that $\dis(\bet^{s'}_{s,k})^{1\leq s'\leq\sig}_{\lis,\ 1\leq k\leq m}$ are positive constants such that
\be\label{bet3}
\bet_{s,k}=\prod_{s'=1}^\sig\bet^{s'}_{s,k}\quad(\lis,\ 1\leq k\leq m)
\quand\bet^{s'}_{s,k}\leq 1\quad(s\neq s'),
\ee
and we set
\be\label{betde}
\bet^s_k:=\sup_{1\leq s'\leq\sig}\bet^s_{s',k}\quad(\lis,\ 1\leq k\leq m)
\quand
\de:=\inf_{1\leq k\leq m}\prod_{s=1}^\sig(1\wedge\bet^s_k).
\ee
(Note the difference in the indices: we use~$\bet^{s'}_{s,k}$ in~\eqref{bet3} and~$\bet^s_{s',k}$ in~\eqref{betde}.) Observe that in the special case that $\hat\rbf_t=\rbf$ for all $t\in\{\circ,1,\ldots,\sig\}$, the theorem simplifies a lot since one can take $\bet^{s'}_{s,k}=1$ for all $s,s',k$ which implies that $\bet^s_k=1$ for all $s,k$ and also $\de=1$.

%Note that below, because of the assumption that $\hat\rbf_\circ=\rbf$, one has $\bet_{\circ,k}=1$ for all $1\leq k\leq m$.

\bt[Bound based on a general polar function]
Assume\label{T:conbd} that $\hat\rbf_\circ=\rbf$ and set
\be\ba{ll}\label{betgast}
\dis\ga^\bullet_s:=\sup_{1\leq k\leq m}\bet^s_k\al^\bullet_{s,k},\quad&\dis\ga^\circ_s:=\sup_{1\leq k\leq m}\al^\circ_{s,k}\qquad(\lis),\\[5pt]
\dis C_\bullet:=B_\bullet\prod_{s=1}^\sig\ga^\bullet_s,\quad&\dis C_\circ:=\de^{-1}B_\circ\prod_{s=1}^\sig\ga^\circ_s.
\ec
Let $\hat p>0$, assume that $p_\circ:=1-(1-\hat p)^{-\sig}C_\bullet>0$, and set $\eps:=\hat p(1-\hat p)p_\circ C_\circ^{-1}$. Then
\be\label{conbd}
1-\ov\rho(p)\leq\frac{p}{\hat p(1-\hat p)}\qquad(p\leq\eps).
\ee
\et

In Subsections \ref{S:some}--\ref{S:coopro} we will show how Theorems \ref{T:cycbd} and \ref{T:conbd} can be applied to obtain Theorems \ref{T:NEC}, \ref{T:triang}, and \ref{T:coop}. These subsections may serve as some sort of manual for how one could apply Theorems \ref{T:cycbd} and \ref{T:conbd} more generally.

\subsection{Discussion and open problems}\label{S:discus}

As we already mentioned in Subsection~\ref{S:edge}, we believe that for monotone cellular automata with intrinsic randomness, it should be possible to define some sort of ``effective edge speeds'' and these should determine stability. As will be explained in Section~\ref{S:randcont}, this is the leading idea behind the results of the present paper and also behind our earlier result \cite[Thm~9]{SST24}. A further result in this direction has earlier been proved by Gray \cite{Gra99}. In the present subsection, we discuss the relation between our work and Gray's and mention some open problems.

We will be interested in the case when $m=2$ with $\phi_1=:\phi$ an eroder and $\phi_2=\phi^{\rm id}$ the \emph{identity map} $\phi^{\rm id}(x):=x(0)$. In Subsection~\ref{S:coop}, we looked at the special case that $\phi=\phi^{\rm coop}$, the cooperative branching map defined as
\be
\phi^{\rm coop}(x):=x(0,0)\vee[x(0,1)\wedge x(1,0)]\qquad\big(x\in\{0,1\}^{\Z^2}\big).
\ee
Inspired by Gray \cite{Gra99}, let us say that a monotone local map $\phi:\{0,1\}^{\Z^d}\to\{0,1\}$ is a \emph{shrinker} if there exists a linear polar function $\Ll$ of dimension $\sig\geq 2$ such that the edge speeds defined in (\ref{edge}) satisfy
\be\label{shrinker}
{\rm(i)}\ \sum_{s=1}^\sig\eps_\phi(\Ll_s)>0,\qquad{\rm(ii)}\ \eps_\phi(\Ll_s)\geq 0\quad(\lis).
\ee
We say that the ``space-time'' map $\phh$ is a shrinker, if the associated ``space'' map $\phi$ is a shrinker. By Lemma~\ref{L:erode}, condition~(i) implies that each shrinker is an eroder. In view of condition~(ii), not every eroder is a shrinker. It is easy to check that the cooperative branching map $\phi^{\rm coop}$ (and thus $\phh^{\rm coop}$) is a shrinker using the linear polar function of dimension two defined as $\Ll_1(z):=z_1+z_2$ and $\Ll_2(z):=-\Ll_1(z)$. The following conjecture is a sweeping generalisation of Theorem~\ref{T:coop}~(i).

\bcon[Shrinkers]
Assume\label{C:shrink} (\ref{Zd}). Let $\Phi^0$ be a monotone cellular automaton that applies a shrinker $\phh$ and the identity map $\phh^{\rm id}$ with probabilities $r$ and $1-r$, respectively. Then $\Phi^0$ is stable for each $r>0$.
\econ

Similarly, we conjecture that a monotone interacting particle system that applies a shrinker with rate one is stable against applications of the zero map $\phh^0$ with small positive rates. This latter conjecture is more or less proved, in fact, in \cite[Thm~18.3.1]{Gra99}, except that Gray uses a somewhat different definition of a shrinker than we do. We suspect the two definitions are equivalent, but unfortunately, due to the complicated nature of Gray's definition, this is far from obvious.

Unfortunately, Conjecture~\ref{C:shrink} cannot be proved with the Peierls argument of the present article. In Subsection~\ref{S:infin} below, we will give an example of a shrinker so that for certain values of $0<r<1$, the Peierls sum is infinite for all $p>0$ small enough. The reasons for this failure are discussed in more detail in Subsection~\ref{S:infin} below. In short, it seems that a drawback of the methods of the present article is that they are \emph{annealed} methods. For monotone cellular automata with intrinsic randomness, there are two sources of randomness: the intrinsic randomness that determines which of the maps $\phh_1,\ldots,\phh_m$ is applied at a given space-time point for the unperturbed evolution, and the extra randomness used to perturb this evolution by adding defective sites. Our Peierls argument is based on estimating the expected number of Toom contours rooted at the origin that are present in $\Phi^p$. Instead of looking at the unconditional expectation of the number of contours, it is possible that better results can be obtained with a \emph{quenched} approach that looks at the conditional expectation given the intrinsic randomness.

\section{Random Toom contours}\label{S:randcont}

\subsection{Toom contours}\label{S:contour}

In this section, we explain our methods and discuss their limitations. We start by recalling Toom's Peierls argument in the reformulation of \cite{SST24}. Before we can state the main theorems, we have no choice but to go through a fairly large number of definitions that we cite almost verbatim from \cite{SST24}.

Let $A$ and $B$ be finite sets. By definition, a \emph{typed directed graph} with \emph{vertex set} $V$, \emph{vertex type set} $A$, and \emph{edge type set} $B$ is a pair $(\Vi,\Ei)$ where $\Vi$ is a subset of $V\times A$ and $\Ei$ is a subset of $V\times V\times B$, such that
\be\label{typedef}
\forall v\in V\ \exists a\in A\mbox{ s.t.\ }(v,a)\in\Vi.
\ee
For each $a\in A$ and $b\in B$, we call
\be\label{typedge}
V_a:=\big\{v:(v,a)\in\Vi\big\}\quand
\vec E_b:=\big\{(v,w):(v,w,b)\in\Ei\big\}
\ee
the set of vertices of type $a$ and the set of directed edges of type $b$, respectively. Note that vertices can have more than one type, i.e., $V_a$ and $V_{a'}$ are not necessarily disjoint for $a\neq a'$, and the same applies to edges. As a consequence, several edges of different types can connect the same two vertices $v,w$, but always at most one of each type. If $(\Vi,\Ei)$ is a typed directed graph, then we let $(V,\vec E)$ denote the directed graph given by
\be
V=\bigcup_{a\in A}V_a\quand\vec E:=\bigcup_{b\in B}\vec E_b,
\ee
where the first equality follows from (\ref{typedef}) and the second equality is a definition. We call $(V,\vec E)$ the \emph{untyped} directed graph associated with $(\Vi,\Ei)$. We also set $E:=\big\{\{v,w\}:(v,w)\in\vec E\big\}$. Then $(V,E)$ is an undirected graph, which we call the undirected graph \emph{associated with} $(V,\vec E)$. We say that a typed directed grap $(\Vi,\Ei)$ or a directed graph $(V,\vec E)$ are \emph{connected} if their associated undirected graph $(V,E)$ is connected. A \emph{rooted} directed graph is a triple $(v_\circ,V,\vec E)$ such that $(V,\vec E)$ is a directed graph and $v_\circ\in V$ is a specially designated vertex, called the \emph{root}. Rooted undirected graphs and rooted typed directed graphs are defined in the same way.

For any directed graph $(V,\vec E)$, we let
\be
\vec E_{\rm in}(v):=\big\{(u,v')\in\vec E:v'=v\big\}
\quand
\vec E_{\rm out}(v):=\big\{(v',w)\in\vec E:v'=v\big\}
\ee
denote the sets of directed edges entering and leaving a given vertex $v\in V$, respectively. Similarly, in a typed directed graph, $\vec E_{b,{\rm in}}(v)$ and $\vec E_{b,{\rm out}}(v)$ denote the sets of incoming or outgoing directed edges of type~$b$ at $v$.

We adopt the following general notation. For any directed graph $(V,\vec E)$, set $\La$, and function $\psi:V\to\La$, we let
\be\label{psiedge}
\psi(V):=\big\{\psi(v):v\in V\big\}\quand
\psi(\vec E):=\big\{\big(\psi(v),\psi(w)\big):(v,w)\in\vec E\big\}
\ee
denote the images of $V$ and $\vec E$ under $\psi$. We can naturally view $\big(\psi(V),\psi(\vec E)\big)$ as a directed graph with set of vertices $\psi(V)$ and set of directed edges $\psi(\vec E)$. We denote this graph by $\psi(V,\vec E):=\big(\psi(V),\psi(\vec E)\big)$. Similarly, if $(\Vi,\Ei)$ is a typed directed graph, then we let $\psi(\Vi,\Ei)$ denote the typed directed graph defined as
\be\ba{r}\label{psiEi}
\dis\psi(\Vi,\Ei):=\big(\psi(\Vi),\psi(\Ei)\big)
\qquad\quad\mbox{with}\quad
\psi(\Vi):=\big\{\big(\psi(v),a\big):(v,a)\in\Vi\big\}\\[5pt]
\dis\quand
\psi(\Ei):=\big\{\big(\psi(v),\psi(w),b\big):(v,w,b)\in\Ei\big\}.
\ec
Also, if $(v_\circ,V,\vec E)$ is a rooted directed graph, then we let $\psi(v_\circ,V,\vec E)$ denote the rooted directed graph $\big(\psi(v_\circ),\psi(V),\psi(\vec E)\big)$, and we use similar notation for rooted typed directed graphs. Two typed directed graphs $(\Vi,\Ei)$ and $(\Wi,\Fi)$ are \emph{isomorphic} if there exists a bijection $\psi:V\to W$ such that $\psi(\Vi,\Ei)=(\Wi,\Fi)$. Similar conventions apply to directed graphs, rooted directed graps, and so on.

\begin{figure}[t]
\begin{center}
\inputtikz{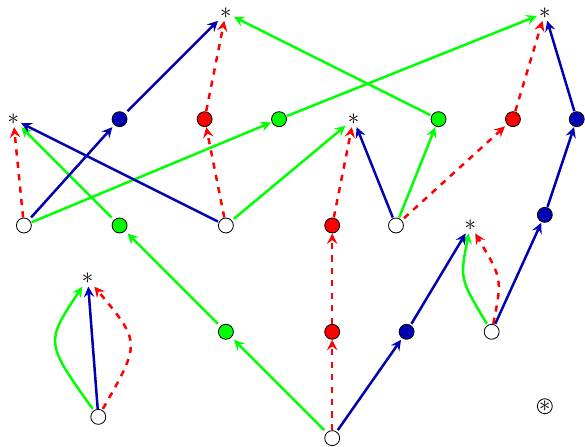}
\caption{Example of a Toom graph with three charges. Sources are indicated with open dots, sinks with asterixes, and internal vertices and edges of the three possible charges with three colours. Note the isolated vertex in the lower right corner, which is a source and a sink at the same time.}
\label{fig:Toomgraph}
\end{center}
\end{figure}

\begin{defi}\label{D:toomgraph}
A \emph{Toom graph} with $\sig\geq 1$ \emph{charges} is a typed directed graph $(\Vi,\Ei)$ with vertex type set $\{\circ,\ast,1,\ldots,\sig\}$ and edge type set $\{1,\ldots,\sig\}$ that satisfies the following conditions:
\begin{enumerate}
\item $|\vec E_{s,{\rm in}}(v)|=0$ $(\lis)$ and $|\vec E_{1,{\rm out}}(v)|=\cdots=|\vec E_{\sig,{\rm out}}(v)|\leq 1$ for all $v\in V_\circ$.
\item $|\vec E_{s,{\rm out}}(v)|=0$ $(\lis)$ and $|\vec E_{1,{\rm in}}(v)|=\cdots=|\vec E_{\sig,{\rm in}}(v)|\leq 1$ for all $v\in V_\ast$.
\item $|\vec E_{s,{\rm in}}(v)|=1=|\vec E_{s,{\rm out}}(v)|$ and $|\vec E_{l,{\rm in}}(v)|=0=|\vec E_{l,{\rm out}}(v)|$ for each $l\neq s$ and $v\in V_s$.
\end{enumerate}
\end{defi}

See Figure~\ref{fig:Toomgraph} for a picture of a Toom graph with three charges.
Toom graphs and the associated Toom contours that will be defined below are designed to make use of the characterisation of eroders in terms of edge speeds from Lemma~\ref{L:erode}. Vertices in $V_\circ,V_\ast$, and $V_s$ are called \emph{sources}, \emph{sinks}, and \emph{internal vertices} with \emph{charge} $s$, respectively. Vertices in $V_\circ\cap V_\ast$ are called \emph{isolated vertices}. With the exception of isolated vertices, the inequalities $\leq 1$ in (i) and (ii) are equalities. Informally, we can imagine that at each source there emerge $\sig$ charges, one of each type, that then travel via internal vertices of the corresponding charge through the graph until they arrive at a sink, in such a way that at each sink there converge precisely $\sig$ charges, one of each type. This informal picture holds even for isolated vertices, if we imagine that in this case, the charges arrive immediately at the sink that is at the same time a source. It is clear from this informal picture that $|V_\circ|=|V_\ast|$, i.e., the number of sources equals the number of sinks. We let $(V,\vec E)$ denote the directed graph associated with $(\Vi,\Ei)$.

Recall that a rooted directed graph is a directed graph with a specially designated vertex, called the root. In the case of Toom graphs, we will always assume that the root is a source.

\begin{defi}\label{D:toomroot}
A \emph{rooted Toom graph} with $\sig\geq 1$ \emph{charges} is a rooted typed directed graph $(v_\circ,\Vi,\Ei)$ such that $(\Vi,\Ei)$ is a Toom graph with $\sig\geq 1$ charges and $v_\circ\in V_\circ$. For any rooted Toom graph $(v_\circ,\Vi,\Ei)$, we write
\be\label{toomroot}
V'_\circ:=V_\circ\beh\{v_\circ\}\quand V'_s:=V_s\cup\{v_\circ\}\quad(\lis).
\ee
\end{defi}

The idea behind (\ref{toomroot}) is that for rooted Toom contours, we view the root more as if it were a collection of internal vertices than as a source. This is reflected in condition~(ii) of the following definition.

\begin{defi}\label{D:embed}
Let $(v_\circ,\Vi,\Ei)$ be a rooted Toom graph and let $\La$ be a countable set. An \emph{embedding} of $(v_\circ,\Vi,\Ei)$ in $\La$ is a map $\psi:V\to\La$ such that:
\begin{enumerate}
\item $\psi(v_1)\neq\psi(v_2)$ for each $v_1\in V_\ast$ and $v_2\in V$ with $v_1\neq v_2$,
\item $\psi(v_1)\neq\psi(v_2)$ for each $v_1,v_2\in V'_s$ with $v_1\neq v_2$ $(\lis)$.
\end{enumerate}
\end{defi}

Condition~(i) says that sinks do not overlap with other vertices and condition~(ii) says that internal vertices do not overlap with other internal vertices of the same charge, where in line with (\ref{toomroot}) we view the root as a collection of internal vertices.

\begin{defi}\label{D:contour}
Let $\La$ be a countable set. A \emph{Toom contour} in $\La$ with $\sig\geq 1$ charges is a quadruple $(v_\circ,\Vi,\Ei,\psi)$, where $(v_\circ,\Vi,\Ei)$ is a rooted connected Toom graph with $\sig$ charges and $\psi$ is an embedding of $(v_\circ,\Vi,\Ei)$ in $\La$. We say that the Toom contour is \emph{rooted} at $i_\circ:=\psi(v_\circ)$. Two Toom contours $(v_\circ,\Vi,\Ei,\psi)$ and $(v'_\circ,\Vi',\Ei',\psi')$ are \emph{isomorphic} if there exists a bijection $\chi:V\to V'$ such that $\chi(v_\circ,\Vi,\Ei)=(v'_\circ,\Vi',\Ei')$ and $\psi(v)=\psi'(\chi(v))$ $(v\in V)$.
\end{defi}

\begin{figure}[htb]
\begin{center}
\inputtikz{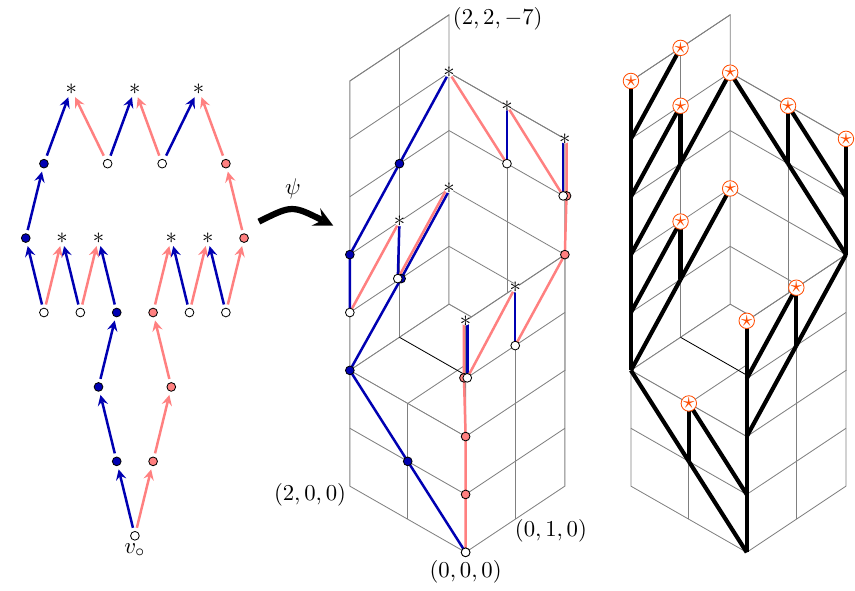}
\caption{The left and middle picture show a Toom contour in $\Z^3$ rooted at $(0,0,0)$, that is present in the typed dependence graph $(\La,\Hi)$ defined in (\ref{Hicoop}). The associated monotone cellular automaton $\Phi^p$ applies the constant map $\phh^0$ in the set $\La_0$ of defective space-time points and the map $\phh^{\rm coop}$ from (\ref{coopid}) in all other points. In the picture on the right defective sites are marked with a star and the black lines help us understand that the maximal trajectory satisfies $\ov X^p(0,0,0)=0$. Sinks of the Toom contour coincide with defective sites. The height function $h(i_1,i_2,i_3):=-i_3$ is plotted upwards in these pictures.}
\label{fig:cycle}
\end{center}
\end{figure}

See Figure~\ref{fig:cycle} for an example of a Toom contour with two charges.

\subsection{Presence of Toom contours}\label{S:presence}

Our next aim is to define when a Toom contour is present in a monotone cellular automaton $\Phq=(\phh_i)_{i\in\La}$. This will require us to make some extra assumptions and equip $\Phq$ with some extra structure. We make use of the general notation for typed directed graphs introduced in (\ref{typedge}).

\begin{defi}\label{D:typdep}
A \emph{typed dependence graph} with $\sig\geq 1$ types of edges is a typed directed graph $(\La,\Hi)$ with vertex type set $\{0,1,\bullet\}$ and edge type set $\{1,\ldots,\sig\}$ such that for $\vec H_s=\{(i, j): (i, j, s)\in\Hi\}$
\begin{enumerate}
\item $\vec H_{s,{\rm out}}(i)=\emptyset$ for all $i\in\La_0\cup\La_1$ and $\lis$,
\item $\vec H_{s,{\rm out}}(i)\neq\emptyset$ for all $i\in\La_\bullet$ and $\lis$,
\end{enumerate}
and its associated untyped directed graph $(\La,\vec H)$ is acyclic. The monotone cellular automaton $\Phq=(\phh_i)_{i\in\La}$ \emph{associated with} the typed dependence graph $(\La,\Hi)$ is defined by
\be\label{phhFi}
\phh_i(x)=\left\{\ba{ll}
\dis\bigvee_{s=1}^\sig\bigwedge_{j:\,(i,j)\in\vec H_s}x(j)\quad&\dis\mbox{if }i\in\La_\bullet,\\[5pt]
\dis r\quad&\dis\mbox{if }i\in\La_r\quad(r=0,1),
\ea\right.
\ee
$\big(i\in\La,\ x\in\{0,1\}^\La\big)$.
\end{defi}

It is useful to see an example. Let $\La:=\Z^3$ and let $(\mu(i))_{i\in\La}$ be i.i.d.\ $\{0,1\}$-valued random variables with $\P[\mu(i)=0]=p$ $(i\in\La)$. Set $\La_0:=\{i\in\La:\mu(i)=0\}$, $\La_1:=\emptyset$, and $\La_\bullet:=\La\beh\La_0$, let $\sig:=2$, and let $\Hi$ be defined by
\bc\label{Hicoop}
\dis\vec H_1&:=&\dis\big\{\big((i_1,i_2,i_3),(i_1+1,i_2,i_3-1)\big),\big((i_1,i_2,i_3),(i_1,i_2+1,i_3-1)\big):i\in\La_\bullet\big\},\\[5pt]
\dis\vec H_2&:=&\dis\big\{\big((i_1,i_2,i_3),(i_1,i_2,i_3-1)\big):i\in\La_\bullet\big\}.
\ec
Then the monotone cellular automaton $\Phq=(\phh_i)_{i\in\La}$ associated with $(\La,\Hi)$ is given by
\be
\phh_i(x)=\left\{\ba{ll}
\dis\phh^{\rm coop}\big((x(i+j))_{j\in\La}\big)\quad&\dis\mbox{if }i\in\La_\bullet,\\[5pt]
\dis\phh^0\big((x(i+j))_{j\in\La}\big)\quad&\dis\mbox{if }i\in\La_0,
\ea\right.
\ee
where $\phh^{\rm coop}$ is defined in (\ref{coopid}). In other words, $\Phq=\Phi^p$ as defined in (\ref{PhH}) for $m=1$ and $\phh_1:=\phh^{\rm coop}$.

It is easy to see that if $(\La,\Hi)$ is a typed dependence graph, $\Phq$ is its associated monotone cellular automaton, and $(\La,\vec H)$ is its associated untyped directed graph, then $(\La,\vec H)$ is the dependence graph of $\Phq$ as defined in Subsection~\ref{S:auto}. In particular, the assumption that $(\La,\vec H)$ is acyclic guarantees that (\ref{phhFi}) defines a cellular automaton. It is clear from (\ref{phhFi}) that $\phh_i$ is monotone for each $i\in\La$ and that $\phh_i$ is one of the constant maps $\phh^r$ $(r=0,1)$ defined in (\ref{phh01}) if and only if $i\in\La_r$ $(r=0,1)$. Elements of $\La_0$, where the constant zero map is applied, are called \emph{defective} sites. Below, we make use of Definition~\ref{D:toomroot}, i.e., we treat the root as if it were a collection of internal vertices.

\begin{defi}\label{D:present}
Let $(\La,\Hi)$ be a typed dependence graph with $\sig\geq 1$ types of edges. We say that a Toom contour $(v_\circ,\Vi,\Ei,\psi)$ with $\sig$ charges is \emph{present} in $(\La,\Hi)$ if:
\begin{enumerate}
\item $\dis\psi(v)\in\La_0$ for all $\dis v\in V_\ast$,
\item $\dis\big(\psi(v),\psi(w)\big)\in\vec H_s$ for all $(v,w)\in\vec E^\bullet_s$ $(\lis$),
\item $\dis\big(\psi(v),\psi(w)\big)\in\vec H$ for all $(v,w)\in\vec E^\circ$,
\end{enumerate}
where for any rooted Toom graph $(v_\circ,\Vi,\Ei)$, we write
\be\ba{l}\label{Ecirc}
\vec E^\bullet:=\bigcup_{s=1}^\sig\vec E^\bullet_s
\quad\mbox{with}\quad
\vec E^\bullet_s:=\big\{(v,w)\in\vec E_s:v\in V'_s\big\}
\quad(\lis),\\[5pt]
\vec E^\circ:=\bigcup_{s=1}^\sig\vec E^\circ_s
\quad\mbox{with}\quad
\vec E^\circ_s:=\big\{(v,w)\in\vec E_s:v\in V'_\circ\big\}
\quad(\lis).
\ec
\end{defi}

Condition (i) says that sinks of the Toom contour correspond to defective sites of the typed dependence graph. Conditions (ii) and (iii) say that directed edges of the Toom graph $(\Vi,\Ei)$ are mapped to directed edges of the typed dependence graph $(\La,\Hi)$, where edges coming out of an internal vertex must be mapped to edges of the corresponding type, and we treat the root as if it were a collection of internal vertices.

We cite the following theorem from \cite[Thm~23]{SST24}. The basic idea goes back to part~3 of the proof of \cite[Thm~1]{Too80}, but the formulation there is quite different from the formulation in \cite{SST24} which we follow here.

\bt[Presence of a Toom contour]
Let\label{T:contour} $(\La,\Hi)$ be a typed dependence graph with $\sig\geq 1$ types of edges, let $\Phq$ be its associated monotone cellular automaton, and let $\ov x$ be its maximal trajectory. If $\ov x(i)=0$ for some $i\in\La$, then a Toom contour $(v_\circ,\Vi,\Ei,\psi)$ rooted at $i$ is present in $(\La,\Hi)$.
\et

Theorem~\ref{T:contour} is demonstrated in Figure~\ref{fig:cycle}. Defective sites are indicated with stars in the picture on the right and the map $\phh^{\rm coop}$ is applied in all non-defective space-time points. This has the consequence that if $\ov X^p(i,t)=0$ in a non-defective space-time point, then $\ov X^p(i,t-1)=0$ and at least one of $\ov X^p(i_1+1,i_2,t-1)$ and $\ov X^p(i_1,i_2+1,t-1)$ must also be zero. This allows us, on the event that $\ov X^p(0,0,0)=0$, to construct a Toom graph rooted at $(0,0,0)$ such that its edges that start at an internal vertex or the root always point from $(i_1,i_2,t)$ to $(i_1+1,i_2,t-1)$ or $(i_1,i_2+1,t-1)$ if the charge is one (blue edges in the picture) and to $(i_1,i_2,t-1)$ if the charge is two (red edges). For edges that start at a source other than the root the rules are less strict and edges do not have to respect their charge. In Figure~\ref{fig:cycle} the height function $h(i_1,i_2,t):=-t$ is plotted upwards.

We note that the converse of Theorem~\ref{T:contour} does not hold, i.e., the presence in $(\La,\Hi)$ of a Toom contour $(v_\circ,\Vi,\Ei,\psi)$ does not imply that $\ov x(i)=0$. This can be seen from Figure~\ref{fig:cycle}. In this example, if there would be no other defective sites apart from the sinks of the Toom contour, then the origin would have the value one.

\subsection{Toom cycles}

As Figure~\ref{fig:cycle} shows, Toom contours with two charges are essentially cycles. It has been shown in \cite[Thm~26]{SST24} that in the special case that $\sig=2$, a stronger version of Theorem~\ref{T:contour} holds that can often be used to obtain better bounds. In the present subsection, we state this result as Theorem~\ref{T:cycle} below. This will later be used in Theorem~\ref{T:cycbd} which in turn is used in Theorem~\ref{T:coop}. We start by defining Toom cycles, which are Toom contours with two charges and certain additional advantageous properties. In what follows, we use the notation
\be
[n]:=\{0,\ldots,n\}\quand\li n\re:=\{1,\ldots,n-1\}.
\ee
Unlike in the previous subsections, where $\La$ could be any countable set, in the present subsection we will assume that $\La$ is a countable group that is equipped with a height function $h$ as defined in Subsection~\ref{S:stable}. These assumptions are not essential but they will simplify the exposition. We refer to \cite[Subsection~2.4]{SST24} for a more general treatment of Toom cycles that does not use these assumptions. Before we can define Toom cycles we need a preparatory concept which for lack of a better name we will call a ``height cycle''.

\begin{defi}
Let\label{D:cyc} $(\La,h)$ be a countable group $\La$ equipped with a height function $h$. A \emph{height cycle} in $(\La,h)$ is a pair $(V',\psi)$ where $V'=[n]$ for some $n\in 2\N$ and\footnote{For ease of notation we write $\psi_v$ for $\psi(v)$ in case of Toom cycles.} $\psi:V'\to\La$ is a function such that $\psi_0=\psi_n$ and setting $\nab{(v-1,v)}\psi:=h(\psi_v)-h(\psi_{v-1})$ $(1\leq v\leq n)$ one has that
\be
\nab{(0,1)}\psi=1,\quad\nab{(n-1,n)}\psi=-1,\quand\nab{(v-1,v)}\psi\in\{-1,1\}\mbox{ for all }1<v<n.
\ee
We call $n$ the \emph{length} of the height cycle. If $\psi$ is a height cycle in $(\La,h)$, then we define
\bc\label{VVV}
\dis V'_\circ&:=&\dis\big\{v\in\li n\re:\nab{(v-1,v)}\psi=-1\mbox{ and }\nab{(v,v+1)}\psi=+1\big\},\\[5pt]
\dis V'_\ast&:=&\dis\big\{v\in\li n\re:\nab{(v-1,v)}\psi=+1\mbox{ and }\nab{(v,v+1)}\psi=-1\big\},\\[5pt]
\dis V'_1&:=&\dis\big\{v\in\li n\re:\nab{(v-1,v)}\psi=+1\mbox{ and }\nab{(v,v+1)}\psi=+1\big\}\cup\{0\},\\[5pt]
\dis V'_2&:=&\dis\big\{v\in\li n\re:\nab{(v-1,v)}\psi=-1\mbox{ and }\nab{(v,v+1)}\psi=-1\big\}\cup\{n\}.
\ec
\end{defi}

Using Definition~\ref{D:cyc}, we can now define Toom cycles as follows. In condition (ii) of the following definition, we equip $[n]$ with the natural total order and we equip $\{1,\circ,2\}$ with the total order $1<\circ<2$.

\begin{defi}
Let\label{D:cycle} $(\La,h)$ be a countable group $\La$ equipped with a height function $h$. A \emph{Toom cycle} in $(\La,h)$ is a height cycle $(V',\psi)$ such that
\begin{enumerate}
\item $\psi_v\neq\psi_w$ for each $v\in V'_\ast$ and $w\in V'$ with $v\neq w$,
\item if $\psi_v=\psi_w$ for some $v\in V'_s$ and $w\in V'_t$ with $s,t\in\{1,\circ,2\}$ and $s\leq t$, then $v\leq w$.
\end{enumerate}
We say that the Toom cycle $(V',\psi)$ is \emph{rooted} at $i:=\psi_0=\psi_n$.
\end{defi}

Toom cycles are really special cases of Toom contours with two charges. (For example, the Toom contour in Figure~\ref{fig:cycle} is actually a Toom cycle.). Indeed, each Toom cycle $(V',\psi)$ of length $n\geq 2$ defines a Toom contour $(v_\circ,\Vi,\Ei,\psi)$ with two charges via the formulas $V:=\{0,\ldots,n-1\}$ and
\be\ba{l}\label{cycont}
\dis v_\circ:=0,\quad V_\circ:=V'_\circ\cup\{0\},\quad V_\ast:=V'_\ast,\quad V_1:=V'_1\beh\{0\},\quad V_2:=V'_2\beh\{n\},\\[5pt]
\vec E_1:=\big\{(v-1,v):\nab{(v-1,v)}\psi=+1\big\},\quand\vec E_2:=\big\{(v,v-1):\nab{(v-1,v)}\psi=-1\big\}.
\ec
For Toom cycles of length zero we set $v_\circ:=0$, $V_\circ=V_\ast:=\{0\}$, and $V_1=V_2=\vec E_1=\vec E_2:=\emptyset$. Condition~(i) of Definition~\ref{D:cycle} clearly implies condition~(i) of Definition~\ref{D:embed} and likewise, condition~(ii) of Definition~\ref{D:cycle} can be seen to imply the corresponding condition of Definition~\ref{D:embed}. This is explained in more detail in \cite[Section~2.4]{SST24}, as well as the fact that Definition~\ref{D:cycle} is stronger than Definition~\ref{D:embed} so that not every Toom contour with two charges can be obtained from a Toom cycle as in (\ref{cycont}). We next define what it means for a Toom contour to be present in a typed dependence graph $(\La,\Hi)$. For simplicity we will assume that $\La$ is a countable group that is equipped with a height function $h$ that is compatible with $\Hi$ in the sense that
\be\label{Hcompat}
h(j)-h(i)=1\mbox{ for all }(i,j)\in\vec H.
\ee
Note that this is equivalent to saying that the monotone cellular automaton $\Phq$ associated with $(\La,\Hi)$ (in the sense of Definition~\ref{D:typdep}) satisfies assumption (\ref{phhDe}) from Subsection~\ref{S:stable}. The definition below uses the definitions of $\vec E_1$ and $\vec E_2$ in (\ref{cycont}).

\begin{defi}
Let\label{D:cycpres} $(\La,\Hi)$ be a typed dependence graph with 2 types of edges as in Definition~\ref{D:typdep}. Assume that $\La$ is a countable group that is equipped with a height function $h$ that is compatible with $\Hi$ in the sense of (\ref{Hcompat}). We say that a Toom cycle $(V',\psi)$ of length $n\geq 2$ is present in $(\La,\Hi,h)$ if:
\begin{enumerate}
\item $\dis\psi_v\in\La_0$ for all $v\in V'_\ast$,
\item $\dis(\psi_v,\psi_w)\in\vec H_s$ for all $(v,w)\in\vec E_s$ with $v\in V'_s$ $(s=1,2)$,
\item $\dis(\psi_v,\psi_w)\in\vec H_{3-s}$ for all $(v,w)\in\vec E_s$ with $v\in V'_\circ$ $(s=1,2)$.
\end{enumerate}
For Toom cycles of length 0, only condition~(i) applies.
\end{defi}

Conditions (i) and (ii) of Definition~\ref{D:cycpres} correspond to conditions (i) and (ii) of Definition~\ref{D:present} but condition (iii) of Definition~\ref{D:cycpres} is stronger than condition (iii) of Definition~\ref{D:present}. The following theorem specialises \cite[Thm~26]{SST24} to our present, somewhat more restrictive setting.

\bt[Presence of a Toom cycle]
Let\label{T:cycle} $(\La,\Hi)$ be a typed dependence graph with $2$ types of edges. Assume that $\La$ is a countable group that is equipped with a height function $h$ that is compatible with $\Hi$ in the sense of (\ref{Hcompat}). Let $\Phq$ be monotone cellular automaton associated with $(\La,\Hi)$ and let $\ov x$ be its maximal trajectory. If $\ov x(i)=0$ for some $i\in\La$, then a Toom cycle rooted at $i$ is present in $(\La,\Hi,h)$.
\et

\subsection{The Peierls argument}\label{S:peierlsarg}

In formula (\ref{Peierls}) we have already mentioned the fact that for a monotone cellular automaton $\Phi^p$ of the form (\ref{PhH}), the probability that the maximal trajectory $\ov X^p$ has a zero at the origin can be estimated from above by a Peierls sum. We now give a precise definition of this Peierls sum. We work in the set-up of Subsections \ref{S:auto} and \ref{S:stable}. In particular, $\Phi^p$ is the random monotone cellular automaton defined in (\ref{PhH}). Recalling that $\Ai(\phh)$ denotes the set of all minimal one-sets of $\phh$, we choose sets
\be\label{Askchoice}
A_{s,k}\in\Ai(\phh_k)\qquad(\lis,\ 1\leq k\leq m),
\ee
and we use these to define a typed dependence graph $(\La,\Hi^p)$ by
\be\ba{c}\label{LaHi}
\dis\La_0:=\big\{i\in\La:\mu(i)=0\big\},\quad\La_1:=\emptyset,\quad\La_\bullet:=\La\beh\La_0,\\[5pt]
\dis\vec H^p_s:=\big\{(i,j):i\in\La_\bullet,\ j\in A_{s,\mu(i)}\big\},
\ec
and we let $\ov X^p$ denote the maximal trajectory of $\Phi^p$. The following theorem is a precise formulation of (\ref{Peierls}). We recall that Toom contours and cycles, as well as isomorphisms of the former, are defined in Definitions \ref{D:contour} and \ref{D:cycle}. In formulas (\ref{Peicont}) and (\ref{Peicyc}) below, we use somewhat different definitions of what it means to be present, as described by Definitions \ref{D:present} and \ref{D:cycpres}.

\bt[Peierls bounds]
Let\label{T:Peierls} $\Ti_0$ denote the set of Toom contours rooted at $0$, up to isomorphism. Then
\be\label{Peicont}
\P\big[\ov X^p(0)=0\big]\leq\sum_{T\in\Ti_0}\P\big[T\mbox{ is present in }(\La,\Hi^p)\big].
\ee
In the special case that $\sig=2$, let $\bar\Ti_0$ denote the set of Toom cycles rooted at $0$. Then
\be\label{Peicyc}
\P\big[\ov X^p(0)=0\big]\leq\sum_{T\in\bar\Ti_0}\P\big[T\mbox{ is present in }(\La,\Hi^p)\big].
\ee
\et

\bpro
Let $\Psi^p$ denote the monotone cellular automaton associated with $(\La,\Hi^p)$, in the sense of Definition~\ref{D:typdep}, and let $\ov Y^p$ denote its maximal trajectory. In the proof of \cite[Thm~27]{SST24} it is shown that $\ov X^p\leq\ov Y^p$ (pointwise) a.s. In view of this, formulas (\ref{Peicont}) and (\ref{Peicyc}) follow immediately from Theorems \ref{T:contour} and \ref{T:cycle} and the fact that the probability that a contour is present can be estimated from above by the expected number of contours that are present. For the details, we refer to  the proof of \cite[Thm~27]{SST24}. The latter uses a slightly different concept than isomorphism of Toom contours, called ``equivalence'' of contours, but the proof remains the same.
\epro

Up to this point, in the present section, we have been citing results from \cite{SST24}. The novelty of the present paper, that allows us to prove the abstract theorems from Subsection~\ref{S:abstract} and the concrete new results from Subsections \ref{S:NEC}--\ref{S:coop}, is a novel way of estimating the Peierls sums on the right-hand sides of (\ref{Peicont}) and (\ref{Peicyc}).

\subsection{Decorated contours}

We now start preparing for our upper bounds on the Peierls sums in (\ref{Peicont}) and (\ref{Peicyc}). Recall that the random typed dependence graph $(\La,\Hi^p)$ in (\ref{LaHi}) is defined in terms of the i.i.d.\ random variables $(\mu(i))_{i\in\La}$ from (\ref{Phip}) that take values in $\{0,\ldots,m\}$. By condition~(i) of Definition~\ref{D:present}, if a Toom contour $T$ is present in $(\La,\Hi^p)$, then its sinks correspond to defective sites, i.e., sites $i\in\La$ for which $\mu(i)=0$. Therefore, if for a given Toom contour $T=(v_\circ,\Vi,\Ei,\psi)$, we want to calculate the probability that $T$ is present in $(\La,\Hi^p)$, then we pick up a factor $p$ for each sink $i\in\psi(V_\ast)$. For the remaining sites on the contour, i.e., sites $i\in\psi(V)\beh\psi(V_\ast)$, we must have $\mu(i)\in\{1,\ldots,m\}$ and the value of $\mu(i)$ must be compatible with $T$ in the sense of conditions (ii) and (iii) of Definition~\ref{D:present} and the definition of $\Hi$ in terms of $(\mu(i))_{i\in\La}$ in (\ref{LaHi}). To calculate the probability of this happening, we must sum over all possible ways we can assign values $\mu(i)\in\{1,\ldots,m\}$ to sites $i\in\psi(V)\beh\psi(V_\ast)$. This idea is formalised in the following definition.

Recall that $(A_{s,k})_{\lis,\ 1\leq k\leq m}$ are the sets from (\ref{Askchoice}). In line with assumption~(ii) of Subsection~\ref{S:abstract}, we set
\be
\De_k:=\bigcup_{s=1}^\sig A_{s,k}\qquad(1\leq k\leq m).
\ee
Below, we implicitly use that each directed edge in $\psi(\vec E)$ (defined as in (\ref{psiedge})) can be written in the form $(i,ij)$ for some unique $i,j\in\La$.

\begin{defi}
A\label{D:decor} \emph{decorated Toom contour} is a quintuple $D:=(v_\circ,\Vi,\Ei,\psi,\kappa)$ where the quadruple $(v_\circ,\Vi,\Ei,\psi)$ is a Toom contour and $\kappa:\psi(V)\to\{0,\ldots,m\}$ is a function such that
\begin{enumerate}
\item if $i\in\psi(V_\ast)$, then $\kappa(i)=0$,
\item if $(i,ij)\in\psi(\vec E^\bullet_s)$, then $j\in A_{s,\kappa(i)}$ $(\lis)$,
\item if $(i,ij)\in\psi(\vec E^\circ)$, then $j\in\De_{\kappa(i)}$.
\end{enumerate}
We set
\be\label{nast}
n_\ast(D):=|V_\ast|\quand n_\dia(D):=\big|\psi(V)\beh\psi(V_\ast)\big|.
\ee
We call $\kappa:\psi(V)\to\{0,\ldots,m\}$ the \emph{decoration function} and we call $\kappa(i)$ the \emph{decoration} of a vertex $i\in\psi(V)$. The \emph{presence-probability} of the decorated Toom contour $D=(v_\circ,\Vi,\Ei,\psi,\kappa)$ is the quantity
\be\label{prespro}
\pi_{p,\rbf}(D):=p^{n_\ast(D)}\,(1-p)^{n_\dia(D)}\prod_{i\in\psi(V)\beh\psi(V_\ast)}\rbf\big(\kappa(i)\big).
\ee
Two decorated Toom contours are isomorphic if the associated Toom contours are isomorphic and the decoration functions are the same. We let $\Di_0$ denote the set of decorated Toom contours rooted at $0\in\La$ (up to isomorphism).

Decorated Toom cycles are defined analogously, except that condition~(iii) is replaced by
\begin{itemize}
\item[{\rm(iii)'}] If $(i,ij)\in\psi(\vec E^\circ_s)$, then $j\in A_{3-s,\kappa(i)}$ $(s=1,2)$.
\end{itemize}
We let $\bar\Di_0$ denote the set of decorated Toom cycles rooted at $0\in\La$.
\end{defi}

We note that since (iii)' is stronger than (iii), in the special case that $\sig=2$, we can view $\bar\Di_0$ as a subset of $\Di_0$. Elements of $\Di_0$ are defined only up to isomorphism in the sense of Definition~\ref{D:contour}, i.e., up to a reparametrisation of the set $V$. For Toom cycles of length $n$, there is a standard way to enumerate the elements of $V$ from $0$ to $n-1$ (see (\ref{cycont})), so there is no need to discuss isomorphism. The first step in the proof of Theorems \ref{T:cycbd} and \ref{T:conbd} is the following simple observation.

\bl[Peierls sum]
The\label{L:Pei} Peierls sum on the right-hand side of (\ref{Peicont}) is given by
\be\label{Pei}
\sum_{T\in\Ti_0}\P\big[T\mbox{ is present in }(\La,\Hi^p)\big]=\sum_{D\in\Di_0}\pi_{p,\rbf}(D),
\ee
and the Peierls sum on the right-hand side of (\ref{Peicyc}) is given by
\be\label{Peic}
\sum_{T\in\bar\Ti_0}\P\big[T\mbox{ is present in }(\La,\Hi^p)\big]=\sum_{D\in\bar\Di_0}\pi_{p,\rbf}(D).
\ee
\el

\bpro
Let $T=(v_\circ,\Vi,\Ei,\psi)$ be a Toom contour rooted at $0$ and let $\mu\big|_T$ denote the restriction of the function $i\mapsto\mu(i)$ to $\psi(V)$. Comparing Definitions \ref{D:present} and \ref{D:decor}, we see that $T$ is present in $(\La,\Hi^p)$ if and only if $(T,\mu\big|_T)$ is a decorated Toom contour. In view of this, we can rewrite the right-hand side of (\ref{Peicont}) as
\be
\sum_{T\in\Ti_0}\P\big[(T,\mu\big|_T)\mbox{  is a decorated Toom contour}\big]
=\sum_{T\in\Ti_0}\sum_{\kappa\in\Ki_T}\P\big[\mu\big|_T=\kappa\big],
\ee
where $\Ki_T$ denotes the space of functions $\kappa:\psi(V)\to\{0,\ldots,m\}$. For any $\kappa\in\Ki_T$, the presence-probability of the decorated Toom contour $(T,\kappa)$ is precisely equal to
\be
\pi_{p,\rbf}(T,\kappa)=\P\big[\mu\big|_T=\kappa\big].
\ee
Inserting this into our previous formula, we arrive at (\ref{Pei}). The proof for Toom cycles is the same, with Definition~\ref{D:cycpres} replacing Definition~\ref{D:present} and part~(iii) of Definition~\ref{D:decor} replaced by (iii)'.
\epro

\subsection{Random contours}\label{S:rand}

We now come to the central idea of the whole paper. In order to estimate the Peierls sums in (\ref{Pei}) and (\ref{Peic}), we are going to rewrite these sums as the expectation of a function with respect to a subprobability measure on the space of all decorated Toom contours. We will then use the polar function $L$ and real constant $\la$ from assumption~(i) of Subsection~\ref{S:abstract} to show that this function can be bounded uniformly from above, leading to Theorems \ref{T:cycbd} and \ref{T:conbd}. We work in the set-up of Subsections~\ref{S:auto} and \ref{S:stable}. In addition, assumptions (i)--(iii) of Subsection~\ref{S:abstract} will be in force throughout.

\bl[The number of internal edges]
For\label{L:ncirc} each decorated Toom contour $D=(v_\circ,\Vi,\Ei,\psi,\kappa)$, there exist integers $n_\circ(D)$ and $n_\bullet(D)$ such that
\be\label{ncirc}
n_\circ(D)=|\vec E^\circ_s|
\quand
n_\bullet(D)=|\vec E^\bullet_s|\qquad(\lis).
\ee
The constants $n_\ast(D)$ and $n_\dia(D)$ from (\ref{nast}) satisfy
\be\label{dici}
n_\ast(D)=n_\circ(D)+1\quand n_\dia(D)\leq n_\circ(D)+\sig n_\bullet(D)+1.
\ee
\el

\bpro
Since $n_\ast(D)$ is the number of sinks which equals the number of sources and $|\vec E^\circ_s|$ is the number of sources other than the root it is clear that $n_\circ(D)=|\vec E^\circ_s|=n_\ast(D)-1$. We observe that $|\vec E_s|=|\vec E^\circ_s|+|\vec E^\bullet_s|$ is the number of edges of charge $\lis$. Formula (\ref{ncirc}) now follows from (\ref{Askchoice}) and (\ref{phiDe}) which imply that $A_{s,k}\sub\De(\phh_k)\sub\{i\in\La:h(i)=1\}$ for all $s,k$, and as a result
\be\label{sames}
|\vec E_1|=\cdots=|\vec E_\sig|=\sum_{v\in V_\ast}h\big(\psi(v)\big)-\sum_{v\in V_\circ}h\big(\psi(v)\big).
\ee
The inequality in (\ref{dici}) follows from the fact that $V$ contains precisely $n_\ast(D)$ sources and $n_\bullet(D)$ internal vertices of each charge, and the image of the set of all these vertices under $\psi$ has cardinality $n_\dia(D)$. 
\epro

We recall from (\ref{cycont}) that if $(v_\circ,\Vi,\Ei,\psi)$ is the Toom contour associated with a Toom cycle of length $n$, then $V=\{0,\ldots,n-1\}$ where $0$ is the root. In Subsection~\ref{S:incomp} below, we will define a total order on the set of vertices $V$ of a general Toom contour such that $v_0$ is the minimal element of $V$ in this order and in the special case that $(v_\circ,\Vi,\Ei,\psi)$ is associated with a Toom cycle of length $n$ this coincides with the natural total order on $\{0,\ldots,n-1\}$. We call this the \emph{natural order} on $V$. For the moment, the exact definition of this natural order does not matter, but the interested reader can look it up in formula (\ref{natord}) below.

For each Toom graph $(\Vi,\Ei)$, we define $\tau':V\to\{\circ,\ast,1,\ldots,\sig\}$ by the requirements that (recall \eqref{toomroot})
\be
\tau'(v_\circ):=1\quand v\in V'_{\tau'(v)}\quad(v\neq v_\circ)
\ee
i.e., $\tau'(v)=\circ$ if $v$ is a source other than the root, $\tau'(v)=\ast$ if $v$ is a sink, $\tau'(v)=s$ if $v$ is an internal vertex of charge $s$, and the root has $\tau'(v_\circ):=1$. For a Toom contour $(v_\circ,\Vi,\Ei,\psi)$, for each $i\in\psi(V)$, we let
\be\label{taudef}
\tau(i):=\tau'\big(f(i)\big)\quad\mbox{with}\quad f(i):=\min\big\{v\in V:\psi(v)=i\big\}
\ee
denote the type of the first vertex, with respect to the natural order on $V$, that is embedded at $i$. Note that the root gets type $1$. In the special case that $(v_\circ,\Vi,\Ei,\psi)$ is associated with a Toom cycle, we note that because of condition~(ii) of Definition~\ref{D:cycle},
\be\label{cycord}
\tau(i)=\left\{\ba{ll}
\dis\ast\quad&\mbox{if }i\in\psi(V_\ast),\\[5pt]
\dis\min\big\{\tau'(v):\psi(v)=i\big\}\quad\mbox{w.r.t.\ the order }1<\circ<2\quad&\mbox{if }i\in\psi(V)\beh\psi(V_\ast).
\ea\right.
\ee

Using the notation introduced above, we are able to formulate the following theorem, that describes a class of subprobability measures on the space $\Di_0$ of decorated Toom contours rooted at~$0$. In formula (\ref{DecD}), the subscript $(i,ij)\in\psi(\vec E^\star_s)$ means that we take the product over all $j\in\La$ such that $(i,ij)\in\psi(\vec E^\star_s)$, with the convention that the product over the empty set equals one.

\bt[Random decorated contour]
Let\label{T:randcont} $p_\circ,\hat p\in[0,1]$ and let $\abf$ and $\hat\rbf$ be as in assumptions (ii) and (iii) from Subsection~\ref{S:abstract}. Let $B_\bullet$ and $B_\circ$ be defined as in (\ref{BB}). Then there exists a subprobability measure $\nu_{\hat p,p_\circ,\hat\rbf,\abf}$ on $\Di_0$ that is given by
\bc\label{DecD}
\dis\nu_{\hat p,p_\circ,\hat\rbf,\abf}(D)
&:=&\dis\hat p^{n_\ast(D)}\,(1-\hat p)^{n_\dia(D)}\prod_{i\in\psi(V)\beh\psi(V_\ast)}\hat\rbf_{\tau(i)}\big(\kappa(i)\big)\\[5pt]
&&\dis\times\prod_{\star\in\{\circ,\bullet\}}\prod_{s=1}^\sig\;\prod_{(i,ij)\in\psi(\vec E^\star_s)}\abf^\star_{s,\kappa(i)}(j)\\[25pt]
&&\dis\times B_\circ^{-n_\circ(D)}B_\bullet^{-n_\bullet(D)}p_\circ^{n_\circ(D)}(1-p_\circ)^{n_\bullet(D)}
\ec
for each $D\in\Di_0$ with $D=(v_\circ,\Vi,\Ei,\psi,\kappa)$.
\et

\noi
\textbf{Remark} When the number of charges is two, we can restrict the subprobability measure $\nu_{\hat p,p_\circ,\hat\rbf,\abf}$ to the set $\bar\Di_0$ of Toom contours that are associated with a Toom cycle. This restricted measure is, of course, also a subprobability measure. In this setting, because of condition~(iii)' of Definition~\ref{D:decor}, it is natural to choose the subprobability distributions $\abf^\circ_{s,k}$ so that $\abf^\circ_{1,k}$ is concentrated on $A_{2,k}$ and $\abf^\circ_{2,k}$ is concentrated on $A_{1,k}$.\med

The proof of Theorem~\ref{T:randcont} can crudely be described as follows. We construct a random decorated contour step by step, starting from the root, by adding charged edges one by one. Edges of charges $1,\ldots,\sig-1$ are added in the ``upward'' fashion and edges of charge $\sig$ in the ``downward'' fashion . Whenever we add an edge, we have to choose a direction for it to go, which gives the factor in the second line on the right-hand side of (\ref{DecD}). The order in which we add new vertices corresponds to the natural order on $V$. Whenever a new vertex is embedded in a space-time point that we have not used yet, we have to choose a decoration, which gives the factor in the first line on the right-hand side of (\ref{DecD}). We allow the probability of the chosen decoration to depend on the type of the vertex that is first embedded at this point. For edges of charge $\sig$, we choose with probabilities $p_\circ$ and $1-p_\circ$ wether to add a source or an internal vertex at the bottom of the edge. The factors $1/B_\circ$ and $1/B_\bullet$ find their origin in the fact that at the moment the edge is added, we do not yet know the decoration at the bottom of the edge, which forces us to use a coupling between the subprobability distributions $\abf^\star_{\sig,k}$ for different values of~$k$. Filling in the details is quite lengthy, so we postpone the proof of Theorem~\ref{T:randcont} till Section~\ref{S:nuproof}.

\subsection{Factoring in the polar function}

Toom contours were designed to make use of the characterisation of eroders in terms of edge speeds from Lemma~\ref{L:erode}. At the end of Subsection~\ref{S:presence}, we have already discussed how directed edges starting at an internal vertex or the root must always point in the right direction as prescribed by their charge. In the present subsection, for the first time in Section~\ref{S:randcont}, polar functions will enter the picture. Even though the two are closely related, it is important to distinguish polar functions $\Ll$ that are defined on space only as in Lemma~\ref{L:erode} from space-time polar functions $L$ from Subsection~\ref{S:abstract}. In Subsection~\ref{S:idea} below, we will discuss how to ``lift'' a spatial polar function $\Ll$ to a space-time polar function $L$.

In the present subsection, we rewrite the Peierls sums in (\ref{Pei}) and (\ref{Peic}) as expectations of a function with respect to the subprobability measure defined in Theorem~\ref{T:randcont}. We will use the polar function $L$ and real constant $\la$ from assumption~(i) of Subsection~\ref{S:abstract} to write the function to be integrated in a form that will later allow us to estimate it. We work in the set-up of Subsections~\ref{S:auto} and \ref{S:stable}. In addition, assumptions (i)--(iii) of Subsection~\ref{S:abstract} will be in force throughout. We start with a preparatory lemma.
\bl[Factoring in a polar function]
Let\label{L:polfac} $L:\La\to\R^\sig$ be a $\La$-linear polar function of dimension $\sig$. Then for each Toom contour $(v_\circ,\Vi,\Ei,\psi)$, one has
\be\label{polfac}
\sum_{s=1}^\sig\sum_{(i,ij)\in\psi(\vec E_s)}L_s(j)=0.
\ee
\el

\bpro
It is fairly easy to see that for any polar function $L$ of dimension $\sig$, 
\be\label{zerosum}
\sum_{s=1}^\sig\sum_{(v,w)\in\vec E_s}\big(L_s(\psi(w))-L_s(\psi(v))\big)=0.
\ee
This is \cite[Lemma~42]{SST24} so for the (short) proof we refer to that paper. We claim that we can rewrite (\ref{zerosum}) as
\be\label{zerosum2}
\sum_{s=1}^\sig\sum_{(i,ij)\in\psi(\vec E_s)}\big(L_s(ij)-L_s(i)\big)=0.
\ee
This is not completely obvious, since we need to show that each term in (\ref{zerosum2}) corresponds to precisely one term in (\ref{zerosum}), i.e., it never happens that two edges in the Toom graph $(\Vi,\Ei)$ of the same charge $s$ are mapped onto the same edge $(i,ij)\in\psi(\vec E_s)$. This follows, however, from the fact already observed in \cite[Lemma~19]{SST24} that as a consequence of Definition~\ref{D:embed}, at each space-time point $i\in\La$ there can be at most one incoming edge of each charge.

To complete the argument, we note that if $L$ is $\La$-linear as defined in (\ref{Lalin}), then $0=L_s(0)=L_s(i^{-1}i)=L_s(i^{-1})+L_s(i)$ for each $i\in\La$ and $\lis$, and hence $L_s(j)=L_s(i^{-1}ij)=L_s(ij)-L_s(i)$ for each $i,j$, and $s$.
\epro

We now rewrite the Peierls sums in (\ref{Pei}) and (\ref{Peic}) as expectations of a function with respect to the subprobability measure defined in Theorem~\ref{T:randcont}. For simplicity, we will at this moment only do this in the special case that $\hat p:=p$ and $\hat\rbf_t:=\rbf$ for all $t\in\{\circ,1,\ldots,\sig\}$. We let $\nu_{p,p_\circ,\rbf,\abf}$ denote the subprobability measure from Theorem~\ref{T:randcont} with this special choice of $\hat p$ and $\hat\rbf$.

\bl[Peierls sum as expectation]
The\label{L:Peiex} Peierls sum in (\ref{Pei}) can be rewritten as
\bc\label{Peiex}
\dis\sum_{D\in\Di_0}\pi_{p,\rbf}(D)&=&\dis\sum_{D\in\Di_0}\nu_{p,p_\circ,\rbf,\abf}(D)\;B_\circ^{n_\circ(D)}p_\circ^{-n_\circ(D)}\Big(\prod_{s=1}^\sig\;\prod_{(i,ij)\in\psi(\vec E^\circ_s)}\ex{-\la L_s(j)}/\abf^\circ_{s,\kappa(i)}(j)\Big)\\[10pt]
&&\dis\quad\times B_\bullet^{n_\bullet(D)}(1-p_\circ)^{-n_\bullet(D)}\Big(\prod_{s=1}^\sig\;\prod_{(i,ij)\in\psi(\vec E^\bullet_s)}\ex{-\la L_s(j)}/\abf^\bullet_{s,\kappa(i)}(j)\Big),
\ec
where $p_\circ\in[0,1]$, and restricting this sum to $\bar\Di_0$ gives the Peierls sum in (\ref{Peic}).
\el

\bpro
By grace of Lemma~\ref{L:polfac}, given a $\La$-linear polar function $L$ and constant $\la\in\R$, we can rewrite formula (\ref{prespro}) for the presence-probability of a decorated Toom contour $D=(v_\circ,\Vi,\Ei,\psi,\kappa)$ as
\be\label{PiL}
\pi_{p,\rbf}(D)=p^{n_\ast(D)}\,(1-p)^{n_\dia(D)}\!\!\prod_{i\in\psi(V)\beh\psi(V_\ast)}\!\rbf\big(\kappa(i)\big)\;\prod_{s=1}^\sig\prod_{(i,ij)\in\psi(\vec E_s)}\ex{-\la L_s(j)},
\ee
where all we have done compared to (\ref{prespro}) is that we have added a factor that is identically one on the space $\Di_0$ of decorated Toom contours by Lemma~\ref{L:polfac}. Comparing this with (\ref{DecD}), using the fact that $\hat p:=p$ and $\hat\rbf_t:=\rbf$ for all $t\in\{\circ,1,\ldots,\sig\}$, we arrive at (\ref{Peiex}).
\epro

Lemma~\ref{L:Peiex} naturally leads to the following upper bound on the Peierls sum.

\bp[Upper bound in terms of random contours]
Let\label{P:randbd} $\al^\star_{s,k}$ and $B_\star$ be defined as in (\ref{alpha}) and (\ref{BB}), and set
\be\label{stardef}
C_\star:=B_\star\prod_{s=1}^\sig\sup_{1\leq k\leq m}\al^\star_{s,k}\qquad\big(\star\in\{\circ,\bullet\}\big).
\ee
Then for any $p_\circ\in[0,1]$ the Peierls sum in (\ref{Pei}) can be estimated from above as
\be\label{albd}
\sum_{D\in\Di_0}\pi_{p,\rbf}(D)\leq\sum_{D\in\Di_0}\nu_{p,p_\circ,\rbf,\abf}(D)\big(C_\circ p_\circ^{-1}\big)^{n_\circ(D)}\big(C_\bullet(1-p_\circ)^{-1}\big)^{n_\bullet(D)}.
\ee
\ep

\bpro
This follows from Lemma~\ref{L:Peiex} by estimating $\ex{-\la L_s(j)}/\abf^\bullet_{s,\kappa(i)}(j)\leq\al^\bullet_{s,k}$ $(j\in A_{s,k})$ and $\ex{-\la L_s(j)}/\abf^\circ_{s,\kappa(i)}(j)\leq\al^\circ_{s,k}$ $(j\in\De_k)$ and using the definitions of $n_\bullet(D)$ and $n_\circ(D)$ in (\ref{ncirc}).
\epro

\subsection{Proof idea of the main results}\label{S:idea}

In this subsection, we explain the ideas behind our main results, which are Theorems \ref{T:cycbd} and \ref{T:conbd}. Before we do this, it is useful to take one step back and look at some of the simplest known Peierls arguments, that are used to prove long-range order for two-dimensional percolation or for the two-dimensional Ising model. (See \cite[Thm (1.10)]{Gri99} for percolation and \cite[Lemma~3.37]{FV18} for the Ising model.) In these arguments, the Peierls sum can be bounded by an expression of the form $\sum_\ga p^{l(\ga)}$ where the sum runs over all contours, $p$ is a small parameter ($p$ is the percolation parameter in the case of percolation and $e^{-2\bet}$ for the Ising model), and $l(\ga)$ is the length of a contour. Rewriting this as $\sum_{n=1}^\infty N_np^n$ where $N_n$ is the number of contours with length $n$, we see that the value of this sum is determined by an entropy-energy balance. For the Ising model or percolation, $N_n$ grows exponentially in $n$ but as long as $p$ is small enough this is beaten by the exponential factor $p^n$.

For Toom's Peierls argument, the situation is similar with one important complication. We can naturally define the length of a decorated Toom contour $D$ as its total number of edges, which is $\sig(n_\circ(D)+n_\bullet(D))$. As shown in \cite[Lemma~40]{SST24}, one can then show that the number $N_n$ of decorated Toom contours with length $n$ grows exponentially in $n$. However, the presence-probability $\pi_{p,\rbf}(D)$ of a decorated Toom contour $D$ only contains a factor $p^{n_\ast(D)}$, where $n_\ast(D)$ is the number of sinks. In other words, there is an entropic factor for each edge, but we only pay an energy cost for each sink of the contour.

To resolve this, we need a way to redistribute the entropy factors associated with the edges to the sinks, where they can be beaten by the energy costs we pay for the sinks. This is where the edge speeds come in. At the end of Subsection~\ref{S:presence}, we already discussed how in the example of Figure~\ref{fig:cycle}, directed edges coming out of an internal vertex or the root always have to point in the ``right'' direction that corresponds to their charge, while the other directed edges, that come out of a source other than the root, are allowed to point in the ``wrong'' direction. For a suitable choice of the polar function $L$ and sets $A_{s,k}$ from assumption~(i) of Subsection~\ref{S:abstract} and from (\ref{Askchoice}), we can make sure that directed edges of charge $s$ coming out of an internal vertex or the root typically point in a direction in which the function $L_s$ increases. Since in the end, one edge of each charge arrives at each sink, this has to be compensated by edges coming out of a source other than the root, that are allowed to point in the ``wrong'' direction. Since the number of sources equals the numer of the sinks, this allows us to redistribute the entropy factors associated with the edges coming out of internal vertices, for which we would otherwise not pay any cost, to the sinks, for which we pay an energy cost.

To demonstrate this, it is instructive to first look at Theorem~\ref{T:thm9}, which has been proved in \cite[Thm~9]{SST24} by different means. We will give an alternative proof here based on Proposition~\ref{P:randbd}. Note that this result is in the setting of~\eqref{Zd}, that is $\La=\Z^d\times\Z$.

We start with a preparatory lemma. Generalising (\ref{edge}), for each $v\in\R^d$, we define a \emph{compensated edge speed} by
\be\label{compedge}
\eps^v_\phi(\ell):=\sup_{\Aa\in\Ai(\phi)}\inf_{i\in\Aa}\ell(i-v).
\ee
Roughly speaking, for fixed $\phi$ and $v$, these are the edge speeds of the deterministic cellular automaton whose evolution has the following description: in each step, we first apply the map $\phi$ in each lattice point, and then shift the whole lattice over a distance $v$. The following lemma, the proof of which we postpone till Subsection~\ref{S:compspeed}, says that under the assumptions of Theorem~\ref{T:thm9}, by subtracting a ``shift'' $v$, we can make all edge speeds strictly positive. This is closely related to an observation of Gray, who proved that each eroder is ``shift-equivalent'' to a shrinker \cite[Thm~18.2.1]{Gra99}.

\bl[Compensated speeds]
Under\label{L:compspeed} the assumptions of Theorem~\ref{T:thm9} there exists a linear polar function $\Ll':\R^d\to\R^{\sig'}$ of dimension $\sigma'\geq 2$ and a $v \in \R^d$ such that $\eps_{\phi_k}^v(\Ll'_s)>0$ for all $1\leq s\leq\sig'$ and $1\leq k\leq m$.
\el

\bpro[of Theorem~\ref{T:thm9}]
The definition of $\ov\rho(p)$ in (\ref{ovrho}), Theorem~\ref{T:Peierls}, and Lemma~\ref{L:Pei} tell us that
\be\label{onemin}
1-\ov\rho(p)\leq\sum_{D\in\Di_0}\pi_{p,\rbf}(D).
\ee
To estimate the Peierls sum on the right-hand side from above, we will use the estimate (\ref{albd}) from Proposition~\ref{P:randbd} for a clever choice of the polar function $L$, constant $\la$, sets $A_{s,k}$ and subprobability measures $\abf^\star_{s,k}$ from assumptions (i) and (ii) of Subsection~\ref{S:abstract} and from (\ref{Askchoice}). We first explain how to choose $L$ and the sets $A_{s,k}$.

Let $\Ll'$ and $v$ be as in Lemma~\ref{L:compspeed}. Since the original polar function from Theorem~\ref{T:thm9} will no longer be needed, to simplify notation, we write $\Ll:=\Ll'$ and $\sig:=\sig'$. For each $\lis$ and $1\leq k\leq m$, we choose $\Aa_{s,k}\in\Ai(\phi_k)$ such that the supremum in the definition of the compensated edge speed $\eps^v_{\phi_k}(\Ll_s)$ in (\ref{compedge}) is attained in $\Aa_{s,k}$, i.e.,
\be
\eps^v_\phi(\Ll_s)=\sup_{\Aa\in\Ai(\phi)}\inf_{i\in\Aa}\Ll_s(i-v)=:\inf_{i\in\Aa_{s,k}}\Ll_s(i-v)\qquad(\lis,\ 1\leq k\leq m).
\ee
Using the drift $v$, we lift the linear polar function $\Ll:\R^d\to\R^\sig$ to space-time by setting
\be\label{drift}
L(i,t):=\Ll(i+vt)\qquad(i\in\Z^d,\ t\in\Z),
\ee
and we let $A_{s,k}:=\big\{(i,-1):i\in\Aa_{s,k}\big\}$ denote the space-time set associated with $\Aa_{s,k}$. Then
\be\label{posdrift}
\inf_{(i,t)\in A_{s,k}}L_s(i,t)=\inf_{i\in\Aa_{s,k}}\Ll_s(i-v)=\eps^v_\phi(\Ll_s)>0\qquad(\lis,\ 1\leq k\leq m).
\ee
It remains to choose the subprobability measures $\abf^\star_{s,k}$ from assumption~(viii) and the constant $\la$ from assumption~(vii). For concreteness, we choose for $\abf^\bullet_{s,k}$ the uniform distribution on $A_{s,k}$ and for $\abf^\circ_{s,k}$ the uniform distribution on $\De_k$, although, in fact, any strictly positive probability distributions on these sets will do. It follows from (\ref{posdrift}) that by making the constant $\la$ in (\ref{alpha}) big enough, we can make the constants $\al^\bullet_{s,k}$ as small as we wish. Since the definition of $B_\bullet$ in (\ref{BB}) does not depend on $\la$, this means that we can make the constant $C_\bullet$ from Proposition~\ref{P:randbd} as small as we wish.

We now apply the bound (\ref{albd}) from Proposition~\ref{P:randbd}. Since we can make the constant $C_\bullet$ as small as we wish, the factor with the power $n^\bullet(D)$ is harmless. To deal with the other factor that is raised to the power $n^\circ(D)$, we use that by Lemma~\ref{L:ncirc} and Theorem~\ref{T:randcont}, for any $0<\hat p<1$,
\bc\label{ptilt}
\dis\nu_{p,p_\circ,\rbf,\abf}(D)
&=&\dis\nu_{\hat p,p_\circ,\rbf,\abf}(D)\left(\frac{p}{\hat p}\right)^{n_\circ(D)+1}\left(\frac{1-p}{1-\hat p}\right)^{n_\dia(D)}\\[15pt]
&\leq&\dis\nu_{\hat p,p_\circ,\rbf,\abf}(D)\left(\frac{p}{\hat p}\right)^{n_\circ(D)+1}(1-\hat p)^{-n_\circ(D)-\sig n_\bullet(D)-1}.
\ec
Combining this with (\ref{albd}), we obtain the bound
\be\ba{r@{\,}l}\label{albdhat}
\dis\sum_{D\in\Di_0}\pi_{p,\rbf}(D)\leq\frac{p}{\hat p(1-\hat p)}\sum_{D\in\Di_0}\nu_{\hat p,p_\circ,\rbf,\abf}(D)&\dis\Big(C_\circ p_\circ^{-1}\frac{p}{\hat p(1-\hat p)}\Big)^{n_\circ(D)}\\[15pt]
&\dis\times\Big(C_\bullet(1-p_\circ)^{-1}(1-\hat p)^{-\sig}\Big)^{n_\bullet(D)},
\ec
that holds for any $p_\circ\in[0,1]$. We have already observed that by making $\la$ large enough this means that we can make the constant $C_\bullet$ as small as we wish. In particular, for any $0<\hat p<1$, we can make $C_\bullet$ small enough so that $1-(1-\hat p)^{-\sig}C_\bullet>0$. This allows us to choose $p_\circ:=1-(1-\hat p)^{-\sig}C_\bullet$ so that the factor with the power $n_\bullet(D)$ in (\ref{albdhat}) equals one. Setting $\eps:=\hat p(1-\hat p)p_\circ C_\circ^{-1}$, we then see that for all $p\leq\eps$ the factor with the power $n_\circ(D)$ in (\ref{albdhat}) is $\leq 1$ and hence, using the fact that $\nu_{\hat p,p_\circ,\rbf,\abf}$ is a subprobability measure,
\be
\sum_{D\in\Di_0}\pi_{p,\rbf}(D)\leq\frac{p}{\hat p(1-\hat p)}\qquad(0\leq p\leq\eps).
\ee
Combining this with (\ref{onemin}) we see that $\lim_{p\to 0}\ov\rho(p)=1$ so $\Phi^0$ is stable.
\epro

We now explain the main ideas behind the proofs of Theorems \ref{T:cycbd} and \ref{T:conbd}. The basic idea is the same as for the proof of Theorem~\ref{T:thm9} we have just given, but we can no longer guarantee that directed edges of charge $s$ coming out of an internal vertex or the root \emph{always} point in a direction in which the $s$-th coordinate $L_s$ of the polar function increases. Instead, we must use that they \emph{typically} point in a direction in which $L_s$ increases. In other words, we must distinguish directed edges of charge $s$ coming out of internal vertices into ``good'' edges, for which $L_s$ increases and whose entropic factors we can redistribute to the sinks, and ``bad'' edges. For the sinks, we pay an energy cost because $p$ is small and for the bad edges we pay an energy cost because they are unlikely under the probability distribution $\rbf$.

To quantify the energy cost we pay for sinks, in (\ref{ptilt}) we compared two subprobability measures on the set of decorated Toom contours: one with the original parameter $p$ and one with a ``tilted'' parameter $\hat p$. This allowed us in (\ref{albdhat}) to obtain a factor $p/\hat p(1-\hat p)$ with a power $n_\circ(D)$. In a similar way, we can ``tilt'' the probability measure $\hat\rbf$ to quantify the energy cost we pay for ``bad'' edges. We can think of the charges of a Toom contour as moving according to random walks. If we have chosen the polar function $L$ and the sets $A_{s,k}$ from assumption~(i) of Subsection~\ref{S:abstract} and from (\ref{Askchoice}) correctly, then each charge $s$ will typically move in the direction in which the corresponding coordinate $L_s$ of the polar function increases. We are interested in exceptional behaviour of this random walk in the sense of large deviation theory. We will use the ``tilted''  probability measures $\hat\rbf_t$ from assumption~(iii) of Subsection~\ref{S:abstract} to quantify the energy cost we pay for having too many ``bad'' edges, in a way that is similar to the use of the Cram\'er transform in large deviation theory (compare \cite[Lemmas~I.8--I.10]{Hol00}).

When we compare the subprobability measures $\nu_{p,p_\circ,\rbf,\abf}$ and $\nu_{p,p_\circ,\hat\rbf,\abf}$ (with $\rbf$ replaced by $(\hat\rbf)_{t\in\{\circ,1,\ldots,\sig\}}$), there is one difficulty we have to deal with: unlike sinks, internal vertices can overlap with other vertices. This is why we introduced the function $\tau$ in (\ref{taudef}) that for overlapping vertices gives the type of the first vertex (according to the natural order on the Toom graph) that is embedded at a given space-time point. When we compare $\nu_{p,p_\circ,\rbf,\abf}$ with $\nu_{p,p_\circ,\hat\rbf,\abf}$ similar to what we did in (\ref{ptilt}), we obtain a factor $\rbf(\mu(i))/\hat\rbf_{\tau(i)}(\mu(i))$ only for the first vertex that is embedded at $i$. For Toom cycles, we can use (\ref{cycord}) which says that vertices of type $1$ are embedded before vertices of types $\circ$ and $2$. In this case, it makes sense to tilt only $\hat\rbf_1$ and choose $\hat\rbf_\circ=\hat\rbf_2=\rbf$ as we do in Theorem~\ref{T:cycbd}. For general Toom contours, we have to do a worst-case bound since for overlapping vertices of different types we do not know which is the first in the natural order on the Toom graph. This is the reason why we introduce the constants $\bet^{s'}_{s,k}$ in (\ref{bet3}). Although this is a somewhat crude method, it is still good enough to prove nontrivial results such as Theorem~\ref{T:triang} that at present cannot be obtained otherwise.

\subsection{Infinity of the Peierls sum}\label{S:infin}

As explained at the end of the last subsection, we tilt the parameter $p$ of the subprobability measure $\nu_{p,p_\circ,\rbf,\abf}$ to estimate the energy cost of sinks and we tilt the probability distribution $\rbf$ to estimate the energy cost of charged edges that point in the ``wrong'' direction, compared to the typical edge speed of the process. Since we are interested in the limit $p\to 0$, we can choose $p$ as small as we wish to make the energy cost for sinks as large as we wish. By contrast, the probability distribution $\rbf$ is fixed (for example, in Theorem~\ref{T:triang} it is the uniform distribution on $\{1,2,3\}$), which means that there is an upper bound on how large we can make the energy cost of charged edges that point in the wrong direction. In the present subsection we show that this is not a technical detail but a fundamental limitation of our method.

Throughout this subsection we are in the setting of~\eqref{Zd}, that is $\La=\Z^d\times\Z$. We will consider a monotone cellular automaton $\Phi^{p,r}$ that applies the zero map with probabilty $p$ and a shrinker $\phh$ and the identity map $\phh^{\rm id}$ with probabilities $(1-p)r$ and $(1-p)(1-r)$, respectively. The situation is the same as in Theorem~\ref{T:coop} except that instead of $\phh^{\rm coop}$ we will use a different shrinker. In line with Conjecture~\ref{C:shrink}, we expect $\Phi^{0,r}$ to be stable for each $r>0$. Nevertheless, for a natural choice of the typed dependence graph $(\La,\Hi^{p,r})$, for $r$ small enough, we will show that the Peierls sum over Toom cycles in (\ref{Peicyc}) is infinite for all $p>0$. This shows that Conjecture~\ref{C:shrink} cannot be proved with the methods of the present paper.

To describe the shrinker we will be interested in, we define sets $\Aa_{s,k}\sub\Z^2$ $(1\leq s,k\leq 2)$ by
\be\label{Acount}
\Aa_{1,1}:=\big\{(0,1),(1,0)\big\},\quad
\Aa_{2,1}:=\big\{(-1,1),(0,0),(1,-1)\big\},\quad \Aa_{1,2}=\Aa_{2,2}:=\big\{(0,0)\big\}.
\ee
We define ``space'' maps
\be\label{phicc}
\phi^{\rm cc}(x):=\bigvee_{s=1}^2\bigwedge_{i\in\Aa_{s,1}}x(i)
\quand
\phi^{\rm id}(x):=\bigvee_{s=1}^2\bigwedge_{i\in\Aa_{s,2}}x(i),
\ee
whose associated ``space-time'' maps as defined in (\ref{phhphi}) are $\phh^{\rm cc}$ and $\phh^{\rm id}$ respectively. Then $\phi^{\rm id}$ is the identity map that just copies the state that was previously at $0$. We claim that $\phi^{\rm cc}$ is a shrinker. To see this, consider the linear polar function $\Ll:\R^2\to\R^2$ defined as 
\be\label{Lcoop2}
\Ll_1(z):=z_1+z_2\quand \Ll_2(z):=-z_1-z_2\qquad\big(z=(z_1,z_2)\in\R^2\big).
\ee
It is straightforward to check that the edge speeds from (\ref{edge}) are
\be
\eps_{\phi^{\rm cc}}(\Ll_1)=1\quand\eps_{\phi^{\rm cc}}(\Ll_2)=0,
\ee
which by (\ref{shrinker}) shows that $\phi^{\rm cc}$ is a shrinker.

We set $\La:=\Z^3$ and let $A_{s,k}:=\big\{(i_1,i_2,-1):(i_1,i_2)\in\Aa_{s,k}\big\}$ denote the space-time set associated with $\Aa_{s,k}$ $(1\leq s,k\leq 2)$. For each $p,r\in[0,1]$, we let $(\mu(i))_{i\in\La}$ be i.i.d.\ random variables with values in $\{0,1,2\}$ such that
\be
\P\big[\mu(i)=0\big]=p,\quad\P\big[\mu(i)=1\big]=(1-p)r,\quand\P\big[\mu(i)=2\big]=(1-p)(1-r),
\ee
and we let $(\La,\Hi^{p,r})$ denote the typed dependence graph defined in terms of the sets $A_{s,k}$ and random variables $\mu(i)$ as in (\ref{LaHi}). We let $\Phi^{p,r}$ denote the monotone cellular automaton associated with $(\La,\Hi^{p,r})$. In view of (\ref{phicc}), this is the monotone cellular automaton that applies the zero map with probability $p$ and the maps $\phh^{\rm cc}$ and $\phh^{\rm id}$ with probabilities $(1-p)r$ and $(1-p)(1-r)$, respectively. In Subsection~\ref{S:diverge} we will prove the following result. Recall from Theorem~\ref{T:Peierls} that $\bar\Ti_0$ denotes the set of Toom cycles rooted at $0$.

\bp[Infinity of the Peierls sum]
There\label{P:Peinf} exists an $r'>0$ such that for all $0<r\leq r'$, there exists a $p'>0$ such that
\be\label{Peinf}
\sum_{T\in\bar\Ti_0}\P\big[T\mbox{ is present in }(\La,\Hi^{p,r})\big]=\infty\quad\forall 0<p\leq p'.
\ee
\ep

It seems plausible that (\ref{Peinf}) actually holds for all $0<p<1$ but our proof only yields the statement for $p$ small enough and there seems to be no obvious monotonicity in $p$ that would allow us to extend (\ref{Peinf}) to all $0<p<1$. In fact, (\ref{Peinf}) does not hold for $p=1$ since in this case there is a.s.\ a single Toom cycle present in $(\La,\Hi^{1,r})$.

The intuition behind Proposition~\ref{P:Peinf} has already been explained at the beginning of this subsection. As long as $r>0$, we can think of the 1-charges move as drifted random walks that tend to move in the direction where the function $L_1$ is increasing. As a result, the probability of finding a large contour in which the 1-charges have not moved much is exponentially small. For small $r$ this exponential factor is beaten, however, by the combinatorical factor that comes from counting all contours with this property, which results in the Peierls sum being infinite for all $p>0$.

There is still a considerable gap between, on the one hand, cellular automata to which Theorems \ref{T:cycbd} and \ref{T:conbd} are applicable, and, on the other hand, cellular automata such as the one in Proposition~\ref{P:Peinf} for which we can prove that the Peierls sum is infinite for $p>0$ small enough. Recall that $\phh^{\rm id}$ is defined in (\ref{coopid}) and $\phh^{\rm NEC}$ is defined in (\ref{NEC}). In analogy to $\phh^{\rm NEC}$, one can also define a \emph{North West Center} map $\phh^{\rm NWC}$ (compare the \emph{South West Center} map $\phh^{\rm SWC}$ defined in (\ref{phitri})). Examples of monotone cellular automata for which we do not know whether the Peierls sum tends to zero as $p\to 0$ are:
\begin{itemize}
\item The cellular automaton that applies the maps $\phh^{\rm NEC}$ and $\phh^{\rm id}$ with probabilities $r$ and $1-r$, respectively, in the regime where $r$ is small.
\item The cellular automaton that applies the maps $\phh^{\rm NEC}$ and $\phh^{\rm NWC}$ with equal probabilities.
\end{itemize}
Our understanding is better for Toom cycles than for general Toom contours, thanks to (\ref{cycord}). The fact that internal vertices of different charges can overlap causes much trouble in proving positive results but seems of little help in proving negative results.

It is worth noting that Proposition~\ref{P:Peinf} is, in a sense, an annealed result. We can view monotone cellular automata with intrinsic randomness a having two sources of randomness: one that describes the unperturbed cellular automaton $\Phi^0$ by deciding at which space-time points each of the non-constant maps $\phh_1,\ldots,\phh_m$ is applied, and another that perturbs $\Phi^0$ by making each space-time point defective with a small probability $p$. Similarly, one can view an associated typed dependence graph $(\La,\Hi^p)$ as being obtained from $(\La,\Hi^0)$ by the addition of some extra randomness that makes certain points defective. One can then define an \emph{annealed} and \emph{quenched} Peierls sum, which give the expected number of Toom contours rooted at $0$ that are present in $(\La,\Hi^p)$ and their conditional expectation given $(\La,\Hi^0)$, respectively. In this language, Proposition~\ref{P:Peinf} says that if $r$ is small enough, then the annealed Peierls sum is infinite for all $p$ small enough. The proof of Proposition~\ref{P:Peinf} makes essential use of configurations that are exponentially unlikely under the unperturbed probability law. This leaves open the possibility that almost surely, the quenched Peierls sum is finite for small $p$ and tends to zero as $p\to 0$, which by Theorem~\ref{T:contour} would still be enough to prove stability of $\Phi^0$.

\section{Random contours}\label{S:nuproof}

\subsection{Inductive construction of Toom graphs}\label{S:incomp}

In this section we prove Theorem~\ref{T:randcont}. We will construct a random variable $\Dec$ with values in $\Di_0\cup\{\dgg\}$, where $\dgg$ is a ``cemetery state'' not included in $\Di_0$, so that its law restricted to $\Di_0$ is precisely the subprobability measure defined in (\ref{DecD}). We will construct the random decorated Toom contour $\Dec$ step by step, by adding edges and vertices to an incomplete decorated Toom contour, one at a time, in a Markovian way. The construction terminates as soon as we have a complete decorated Toom contour. Each step can fail with a certain probability, in which case we end up in the cemetery state $\dgg$. The Markovian construction may also run for infinite time without ever terminating, in which case we also end up in $\dgg$. Although the basic idea of the construction is not difficult, the proof of Theorem~\ref{T:randcont} requires a rather unwieldy lot of definitions. The charge $\sig$ will play a special role in our contruction. Edges of charges $1,\ldots,\sig-1$ will be added in an ``upward'' fashion while edges of charge $\sig$ will be added in a ``downward'' fashion.

\begin{figure}[t]
\begin{center}
\inputtikz{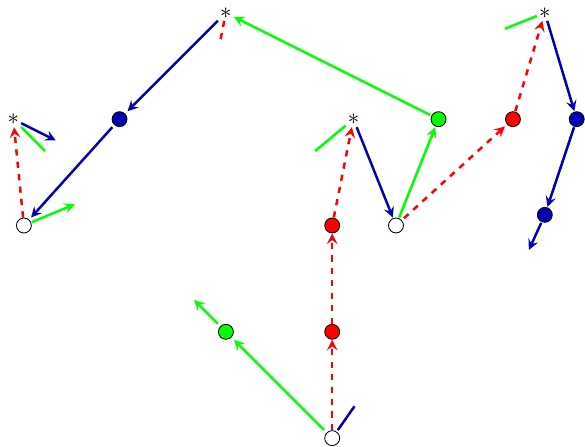}
\caption{An incomplete Toom graph with charges $1=$ green, $2=$ red, and $3=$ blue. To stress the special role of charge 3, the direction of edges with this charge has been reversed in the picture. Loose ends are indicated by small arrows and dead ends by small line segments.}
\label{fig:Toomincomp}
\end{center}
\end{figure}

We start by defining incomplete Toom graphs. We only need rooted connected Toom graphs, so we make the rootedness and connectedness part of our definition. It will also be convenient to assume that the vertex set is of the form $V=\{0,\ldots,N\}$ for some $N\geq 0$, where $0$ is the root. By definition, the \emph{trivial Toom graph} with $\sig$ charges is the Toom graph $(\Vi,\Ei)$ defined as
\be\label{triv}
V_\circ=V_\ast=\{0\},\ V_s=\emptyset=\vec E_s\quad(\lis).
\ee
The following definitions are demonstrated in Figure~\ref{fig:Toomincomp}.

\begin{defi}
An\label{D:Toomincomp} \emph{incomplete Toom graph} with $\sig\geq 2$ charges is a connected typed directed graph $(\Vi,\Ei)$ with vertex type set $\{\circ,\ast,1,\ldots,\sig\}$, edge type set $\{1,\ldots,\sig\}$, and vertex set of the form $V=\{0,\ldots,N\}$ for some $N\geq 0$, with $0\in V_\circ$, that satisfies the conditions
\begin{enumerate}
\item $|\vec E_{s,{\rm in}}(v)|=0$ and $|\vec E_{s,{\rm out}}(v)|\leq 1$ $(\lis)$ for all $v\in V_\circ$,
\item $|\vec E_{s,{\rm out}}(v)|=0$ and $|\vec E_{s,{\rm in}}(v)|\leq 1$ $(\lis)$ for all $v\in V_\ast$,
\item $|\vec E_{s,{\rm in}}(v)|\leq 1$, $|\vec E_{s,{\rm out}}(v)|\leq 1$, and $|\vec E_{l,{\rm in}}(v)|=0=|\vec E_{l,{\rm out}}(v)|$ for each $\lis$, $l\neq s$ and $v\in V_s$,
\end{enumerate}
and that is either the trivial Toom graph or satisfies the additional conditions
\begin{enumerate}\addtocounter{enumi}{3}
\item $|\vec E_{\sig,{\rm out}}(v)|=1$ for each $v\in(V_\sig\cup V_\circ)\beh\{0\}$,
\item $|\vec E_{s,{\rm in}}(v)|=1$ for each $1\leq s<\sig$ and $v\in V_s$,
\item for each $v\in V_\ast$, there exists an $1\leq s<\sig$ such that $|\vec E_{s,{\rm in}}(v)|=1$.
\end{enumerate}
We call $0$ the \emph{root} of $(\Vi,\Ei)$ and in line with (\ref{toomroot}) we set $V'_\circ:=V_\circ\beh\{0\}$ and $V'_s:=V_s\cup\{0\}$.
\end{defi}

\begin{defi}\label{d:loose}
A\label{D:loose} \emph{loose end} of an incomplete Toom graph $(\Vi,\Ei)$ is a pair $(v,s)$ where $v\in V$ and $\lis$ satisfy one of the following conditions:
\begin{enumerate}
\item $v\in V_\circ\cup V_s$, $1\leq s<\sig$, and $|\vec E_{s,{\rm out}}(v)|=0$,
\item $v\in V_\ast\cup V_\sig$, $s=\sig$, and $|\vec E_{\sig,{\rm in}}(v)|=0$.
\end{enumerate}
We call $s$ the \emph{charge} of the loose end. The \emph{most urgent} loose end of an incomplete Toom graph is the loose end $(v,s)$ that has the smallest possible charge $s$, and among all loose ends with charge $s$, has the smallest possible vertex $v\in\{0,\ldots,N\}$.
\end{defi}

\bl[Complete Toom graphs]
An\label{L:loose} incomplete Toom graph that is not the trivial Toom graph is a Toom graph if and only if it has no loose ends.
\el

\bpro
It is clear that a Toom graph that is not the trivial Toom graph has no loose ends. (This is the only place in the proof where nontriviality is used.) It remains to show that conversely, if an incomplete Toom graph has no loose ends, then it is a Toom graph. Let us define a \emph{dead end} of an incomplete rooted Toom graph $(\Vi,\Ei)$ to be a pair $(v,s)$ where $v\in V$ and $\lis$ satisfy one of the following conditions:
\begin{enumerate}\addtocounter{enumi}{2}
\item $v=0$, $s=\sig$, and $|\vec E_{s,{\rm out}}(v)|=0$,
\item $v\in V_\ast$, $1\leq s<\sig$, and $|\vec E_{s,{\rm in}}(v)|=0$.
\end{enumerate}
We call $s$ the \emph{charge} of the dead end. Then it is immediate from the definitions that an incomplete rooted Toom graph is a rooted Toom graph if and only if it has no loose ends or dead ends. To complete the proof, it suffices to show that if an incomplete rooted Toom graph has no loose ends, then it has no dead ends either.

To see this, we observe that in each source $v$, for each charge $1\leq s<\sig$, there starts a unique directed path of edges of charge $s$ that has length zero if $(v,s)$ is a loose end and otherwise either ends in an internal vertex $w$ such that $(w,s)$ is a loose end, or in a sink. Let say that $v$ is \emph{$s$-satisfied} if this path ends in a sink. Similarly, each sink $v$ is the endvertex of a unique directed path of edges of charge $\sig$ that has length zero if $(v,\sig)$ is a loose end and otherwise either starts in an internal vertex $w$ such that $(w,\sig)$ is a loose end, or in a source. Let us say that $v$ is $\sig$-\emph{satisfied} if this path starts in a source.

Then the number of dead ends of charge $1\leq s<\sig$ equals the number of sinks minus the number of $s$-satisfied sources, and the number of dead ends of charge $\sig$ (which can be either zero or one by condition (iii) in the definition of a dead end) equals the number of sources minus the number of $\sig$-satisfied sinks. If there are no loose ends, then every source is $s$-satisfied for each $1\leq s<\sig$, and every sink is $\sig$-satisfied. So in this case, the number of dead ends of charge $1\leq s<\sig$ equals the number of sinks minus the number of sources, and the number of dead ends of charge $\sig$ equals the number of sources minus the number of sinks. Since both quantities must be nonnegative, we conclude that the number of sources equals the number of sinks and there are no dead ends.
\epro

Recall that formally, $\Ei$ is a set with elements of the form $(v,w,s)$, which have the interpretation of a directed edge of charge $s$ that starts in the vertex $v$ and ends in~$w$.

\begin{defi}
Let $(v,s)$ be a loose end of an incomplete Toom graph $(\Vi,\Ei)$. An \emph{extension} of $(\Vi,\Ei)$ at $(v,s)$ is an incomplete Toom graph $(\Vi',\Ei')$ such that $\Vi\sub\Vi'$, $\Ei\sub\Ei'$, and $\Ei'\beh\Ei$ consists of a single element, which is of the form $(v,w,s)$ if $1\leq s<\sig$ and of the form $(w,v,\sig)$ if $s=\sig$.
\end{defi}

Recall the definition of the trivial Toom graph in (\ref{triv}). By definition, the \emph{trivial incomplete Toom graph} is the incomplete Toom graph $(\Vi,\Ei)$ defined as
\be\label{trivinc}
V_\circ:=\{0\},\ V_\ast:=\emptyset,\mbox{ and }V_s:=\emptyset=:\vec E_s\quad(\lis).
\ee
Note that this differs from the trivial Toom graph in (\ref{triv}) only since $V_\ast\neq\{0\}$.

The following lemma says that each rooted connected Toom graph has a unique inductive construction, starting at the root and adding edges once at a time at the most urgent loose end. Recall the definition in (\ref{psiEi}) of the image $\psi(\Vi,\Ei)$ of a typed directed graph $(\Vi,\Ei)$ under a map $\psi$. The equality $(\Vi(e),\Ei(e)):=\eta^{-1}(\Vi'(e),\Ei'(e))$ below means that $(\Vi(e),\Ei(e))$ is obtained from $(\Vi'(e),\Ei'(e))$ by relabeling the vertices of $(\Vi,\Ei)$ with $\eta^{-1}$.

\bl[Inductive construction]
Let\label{L:induct} $(v_\circ,\Vi,\Ei)$ be a rooted connected Toom graph that is not the trivial Toom graph, let $N:=|V|-1$, and $M:=|\Ei|$. Then there exists a unique bijection $\eta:V\to\{0,\ldots,N\}$ and incomplete Toom graphs
\be\label{indinc}
\big(\Vi'(e),\Ei'(e)\big)_{0\leq e\leq M},
\ee
such that setting $(\Vi(e),\Ei(e)):=\eta^{-1}(\Vi'(e),\Ei'(e))$ $(0\leq e\leq M)$, one has
\begin{enumerate}
\item $\eta(v_\circ)=0$,
\item $(\Vi'(0),\Ei'(0))$ is the trivial incomplete Toom graph and $(\Vi(M),\Ei(M))=(\Vi,\Ei)$,
\item for all $0<e\leq M$, the incomplete Toom graph $(\Vi'(e),\Ei'(e))$ is an extension of\\ $(\Vi'(e-1),\Ei'(e-1))$ at its most urgent loose end.
\end{enumerate}
\el

\bpro
We claim that for each $0\leq\ti M\leq M$, there exists an $0\leq\ti N\leq N$ and an injective function $\eta^{-1}:\{0,\ldots,\ti N\}\to V$ and incomplete Toom graphs
\be
\big(\Vi'(e),\Ei'(e)\big)_{0\leq e\leq\ti M},
\ee
such that setting $(\Vi(e),\Ei(e)):=\eta^{-1}(\Vi'(e),\Ei'(e))$ $(0\leq e\leq\ti M)$, one has
\begin{itemize}
\item[${\rm(i)'}$] $\eta(v_\circ)=0$,
\item[${\rm(ii)'}$] $(\Vi'(0),\Ei'(0))$ is the trivial incomplete Toom graph and $(\Vi(\ti M),\Ei(\ti M))$ is a typed subgraph of $(\Vi,\Ei)$,
\item[${\rm(iii)'}$] for all $0<e\leq\ti M$, the incomplete Toom graph $(\Vi'(e),\Ei'(e))$ is the extension of\\ $(\Vi'(e-1),\Ei'(e-1))$ at its most urgent loose end.
\end{itemize}
We prove our claim by induction on $\ti M$. The statement is clearly true for $\ti M=0$. Assume that it is true for some $0\leq\ti M<M$. We claim that $(\Vi(\ti M),\Ei(\ti M))$ has a loose end. In the opposite case, by Lemma~\ref{L:loose} it would be a complete Toom graph. By our assumption that  $(\Vi,\Ei)$ is connected, it follows that either 1.\ $(\Vi(\ti M),\Ei(\ti M))$ is a trivial Toom graph or 2.\ $(\Vi(\ti M),\Ei(\ti M))=(\Vi,\Ei)$. Condition~(iii)' implies that $|\Ei(\ti M)|=\ti M$, so case~1 can only occur for $\ti M=0$, where it contradicts (ii)', while case~2 contradicts the assumption $\ti M<M$.

Let $(v',s)$ be the most urgent loose end of $(\Vi'(\ti M),\Ei'(\ti M))$ and let $(v,s)$ with $v:=\eta^{-1}(v')$ be the corresponding loose end of $(\Vi(\ti M),\Ei(\ti M))$. By the assumption that $(\Vi,\Ei)$ is a Toom graph, it follows that:
\begin{itemize}
\item If $1\leq s<\sig$, then there exists a unique $w\in V$ such that $(v,w)\in\vec E_s$.
\item If $s=\sig$, then there exists a unique $w\in V$ such that $(w,v)\in\vec E_\sig$.
\end{itemize}
From this, it is easy to see that the induction hypothesis holds for $\ti M+1$. Indeed, we need to choose $\Ei(\ti M+1):=\Ei(\ti M)\cup\{(v,w,s)\}$ if $1\leq s<\sig$ and $:=\Ei(\ti M)\cup\{(w,v,s)\}$ if $s=\sig$. Further, we need $V(\ti M+1):=V(\ti M)\cup\{w\}$. If $w\not\in V(\ti M)$, then we need to assign it the same type as it has in $(\Vi,\Ei)$ and set $\eta(w):=\ti N+1$. This completes the proof of our claim.

To see that our claim for $\ti M=M$ implies the claim of the lemma, it suffices to note that by condition~(iii)' $|\Ei(M)|=M$ and hence by condition~(ii)' $(\Vi(M),\Ei(M))=(\Vi,\Ei)$.
\epro

Using Lemma~\ref{L:induct}, we can unambiguously equip the vertex set $V$ of any rooted connected Toom graph $(v_\circ,\Vi,\Ei)$ with a total order by setting
\be\label{natord}
v\leq w\quad\desd\quad\eta(v)\leq\eta(w)\qquad(v,w\in V).
\ee
We call this the \emph{natural order} on~$V$. We make the following observation.

\bl[The natural order on Toom cycles]
Let\label{L:cycord} $(v_\circ,\Vi,\Ei,\psi)$ be the Toom contour associated with a Toom cycle as in (\ref{cycont}). Then for the rooted connected Toom graph $(v_\circ,\Vi,\Ei)$, the natural order on $V=\{0,\ldots,n-1\}$ is the order $0<1<\cdots<n-1$.
\el

\bpro
Obvious from the fact that an incomplete Toom contour with two charges has at most one loose end.
\epro

\subsection{Random extensions of incomplete contours}\label{S:context}

We extend Definition~\ref{D:embed} of an embedding of a Toom graph and Definition~\ref{D:decor} of a decorated Toom contour in a straightforward way to incomplete Toom graphs. Recall that for incomplete Toom contours we assume that the vertex set is of the form $V=\{0, \dots, N\}$ where 0 is the root.

\begin{defi}\label{D:embed2}
Let $(\Vi,\Ei)$ be an incomplete Toom graph and let $\La$ be a countable set. By definition, an \emph{embedding} of $(\Vi,\Ei)$ in $\La$ is a map $\psi:V\to\La$ such that:
\begin{enumerate}
\item $\psi(v_1)\neq\psi(v_2)$ for each $v_1\in V_\ast$ and $v_2\in V$ with $v_1\neq v_2$,
\item $\psi(v_1)\neq\psi(v_2)$ for each $v_1,v_2\in V'_s$ with $v_1\neq v_2$ $(\lis)$.
\end{enumerate}
\end{defi}

Starting with the following definition, $\La$ is a countable group, $h$ is a height function, and $(A_{s,k})_{\lis,\ 1\leq k\leq m}$, with $\sig\geq 2$ and $m\geq 1$, are finite subsets of $\{i\in\La:h(i)=1\}$.

\begin{defi}\label{D:incomp}
An \emph{incomplete decorated Toom contour} is a quadruple $(\Vi,\Ei,\psi,\kappa)$ where $(\Vi,\Ei)$ is an incomplete Toom graph, $\psi$ is an embedding of $(\Vi,\Ei)$ in $\La$, and $\kappa:\psi(V)\to\{0,\ldots,m\}$ is a function such that:
\begin{enumerate}
\item If $v\in V_\ast$ and $\psi(v)=i$, then $\kappa(i)=0$,
\item If $v\not\in V'_\circ$, $(v,w)\in\vec E_s$, and $\psi(v)=j$, then $\psi(w)=ij$ for some $j\in A_{s,\kappa(i)}$,
\item If $v\in V'_\circ$, $(v,w)\in\vec E_s$, and $\psi(v)=i$, then $\psi(w)=ij$ for some $j\in\De_{\kappa(i)}$.
\end{enumerate}
We say that the contour $(\Vi,\Ei,\psi,\kappa)$ is \emph{rooted} at $i_\circ:=\psi(0)$. We let $\Di'_0$ denote the space of all incomplete decorated Toom contours rooted at $0\in\La$.
\end{defi}

We now define three ways to construct a random extension of an incomplete decorated Toom contour $(\Vi,\Ei,\psi,\kappa)$ at a given loose end $(v,s)$ of $(\Vi,\Ei)$: one that applies when $1\leq s<\sig$, and two that apply when $s=\sig$ and that depend on whether we want to add a new source or not. The constructions depend on a number of parameters that can be freely chosen:
\begin{itemize}
\item a constant $\hat p\in[0,1]$,
\item for each $t\in\{\circ,1,\ldots,\sig\}$, a probability distribution $\hat\rbf_t$ on $\{1,\ldots,m\}$,
\item for each $\lis$ and $1\leq k\leq m$, a subprobability distribution $\abf^\bullet_{s,k}$ on $A_{s,k}$,
\item  for each $1\leq k\leq m$, a subprobability distribution $\abf^\circ_{s,k}$ on $\De_k$.
\end{itemize}
Formally, the definitions below define a random variable $\Dec$ with values in $\Di'_0\cup\{\dgg\}$, where $\dgg$ is a cemetery state not included in $\Di'_0$. We use the informal phrase ``the extension fails'' to say that $\Dec=\dgg$. Sometimes, we say that we draw a certain object at random from a subprobability distribition. In such a case, it is understood that by adding a cemetery state, we extend our subprobability distribition to a probability distribution, we draw the object from this extended probability distribution, and the extension fails if we draw the cemetery state.

\begin{defi}
Let\label{D:forext} $D=(\Vi,\Ei,\psi,\kappa)\in\Di'_0$ be an incomplete decorated Toom contour and let $(v,s)$ be a loose end of charge $1\leq s<\sig$. Let $V=\{0,\ldots,N\}$, $(W,\Fi):=\psi(\Vi,\Ei)$, and $i:=\psi(v)$. By definition, if $v\not\in V'_\circ$, then the \emph{upward extension} of $D$ at $(v,s)$ is constructed as follows. First, we draw $j$ from $A_{s,\kappa(i)}$ from the subprobability law $\abf^\bullet_{s,\kappa(i)}$. If $\vec F_{s,{\rm in}}(ij)\neq\emptyset$, then the extension fails. If $ij=\psi(w)$ for some (necessarily unique) $w\in V_\ast$, then we add the edge $(v,w)$ to $\vec E_s$ and we are done. In the opposite case, we add the new element $w:=N+1$ to the set $V_s$, we add the edge $(v,w)$ to $\vec E_s$, and we set $\psi(w):=ij$. If $ij\in W$, then we are done. In the opposite case, we draw $\kappa(ij)$ from the probability law
\be\label{Pkap}
\P\big[\kappa(ij)=k\big]=\left\{\ba{ll}
\dis\hat p\quad&\mbox{if }k=0,\\[5pt]
\dis(1-\hat p)\hat\rbf_s(k)\quad&\mbox{if }1\leq k\leq m,
\ea\right.
\ee
and we are done. If $v\in V'_\circ$, then we proceed in exactly the same way, except that we draw $j$ from the subprobability law $\abf^\circ_{s,\kappa(i)}$, and in (\ref{Pkap}) we replace $\hat\rbf_s$ by $\hat\rbf_\circ$.
\end{defi}

We next turn our attention to loose ends of charge $\sig$. These are a bit more complicated, since edges need to be added in a downward fashion, and we have a choice whether at the end of such an edge we add an internal vertex or a source. We start by describing how to add internal vertices. The fact that (\ref{musig}) defines a subprobability distribution is not obvious and will be proved below.

\begin{defi}
Let\label{D:downext} $D=(\Vi,\Ei,\psi,\kappa)\in\Di'_0$ be an incomplete decorated Toom contour and let $(v,\sig)$ be a loose end of charge $\sig$. Let $V=\{0,\ldots,N\}$, $(W,\Fi):=\psi(\Vi,\Ei)$, and $i:=\psi(v)$. We define a subprobability distribution $\mu_\bullet$ on $\La\times\{1,\ldots,m\}$ by
\be\label{musig}
\mu_\bullet(j,k):=\left\{\ba{ll}
\dis B_\bullet^{-1}\abf^\bullet_{\sig,k}(j)\quad&\dis\mbox{if }ij^{-1}\in W,\ k=\kappa(ij^{-1}),\mbox{ and }j\in A_{\sig,k},\\[5pt]
\dis B_\bullet^{-1}\hat\rbf_\sig(k)\abf^\bullet_{\sig,k}(j)\quad&\dis\mbox{if }ij^{-1}\not\in W\mbox{ and }j\in A_{\sig,k},\\[5pt]
\dis 0\quad&\mbox{in all other cases.}
\ea\right.
\ee
By definition, the \emph{downward extension} of $D$ at $(v,s)$ is now constructed as follows. First, we draw $(j,k)$ from the subprobability measure $\mu_\bullet$. If $ij^{-1}\neq 0$, then we add the new element $w:=N+1$ to $V_\sig$, we add the edge $(w,v)$ to $\vec E_\sig$, and we set $\psi(w):=ij^{-1}$. If $ij^{-1}\not\in W$, then we set $\kappa(ij^{-1}):=k$. If there already exists a $w'\in V_\sig\cup V_\ast$ such that $\psi(w')=ij^{-1}$, then the extension fails. If $ij^{-1}=0$, then we proceed similarly, except when $|\vec E_{\sig,{\rm out}}(0)|=0$. In the latter case, we just add the edge $(0,v)$ to $\vec E_\sig$ and we are done.
\end{defi}

To show that this is a good definition, we need to prove that (\ref{musig}) defines a subprobability distribution. To see this, we extend the probability distributions $\abf^\bullet_{\sig,k}$ to $\La$ by setting $\abf^\bullet_{\sig,k}(j):=0$ if $j\not\in A_{\sig,k}$. We claim that for each $j\in\La$,
\be\label{musumj}
\sum_{k=1}^m\mu_\bullet(j,k)\leq B_\bullet^{-1}\sup_{0\leq k\leq m}\abf^\bullet_{\sig,k}(j).
\ee
Indeed, if $ij^{-1}\in W$, then this follows from the fact that $\sum_{k=1}^m\mu_\bullet(j,k)=B_\bullet^{-1}\abf^\bullet_{\sig,\kappa(ij^{-1})}(j)$, and on the other hand, if $ij^{-1}\not\in W$, then this follows from the fact that $\sum_{k=1}^m\mu_\bullet(j,k)=B_\bullet^{-1}\sum_{k=1}^m\hat\rbf_\sig(k)\abf^\bullet_{\sig,k}(j)$ and $\hat\rbf_\sig$ is a probability distribution. Summing (\ref{musumj}) over $j\in\La$, recalling the definition of $B_\bullet$ in (\ref{BB}), we see that $\mu_\bullet$ is a subprobability distribution.

We finally describe how to add new sources to an incomplete contour.

\begin{defi}
Let\label{D:sourext} $D=(\Vi,\Ei,\psi,\kappa)\in\Di'_0$ be an incomplete decorated Toom contour and let $(v,\sig)$ be a loose end of charge $\sig$. Let $V=\{0,\ldots,N\}$, $(W,\Fi):=\psi(\Vi,\Ei)$, and $i:=\psi(v)$. We define a subprobability distribution $\mu_\circ$ on $\La\times\{1,\ldots,m\}$ by
\be\label{mu}
\mu_\circ(j,k):=\left\{\ba{ll}
\dis B_\circ^{-1}\abf^\circ_{\sig,k}(j)\quad&\dis\mbox{if }ij^{-1}\in W,\ k=\kappa(ij^{-1}),\mbox{ and }j\in\De_k,\\[5pt]
\dis B_\circ^{-1}\hat\rbf_\sig(k)\abf^\circ_{\sig,k}(j)\quad&\dis\mbox{if }ij^{-1}\not\in W\mbox{ and }j\in\De_k,\\[5pt]
\dis 0\quad&\mbox{in all other cases.}
\ea\right.
\ee
By definition, the \emph{source extension} of $D$ at $(v,s)$ is now constructed as follows. First, we draw $(j,k)$ from the subprobability measure $\mu_\circ$. We then add the new element $w:=N+1$ to $V_\circ$, we add the edge $(w,v)$ to $\vec E_\sig$ and we set $\psi(w):=ij^{-1}$. If $ij^{-1}\in W$, then we are done. In the opposite case, we set $\kappa(ij^{-1}):=k$ and we are done.
\end{defi}

The proof that (\ref{mu}) defines a subprobability distribution is completely analogous to the proof for (\ref{musig}), so we skip it.

\subsection{Markovian construction of contours}\label{S:Marcon}

Recall from Definitions \ref{D:decor} and \ref{D:incomp} that $\Di_0$ is the set of decorated Toom contours rooted at $0\in\La$ (up to isomorphism), and $\Di'_0$ is the space of all incomplete decorated Toom contours rooted at $0\in\La$. By Lemma~\ref{L:loose}, we can identify  $\Di_0$ with the subspace of $\Di'_0$ consisting of all incomplete decorated Toom contours rooted at $0$ that have no loose ends. As before, we let $\dgg$ denote a cemetery state not included in $\Di'_0$. We define a trivial decorated Toom contour $D^{\rm triv}_0$ and incomplete decorated Toom contours $D^{\rm triv}_1,\ldots,D^{\rm triv}_m$ by
\be\ba{r@{\,}c@{\,}ll}
\dis D^{\rm triv}_0&:=&\dis(\Vi,\Ei,\psi,\kappa)\quad\mbox{with}\quad
&\dis V_\circ=V_\ast=\{0\},\ V_s=\emptyset=\vec E_s\quad(\lis)\\[5pt]
&&&\dis\psi(0)=0,\ \kappa(0)=0,\\[10pt]
\dis D^{\rm triv}_k&:=&\dis(\Vi,\Ei,\psi,\kappa)\quad\mbox{with}\quad
&\dis V_\circ:=\{0\},\ V_\ast:=\emptyset,\ V_s:=\emptyset=:\vec E_s\quad(\lis)\\[5pt]
&&&\dis\psi(0)=0,\ \kappa(0)=k\qquad(1\leq k\leq m).
\ec
Note that $(\Vi,\Ei)$ is the trivial Toom graph defined in (\ref{triv}) in the case of $D^{\rm triv}_0$, and the trivial incomplete Toom graph defined in (\ref{trivinc}) in the case of $D^{\rm triv}_1,\ldots,D^{\rm triv}_m$.

The following definition defines a random variable $\Dec$ with values in $\Di'_0\cup\{\dgg\}$. Its law depends on the parameters used in the previous subsection, i.e., the constant $\hat p\in[0,1]$, the probability distributions $\hat{\rbf}_t$, and the subprobability distributions $\abf^\star_{s,k}$, and also on an additional parameter $p_\circ\in[0,1]$ that determines how frequently we add a new source.

\begin{defi}
The\label{D:Markont} \emph{contour extending Markov chain} is the Markov chain $(\Dec_n)_{n\geq 0}$ with state space $\Di'_0\cup\{\dgg\}$, initial law
\be\label{Deinit}
\P\big[\Dec_0=D^{\rm triv}_0\big]=\hat p
\quand
\P\big[\Dec_0=D^{\rm triv}_k\big]=(1-\hat p)\hat\rbf_1(k)\qquad(1\leq k\leq m),
\ee
and the following transition mechanism. If $\Dec_{n-1}$ is given $(n\geq 1)$, then we construct $\Dec_n$ as follows. If there are no loose ends or if $\Dec_{n-1}=\dgg$, then we set $\Dec_n:=\Dec_{n-1}$. In the opposite case, we find the most urgent loose end $(v,s)$ of $\Dec_{n-1}$. If $1\leq s<\sig$, then we let $\Dec_n$ be the upward extension of $\Dec_{n-1}$ at $(v,s)$, according to Definition~\ref{D:forext}. If $s=\sig$, then with probabilities $1-p_\circ$ and $p_\circ$ we let $\Dec_n$ be the downward or source extension of $\Dec_{n-1}$ at $(v,s)$, according to Definitions \ref{D:downext} and \ref{D:sourext}. We set
\be
\tau:=\inf\big\{n\geq 0:\Dec_n=\dgg\mbox{ or }\Dec_n\mbox{ has no loose ends}\big\},
\ee
with the convention that $\inf\emptyset:=\infty$, and we define a random variable $\Dec$ with values in $\Di_0\cup\{\dgg\}$ by
\be\label{Dec}
\Dec:=\left\{\ba{ll}
\dis\Dec_\tau\quad&\mbox{if }\tau<\infty,\\
\dis\dgg\quad&\mbox{if }\tau=\infty.
\ea\right.
\ee
\end{defi}

Theorem~\ref{T:randcont} follows from the following proposition, which says that for each $D\in\Di_0$, the probability that $\P[\Dec=D]$ is exactly equal to the right-hand side of (\ref{DecD}). This then proves that the right-hand side of (\ref{DecD}) defines a subprobability measure on $\Di_0$.

\bp[The probability of a given contour]
Let\label{P:nuform} $\Dec$ be the random decorated Toom contour defined in (\ref{Dec}). Then
\be
\P\big[\Dec=D\big]=\nu_{\hat p,p_\circ,\hat\rbf,\abf}(D)\qquad(D\in\Di_0),
\ee
where $\nu_{\hat p,p_\circ,\hat\rbf,\abf}$ is defined in (\ref{DecD}).
\ep

\bpro
For brevity, let us write $\nu=\nu_{\hat p,p_\circ,\hat\rbf,\abf}$. It is clear from (\ref{DecD}) and (\ref{Deinit}) that $\P[\Dec=D^{\rm triv}_0]=\hat p=\nu(D^{\rm triv}_0)$ so it remains to show that $\P[\Dec=D]=\nu(D)$ for all nontrivial $D\in\Di_0$. Let $D=(v_\circ,\Vi,\Ei,\psi,\kappa)$ be a nontrivial decorated Toom contour rooted at $0$. We set $N:=|V|-1$ and $M:=|\Ei|$, we define a function $\eta:V\to\{0,\ldots,N\}$ and incomplete Toom graphs
\be
\big(\Vi'(e),\Ei'(e)\big)_{0\leq e\leq M}
\ee
as in Lemma~\ref{L:induct}. We set $(\Vi',\Ei'):=(\Vi'(M),\Ei'(M))$ and $\psi':=\psi\circ\eta^{-1}$. Then $D':=(0,\Vi',\Ei',\psi',\kappa)$ is a decorated Toom contour that is isomorphic to $D$ and we have $\Dec=D$ up to isomorphism if and only if $\Dec$ is precisely equal to $D'$.

By the Markovian nature of our construction, the probability that $\Dec$ is precisely equal to $D'$ is given by the probability that the root gets the right decoration multiplied by the product of the probabilities that in each step $1\leq e\leq M$, we extend the incomplete Toom contour $\big(\Vi'(e-1),\Ei'(e-1)\big)$ in precisely the right way. We claim that this product is precisely the right-hand side of (\ref{DecD}).

Indeed, by Definition~\ref{D:forext}, for the charges $1\leq s<\sig$, each time we add an edge $(v,w)$ to $\vec E^\ast_s$ with $\psi(v,w)=(i,ij)$, we pick up a factor $\abf^\ast_{s,\kappa(i)}$. If the space-time point $ij$ is not yet included in our Toom contour, then in view of (\ref{taudef}) and our definition of the natural order on $V$ in (\ref{natord}), $\tau(i)\in\{\circ,\ast,1,\ldots,\sig\}$ describes the type of the new vertex $w$ we are adding to our incomplete Toom graph. Now if we assign to the new space-time point $ij$ a decoration $\kappa(ij)\in\{0,\ldots,m\}$, then we pick up a factor $\hat p$ for a sink and a factor $(1-\hat p)\hat\rbf_{\tau(ij)}(\kappa(ij))$ for an internal vertex or a source. From the initial condition (\ref{Deinit}), we similarly pick up a factor $(1-\hat p)\hat\rbf_1(\kappa(0))$ for the root to have the correct decoration. Note that this is consistent with the convention introduced in Subsection~\ref{S:rand} that the root has type~1. By Definitions \ref{D:downext} and \ref{D:sourext}, the factors we pick up when we add an edge of charge $\sig$ are the same as for the charges $1\leq s<\sig$, except that we pick up an additional factor $B_\bullet^{-1}(1-p_\circ)$ each time we add an internal vertex or connect to the root, and an additional factor $B_\circ^{-1}p_\circ$ each time we add a source other than the root. Taking into account (\ref{ncirc}), multiplying all factors gives exactly the right-hand side of (\ref{DecD}).
\epro

\noi
\bpro[of Theorem~\ref{T:randcont}]
Immediate from Proposition~\ref{P:nuform} and the remarks above it.
\epro

\section{Evaluation of the Peierls sum}\label{S:evalu}

\subsection{Subprobabilistic comparison}\label{S:subpro}

In this subsection we prove Theorems \ref{T:cycbd} and \ref{T:conbd}. We first introduce some useful notation and make some preliminary remarks. For any decorated Toom contour $D=(v_\circ,\Vi,\Ei,\psi,\kappa)$, we introduce the notation
\be\label{UUU}
U:=\psi(V)\beh\psi(V_\ast),
\ee
and we set
\be
J^\star_s(i):=\big\{j\in\La:(i,ij)\in\psi(\vec E^\star_s)\big\}\qquad\big(\star\in\{\circ,\bullet\},\ \lis,\ i\in U).
\ee
We note that by (\ref{ncirc}), we have $n_\ast(D)=n_\circ(D)+1$ and
\be\label{ncirc2}
\big|\big\{(i,ij):i\in U,\ j\in J^\star_s(i)\big\}\big|=n_\star(D)\qquad\big(\lis,\ \star\in\{\circ,\bullet\}\big).
\ee
We observe that since internal vertices of a given charge do not overlap with each other (by Definition~\ref{D:embed}~(ii)), for each $i\in U$ and $\lis$, the set $J^\bullet_s(i)$ has at most one element. We use this to define
\be\label{Si}
S(i):=\big\{s\in\{1,\ldots,\sig\}:J^\bullet_s(i)\neq\emptyset\big\}
\quand J^\bullet_s(i)=:\{j_s(i)\}\qquad\big(i\in U,\ s\in S(i)\big).
\ee
With the notation we have just introduced, we can rewrite formula (\ref{PiL}) for the presence-probability $\pi_{p,\rbf}(D)$ as
\be\label{piform}
\pi_{p,\rbf}(D)=p^{n_\ast(D)}\,(1-p)^{n_\dia(D)}\underbrace{\Big(\prod_{i\in U}\rbf\big(\kappa(i)\big)\prod_{s\in S(i)}\ex{-\la L_s\big(j_s(i)\big)}\Big)}_{\txt=:F_\bullet}\underbrace{\Big(\prod_{i\in U}\prod_{s=1}^\sig\prod_{j\in J^\circ_s(i)}\ex{-\la L_s(j)}\Big)}_{\txt=:F_\circ}.
\ee
The definitions of the constants $\al^\star_{s,k}$ and $\bet_{s,k}$ in (\ref{alpha}) and (\ref{beta}) allow us to estimate
\bc\label{FS}
\dis F_\bullet&\leq&\dis\prod_{i\in U}\bet_{\tau(i),\kappa(i)}\hat\rbf_{\tau(i)}\big(\kappa(i)\big)\prod_{s\in S(i)}\al^\bullet_{s,\kappa(i)}\abf^\bullet_{s,\kappa(i)}\big(j_s(i)\big),\\[5pt]
\dis F_\circ&\leq&\dis\prod_{i\in U}\prod_{s=1}^\sig\prod_{j\in J^\circ_s(i)}\al^\circ_{s,\kappa(i)}\abf^\circ_{s,\kappa(i)}(j).
\ec
If we are only interested in Toom cycles, then because of Definition~\ref{D:cycpres}~(iii), we can replace the constants $\al^\circ_{s,k}$ by the slightly better constants $\ti\al^\circ_{s,k}$ from (\ref{varalpha}). We need to compare the expression in (\ref{piform}) to $\nu_{\hat p,p_\circ,\hat\rbf,\abf}$ defined in (\ref{DecD}), which using our new notation can be written as
\bc\label{nuform}
\dis\nu_{\hat p,p_\circ,\hat\rbf,\abf}(D)
&=&\dis\Big(\prod_{i\in U}\hat\rbf_{\tau(i)}\big(\kappa(i)\big)\prod_{s\in S(i)}\abf^\bullet_{s,\kappa(i)}\big(j_s(i)\big)\Big)
\Big(\prod_{i\in U}\prod_{s=1}^\sig\prod_{j\in J^\circ_s(i)}\abf^\circ_{s,\kappa(i)}(j)\Big)\\[25pt]
&&\dis\times\hat p^{n_\ast(D)}\,(1-\hat p)^{n_\dia(D)}\, B_\circ^{-n_\circ(D)}B_\bullet^{-n_\bullet(D)}p_\circ^{n_\circ(D)}(1-p_\circ)^{n_\bullet(D)}.
\ec

\bpro[of Theorem~\ref{T:cycbd}]
We will prove that under the assumptions of Theorem~\ref{T:cycbd}
\be\label{nucycbd}
\pi_{p,\rbf}(D)\leq\frac{p}{\hat p(1-\hat p)}\nu_{\hat p,p_\circ,\hat\rbf,\abf}(D)\qquad(p\leq\eps)
\ee
for all $D\in\bar\Di_0$, the set of decorated Toom cycles rooted at~$0$. Summing over $\bar\Di_0$, using the fact that $\nu_{\hat p,p_\circ,\hat\rbf,\abf}$ from Theorem~\ref{T:randcont} is a subprobability measure, the claim of Theorem~\ref{T:cycbd} then follows from Theorem~\ref{T:Peierls}, Lemma~\ref{L:Pei}, and the definition of $\ov\rho(p)$ in (\ref{ovrho}).

Under the assumptions of Theorem~\ref{T:cycbd}, we claim that for any decorated Toom cycle $D=(v_\circ,\Vi,\Ei,\psi,\kappa)$,
\bc\label{Fcyc}
F_\bullet&\leq&\dis\prod_{i\in U}\hat\rbf_{\tau(i)}\big(\kappa(i)\big)\prod_{s\in S(i)}\bet_{s,\kappa(i)}\al^\bullet_{s,\kappa(i)}\abf^\bullet_{s,\kappa(i)}\big(j_s(i)\big)\\[5pt]
&\leq&\dis\big(\prod_{s=1}^\sig\ga^\bullet_s\big)^{n_\bullet(D)}\prod_{i\in U}\hat\rbf_{\tau(i)}\big(\kappa(i)\big)\prod_{s\in S(i)}\abf^\bullet_{s,\kappa(i)}\big(j_s(i)\big).
\ec
Indeed, we recall from~\eqref{cycord} that as a result of Definition~\ref{D:cycle}~(ii) and Lemma~\ref{L:cycord}, we have $\tau(i)=1$ if and only if $1\in S(i)$. Using this, (\ref{FS}), and the fact that $\bet_{t,k}=1$ if $t\neq 1$ (as $\hat\rbf_\circ=\hat\rbf_2=\rbf$ by assumption), we obtain the first inequality in (\ref{Fcyc}). The second inequality follows from the definition of $\ga^\bullet_s$ in (\ref{betcyc}) and (\ref{ncirc2}). Similarly, but more simply, we obtain from (\ref{FS}) that
\be
F_\circ\leq\big(\prod_{s=1}^\sig\ga^\circ_s\big)^{n_\circ(D)}\prod_{i\in U}\prod_{s=1}^\sig\prod_{j\in J^\circ_s(i)}\abf^\circ_{s,\kappa(i)}(j).
\ee
Inserting these estimates into (\ref{piform}), throwing away the factors $(1-p)$, recalling the definition of $B_\ast$ from~\eqref{BB} and comparing with (\ref{nuform}), we see that
\bc
\dis\pi_{p,\rbf}(D)&\leq&\dis p^{n_\ast(D)}\big(B_\bullet\prod_{s=1}^\sig\ga^\bullet_s\big)^{n_\bullet(D)}\big(B_\circ\prod_{s=1}^\sig\ga^\circ_s\big)^{n_\circ(D)}\\[5pt]
&&\times\hat p^{-n_\ast(D)}\,(1-\hat p)^{-n_\dia(D)}p_\circ^{-n_\circ(D)}(1-p_\circ)^{-n_\bullet(D)}\nu_{\hat p,p_\circ,\hat\rbf,\abf}(D).
\ec
In the factors involving $B_\bullet$ and $B_\circ$ we recognise the constants $C_\bullet$ and $C_\circ$ defined in Theorem~\ref{T:cycbd}. Using the fact that $n_\ast(D)=n_\circ(D)+1$ and (\ref{dici}), we obtain
\be\label{lapi}
\pi_{p,\rbf}(D)\leq\frac{p}{\hat p(1-\hat p)}\Big((1-p_\circ)^{-1}(1-\hat p)^{-\sig}C_\bullet\Big)^{n_\bullet(D)}\Big(\frac{p}{\hat p}p_\circ^{-1}(1-\hat p)^{-1}C_\circ\Big)^{n_\circ(D)}\nu_{\hat p,p_\circ,\rbf,\abf}(D).
\ee
Using the assumption that $p_\circ:=1-(1-\hat p)^{-\sig}C_\bullet>0$, we see that the factor with the power $n_\bullet(D)$ in (\ref{lapi}) equals one. Provided that
\be
p\leq\eps:=\hat p(1-\hat p)p_\circ C_\circ^{-1},
\ee
the factor with the power $n_\circ(D)$ in (\ref{lapi}) can also be estimated from above by one, and (\ref{nucycbd}) follows.
\epro

\bpro[of Theorem~\ref{T:conbd}]
We will prove that under the assumptions of Theorem~\ref{T:conbd}
\be\label{nuconbd}
\pi_{p,\rbf}(D)\leq\frac{p}{\hat p(1-\hat p)}\nu_{\hat p,p_\circ,\hat\rbf,\abf}(D)\qquad(p\leq\eps)
\ee
for all $D\in\Di_0$, the set of decorated Toom contours rooted at~$0$. Summing over $\Di_0$, using the fact that $\nu_{\hat p,p_\circ,\hat\rbf,\abf}$ from Theorem~\ref{T:randcont} is a subprobability measure, the claim of Theorem~\ref{T:conbd} then follows from Theorem~\ref{T:Peierls}, Lemma~\ref{L:Pei}, and the definition of $\ov\rho(p)$ in (\ref{ovrho}).

If $\tau(i)\neq\circ$, then we can use the constants $\bet^{s'}_{s,k}$ from (\ref{bet3}) to estimate
\be\ba{l}
\dis\bet_{\tau(i),\kappa(i)}\hat\rbf_{\tau(i)}\big(\kappa(i)\big)\prod_{s\in S(i)}\al^\bullet_{s,\kappa(i)}\abf^\bullet_{s,\kappa(i)}\big(j_s(i)\big)\\[5pt]
\dis\quad=\Big(\prod_{s'=1}^\sig\bet^{s'}_{\tau(i),\kappa(i)}\Big)\hat\rbf_{\tau(i)}\big(\kappa(i)\big)\prod_{s\in S(i)}\al^\bullet_{s,\kappa(i)}\abf^\bullet_{s,\kappa(i)}\big(j_s(i)\big)\\[5pt]
\dis\quad\leq\hat\rbf_{\tau(i)}\big(\kappa(i)\big)\prod_{s\in S(i)}\bet^s_{\tau(i),\kappa(i)}\al^\bullet_{s,\kappa(i)}\abf^\bullet_{s,\kappa(i)}\big(j_s(i)\big)\\[5pt]
\dis\quad\leq\hat\rbf_{\tau(i)}\big(\kappa(i)\big)\prod_{s\in S(i)}\bet^s_{\kappa(i)}\al^\bullet_{s,\kappa(i)}\abf^\bullet_{s,\kappa(i)}\big(j_s(i)\big).
\ec
Here in the first inequality we have used that if $\tau(i)\neq\circ$, then $\tau(i)\in S(i)$, as well as the fact that $\bet^{s'}_{s,k}\leq 1$ for $s\neq s'$, and in the second inequality the definition of $\beta_k^s$ from~\eqref{betde}. When $\tau(i)=\circ$, we use the facts that $\hat\rbf_\circ=\rbf$ and hence $\bet_{\circ,k}=1$ to estimate more simply
\be\ba{l}
\dis\bet_{\tau(i),\kappa(i)}\hat\rbf_{\tau(i)}\big(\kappa(i)\big)\prod_{s\in S(i)}\al^\bullet_{s,\kappa(i)}\abf^\bullet_{s,\kappa(i)}\big(j_s(i)\big)\\[5pt]
\dis\quad\leq\Big(\prod_{s=1}^\sig(1\wedge\bet^s_{\kappa(i)})\Big)^{-1}\hat\rbf_{\tau(i)}\big(\kappa(i)\big)\prod_{s\in S(i)}\bet^s_{\kappa(i)}\al^\bullet_{s,\kappa(i)}\abf^\bullet_{s,\kappa(i)}\big(j_s(i)\big).
\ec
Using the definitions of $\de$ and $\ga^\bullet_s$ in (\ref{betde}) and (\ref{betgast}) as well as (\ref{FS}), we can then estimate
\bc
F_\bullet &\leq&\dis\de^{-n_\circ(D)}\prod_{i\in U}\hat\rbf_{\tau(i)}\big(\kappa(i)\big)\prod_{s\in S(i)}\ga^\bullet_s\abf^\bullet_{s,\kappa(i)}\big(j_s(i)\big)\\[5pt]
&=&\dis\de^{-n_\circ(D)}\big(\prod_{s=1}^\sig\ga^\bullet_s\big)^{n_\bullet(D)}\prod_{i\in U}\hat\rbf_{\tau(i)}\big(\kappa(i)\big)\prod_{s=1}^\sig\prod_{j\in J^\bullet_s(i)}\abf^\bullet_{s,\kappa(i)}(j).
\ec
The rest of the proof is identical to the proof of Theorem~\ref{T:cycbd}.
\epro

\subsection{Toom's rule}\label{S:some}

In this subsection, we prove Theorem~\ref{T:NEC}.\med

\bpro[of Theorem~\ref{T:NEC}.]
We apply Theorem~\ref{T:conbd} with $\La:=\Z^3$, $m:=1$, and $\phh_1=\phh$ the map defined in (\ref{Toomrule}). The height function is $h(i):=i_1+i_2+i_3$. We use the $\La$-linear polar function of dimension $\sig:=3$ defined by
\be
L_1(i):=-2i_1+i_2+i_3,\quad L_2(i):=i_1-2i_2+i_3,\quad L_3(i):=i_1+i_2-2i_3\qquad(i\in\La),
\ee
and for the sets $A_s=A_{s,1}\in\Ai(\phh_1)$ $(\lis)$ we choose
\be
A_1:=\big\{(0,0,1),(0,1,0)\big\}\quad A_2:=\big\{(0,0,1),(1,0,0)\big\}\quad A_3:=\big\{(0,1,0),(1,0,0)\big\},
\ee
which means that the set $\De=\De_1$ from assumption~(ii) of Subsection~\ref{S:abstract} is given by
\be
\De:=\bigcup_{s=1}^3A_s=\big\{(0,0,1),(0,1,0),(1,0,0)\big\}.
\ee
We set $\tet:=e^{-3\la}$ with $\la$ the constant from assumption~(i) of Subsection~\ref{S:abstract}. This has the effect that the functions $e^{-\la L_s}$ $(\lis)$ take the following values on $\De$:

\begin{center}
\inputtikz{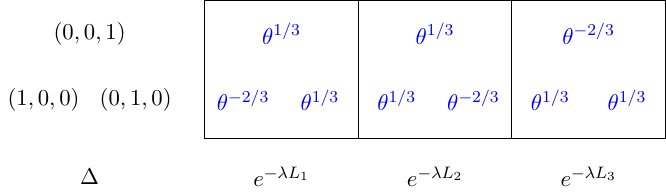}
\end{center}

\noi
We choose for $\abf^\bullet_s=\abf^\bullet_{s,1}$ the uniform distribution on $A_s$ $(\lis)$, which results in the constants $\al^\bullet_s$ from (\ref{alpha}) taking the values
\be
\al^\bullet_1=\al^\bullet_2=\al^\bullet_3=2\tet^{1/3}.
\ee
For $\abf^\circ_s=\abf^\circ_{s,1}$ $(\lis)$ we choose the measure on $\De$ defined by normalising the function $e^{-\la L_s}$ so that it becomes a probability measure. The result is that the constants $\al^\circ_s$ from (\ref{alpha}) are given by
\be
\al^\circ_1=\al^\circ_2=\al^\circ_3=\tet^{-2/3}(1+2\tet).
\ee
Since there is no intrinsic randomness, $\rbf$ and the measures $\hat\rbf_t$ are all equal to the trivial probability distribution on $\{1\}$. This means that the constants $\bet_{t,k}$ from (\ref{beta}) are all equal to one, so we can take all constants $\bet^{s'}_{s,k}$ from (\ref{bet3}) equal to one which has the effect that the constants $\bet^s_k$ and $\de$ from (\ref{betde}) are equal to one. Since there is no $k$-dependence the constants $B_\star$ from (\ref{BB}) are also equal to one and (\ref{betgast}) simplifies to
\be\label{Csimp}
C_\star=\prod_{s=1}^\sig\al^\star_s\qquad\big(\star\in\{\circ,\bullet\}\big).
\ee
This yields
\be
C_\bullet=8\tet\quand C_\circ=\tet^{-2}(1+2\tet)^3.
\ee
Choosing $\tet<1/8$ allows us to choose $\hat p>0$ such that $p_\circ:=1-(1-\hat p)^{-\sig}C_\bullet>0$. Theorem~\ref{T:conbd} then implies that
\be\label{ovlow}
\ov\rho(p)\geq 1-\frac{p}{\hat p(1-\hat p)}\qquad(p\leq\eps),
\ee
with
\be
\eps:=\hat p(1-\hat p)p_\circ C_\circ^{-1}=\frac{\hat p(1-\hat p)\big(1-8\tet(1-\hat p)^{-3}\big)}{\tet^{-2}(1+2\tet)^3}.
\ee
At this point, we resort to numerical calculations, which tell us that setting $\hat p:=0.14$ and $\tet:=0.05$ is approximately optimal. Using these values, we obtain $\eps\approx 8.3928\cdot 10^{-5}$ which is a little bit larger than $1/12000$. Setting $p=1/12000$ and $\hat p=0.14$ in (\ref{ovlow}) gives a lower bound on $\ov\rho(p)$ of approximately $0.9993$, so we conclude that $p_{\rm c}\geq 1/12000$.
\epro

\subsection{Majority rules on the triangular lattice}\label{S:tripro}

In this subsection we show how Theorem~\ref{T:conbd} implies Theorem~\ref{T:triang}.\med

\bpro[of Theorem~\ref{T:triang}]
We have $\La:=\Z^3$, $m:=3$, $\rbf$ is the uniform distribution on $\{1,2,3\}$, and the local maps $\phh_1=\phh^{\rm SWC},\phh_2=\phh^{\rm SUC},\phh_3=\phh^{\rm WUC}$ the maps defined in~(\ref{phitri}). Let $\phi_1$, $\phi_2$, and $\phi_3$ be the spatial maps associated with $\phh_1$, $\phh_2$, and $\phh_3$ in the sense of (\ref{phhphi}), and let $\Ll$ be the linear polar function with dimension $\sig=3$ defined by
\be\label{triL}
\Ll_1(z):=-z_1-z_2,\quad \Ll_2(z):=2z_1-z_2,\quand \Ll_3(z):=2z_2-z_1\quad(z=(z_1, z_2)\in\R^2).
\ee
It is straightforward to check that the edge speeds defined in (\ref{edge}) are given by
\be\label{triedge}
\eps_{\phi_k}(\Ll_s):=\sup_{A\in\Ai(\phi_k)}\inf_{i\in A}\Ll_s(i)=\left\{\ba{ll}
1\quad&\mbox{if }s=k,\\[5pt]
0\quad&\mbox{if }s\neq k.
\ea\right.
\ee
For $1\leq s,k\leq 3$, we let ${\rm A}_{s,k}$ be the unique element of $\Ai(\phi_k)$ for which the supremum in (\ref{triedge}) is attained. In a picture (compare Figure~\ref{fig:triang}), these are the following sets:

\begin{center}
\inputtikz{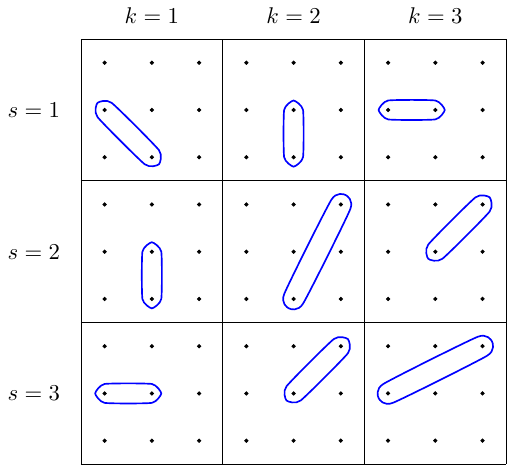}
\end{center}

\noi
We lift these sets to space-time by setting $A_{s,k}:=\big\{(i_1,i_2,-1):(i_1,i_2)\in\Aa_{s,k}\big\}$ $(1\leq s,k\leq 3)$ and we lift the spatial polar function $\Ll$ to space-time by setting $L(i_1,i_2,i_3):=\Ll(i_1,i_2)$ $(i=(i_1,i_2,i_3)\in\La)$, which in the language of (\ref{drift}) means that we choose the drift $v$ to be zero. Setting $\tet:=e^{-\la}$, we observe that the functions $e^{-\la L_s}$ take the following values on $A_{s,k}$:

\begin{center}
\inputtikz{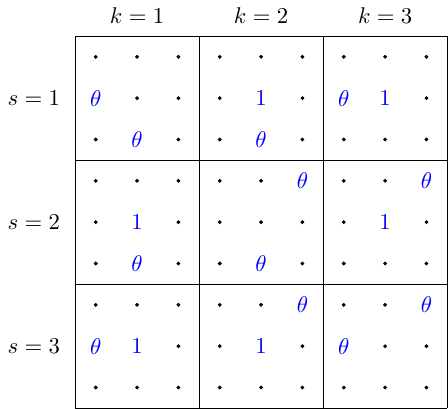}
\end{center}

\noi
This motivates us to choose the subprobability distributions $\abf^\bullet_{s,k}$ as follows.

\begin{center}
\inputtikz{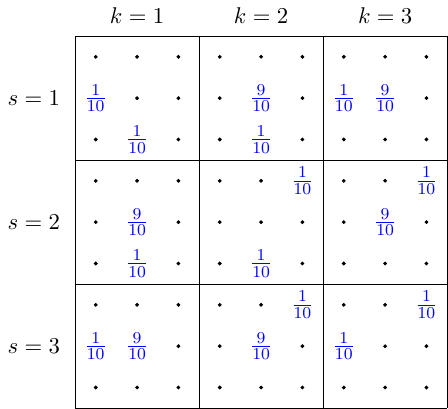}
\end{center}

\noi
Note that for $s=k$ these are strict subprobability distributions in the sense that they sum up to a value strictly less than one. The motivation for this is that we want to keep the constant $B_\bullet$ from (\ref{BB}) small. Indeed, four our choice of $\abf^\bullet_{s,k}$, we have
\be\label{Bbul}
B_\bullet=\ffrac{11}{10}.
\ee
For the moment, we do not care yet about the precise choice of the subprobability distributions $\abf^\circ_{s,k}$ and the precise value of $B_\circ$.

By choosing $\la$ large enough we can make $\tet$ as small as we wish. In particular, we can make $\tet$ small enough such that $10\,\tet\leq 10/9$, so that the constants $\al^\bullet_{s,k}$ from (\ref{alpha}) are given by
\be\label{ask}
\al^\bullet_{s,k}=\left\{\ba{ll}
10\,\tet\quad&\mbox{if }s=k,\\[5pt]
\ffrac{10}{9}\quad&\mbox{if }s\neq k.
\ea\right.
\ee
We define probability distributions $\hat\rbf_s$ $(1\leq s\leq 3)$ by
\be
\hat\rbf_s(k):=\left\{\ba{ll}
\ffrac{1}{9}\quad&\mbox{if }s=k,\\[5pt]
\ffrac{4}{9}\quad&\mbox{if }s\neq k.
\ea\right.
\ee
Since $\rbf(1)=\rbf(2)=\rbf(3):=1/3$, it follows that the constants $\bet_{s,k}$ from (\ref{beta}) are given by
\be\label{rhatr}
\bet_{s,k}=
\left\{\ba{ll}
3\quad&\mbox{if }s=k,\\[5pt]
\ffrac{3}{4}\quad&\mbox{if }s\neq k.
\ea\right.
\ee
As in (\ref{bet3}), we want to write $\bet_{s,k}=\prod_{s'=1}^3\bet^{s'}_{s,k}$ with $\bet^{s'}_{s,k}\leq 1$ if $s\neq s'$. The philosophy behind this is as follows: if $\tet$ is small enough, then formula (\ref{ask}) tells us that for each charge $s$, there is one ``good'' decoration $k=s$ for which $\al^\bullet_{s,k}<1$ while the other two decorations $k\neq s$ are ``bad'' in the sense that $\al^\bullet_{s,k}>1$. We want to multiply these ``bad'' values with a factor $\bet^s_k$ so that their product is $<1$. In view of these considerations, we choose the constants $\bet^s_{s',k}$ with $1\leq s,s',k\leq 3$ as follows:
\be\label{bsk}
\bet^{s'}_{s,k}:=
\left\{\ba{ll}
\dis(\ffrac{3}{4})^{1/2}\quad&\mbox{if }s'\neq k,\\[5pt]
\dis 3(\ffrac{3}{4})^{-1}\quad&\mbox{if }s'=k=s,\\[5pt]
\dis 1\quad&\mbox{if }s'=k\neq s.
\ea\right.
\ee
Note that $\bet_{s,k}=\prod_{s'=1}^3\bet^{s'}_{s,k}$ and $\bet^s_{s',k}\leq 1$ if $s\neq s'$ as required in (\ref{bet3}). The constants in (\ref{betde}) are now given by
\be\label{eq:delmaj}
\bet^s_k:=\sup_{1\leq s'\leq 3}\bet^s_{s',k}=\left\{\ba{ll}
\dis(\ffrac{3}{4})^{1/2}\quad&\mbox{if }s\neq k,\\[5pt]
\dis 4\quad&\mbox{if }s=k,
\ea\right.
\quad\quand\de:=\inf_{1\leq k\leq 3}\prod_{s=1}^3(1\wedge\bet^s_k)=\ffrac{3}{4}
\ee
We can now calculate the constants in (\ref{betgast}). We do not care about the precise values of $\ga^\circ_s$ and $C_\circ$. Assuming that $\tet$ is small enough so that $40\tet\leq(\ffrac{3}{4})^{1/2}\ffrac{10}{9}$, we see that the other two constants are given by
\be
\ga^\bullet_s:=\sup_{1\leq k\leq 3}\bet^s_k\al^\bullet_{s,k}=(\ffrac{3}{4})^{1/2}\ffrac{10}{9}\quad(1\leq s\leq 3)
\ee
and, using (\ref{Bbul}),
\be
C_\bullet:=B_\bullet\prod_{s=1}^\sig\ga^\bullet_s=\ffrac{11}{10}(\ffrac{3}{4})^{3/2}(\ffrac{10}{9})^3=\ffrac{275}{243}\sqrt{3/4}\approx0.98<1.
\ee
Since $C_\bullet<1$, we can choose $\hat p>0$ such that $p_\circ:=1-(1-\hat p)^{-\sig}C_\bullet>0$, and use the inequality (\ref{conbd}) of Theorem~\ref{T:conbd} to conclude that $\ov\rho(p)\to 1$ as $p\to 0$, which means that the cellular automaton of Theorem~\ref{T:triang} is stable.

To also prove the bound $p_{\rm c}\geq 7.7\cdot 10^{-13}$ we have to work a bit more. Modifying our previous argument, it will be convenient to choose the subprobability distributions $\abf^\bullet_{s,k}$ a bit differently. We choose:

\begin{center}
\inputtikz{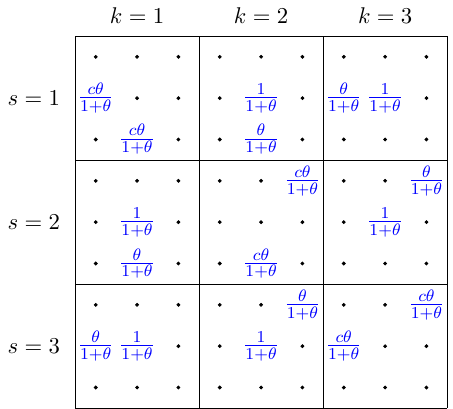}
\end{center}

\noi
Here $c>1$ is a constant to be chosen later. The effect of this is that
\be\label{ask2}
B_\bullet=\frac{1+2c\tet}{1+\tet}\quand\al^\bullet_{s,k}=\left\{\ba{ll}
(1+\tet)/c\quad&\mbox{if }s=k,\\[5pt]
1+\tet\quad&\mbox{if }s\neq k.
\ea\right.
\ee
Although this is probably not optimal, for simplicity, we stick to our previous choice of the measures $\hat\rbf_s$, which means that also the constants $\bet^s_k$ and $\de$ remain unchanged. Then
\be
\bet^s_k\al^\bullet_{s,k}=\left\{\ba{ll}
4(1+\tet)/c\quad&\mbox{if }s=k,\\[5pt]
(\ffrac{3}{4})^{1/2}(1+\tet)\quad&\mbox{if }s\neq k.
\ea\right.
\ee
This motivates us to choose
\be
c:=4(\ffrac{3}{4})^{-1/2},
\ee
which has the effect that
\be
\ga^\bullet_s=(\ffrac{3}{4})^{1/2}(1+\tet)\quad(1\leq s\leq 3)
\ee
and
\be\label{cbullmaj}
C_\bullet=\frac{1+2c\tet}{1+\tet}\prod_{s=1}^3\ga^\bullet_s=(1+16\cdot 3^{-1/2}\tet)(\ffrac{3}{4})^{3/2}(1+\tet)^2.
\ee
We also have to make a choice for the subprobability distributions $\abf^\circ_{s,k}$ on $\De_k=\bigcup_{s=1}^3{\rm A}_{s,k}$. We set $Z:=1+\tet+\tet^{-2}$ and choose $\abf^\circ_{s,k}$ as follows:

\begin{center}
\inputtikz{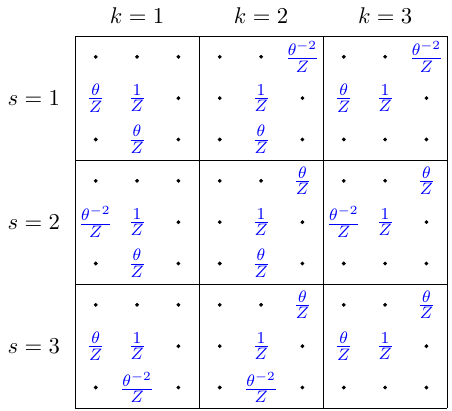}
\end{center}

\noi
The consequence of this is that $\al^\circ_{s,k}=Z$ for all $1\leq s,k\leq 3$ hence $\ga^\circ_s=Z$ $(1\leq s\leq 3)$, while $B_\circ$ defined in (\ref{BB}) is given by
\be
B_\circ=\frac{1+2\tet+\tet^{-2}}{1+\tet+\tet^{-2}}.
\ee
Using the fact that $\de=3/4$ by \eqref{eq:delmaj}, it follows that
\be
C_\circ:=\de^{-1}B_\circ\prod_{s=1}^3\ga^\circ_s
=\frac{4}{3}\cdot \frac{1+2\tet+\tet^{-2}}{1+\tet+\tet^{-2}}Z^3
=\ffrac{4}{3}(1+2\tet+\tet^{-2})(1+\tet+\tet^{-2})^2.
\ee
We have from \eqref{cbullmaj} that $C_\bullet<1$ for $\tet$ small enough which allows us to choose $\hat p>0$ such that $p_\circ:=1-(1-\hat p)^{-\sig}C_\bullet>0$. Theorem~\ref{T:conbd} then allows us to conclude that
\be
\ov\rho(p)\geq 1-\frac{p}{\hat p(1-\hat p)}\qquad(p\leq\eps),
\ee
with 
\bc
\dis\eps&:=&\dis\hat p(1-\hat p)p_\circ C_\circ^{-1}\\[5pt]
&=&\dis\hat p(1-\hat p)\frac{1-(1-\hat p)^{-3}(1+16\cdot 3^{-1/2}\tet)(\ffrac{3}{4})^{3/2}(1+\tet)^2}{\ffrac{4}{3}(1+2\tet+\tet^{-2})(1+\tet+\tet^{-2})^2}.
\ec
Numerically, it seems that the maximal possible value of $\eps$ is attained approximately in $\tet=0.033$ and $\hat p=0.016$. For these values we find $\eps>7.7\cdot 10^{-13}$ and $\ov\rho(\eps)>1-5\cdot 10^{-11}$, which allows us to conclude that $p_{\rm c}\geq 7.7\cdot 10^{-13}$.
\epro

\subsection{Cooperative branching and the identity map}\label{S:coopro}

In this subsection we derive Theorem~\ref{T:coop} from Theorem~\ref{T:cycbd}.\med

\bpro[of Theorem~\ref{T:coop}]
We have $\La=\Z^3$, $m:=2$, $\phh_1:=\phh^{\rm coop}$ and $\phh_2:=\phh^{\rm id}$ as defined in (\ref{coopid}), and $\rbf$ is the probability distribution given by $\rbf(1):=r$ and $\rbf(2):=1-r$. Let $\phi_1$ and $\phi_2$ be the spatial maps associated with $\phh_1$ and $\phh_2$ in the sense of (\ref{phhphi}), and let $\Ll$ be the linear polar function with dimension $\sig=2$ defined by
\be\label{Lcoopp}
\Ll_1(z):=z_1+z_2\quand\Ll_2(z):=-z_1-z_2\qquad\big(z=(z_1,z_2)\in\R^2\big).
\ee
Then the edge speeds defined in (\ref{edge}) are
\be\ba{l@{\quad}l}
\dis\eps_{\phi_1}(\Ll_1)=1,
&\dis\eps_{\phi_1}(\Ll_2)=0,\\[5pt]
\dis\eps_{\phi_2}(\Ll_1)=0,
&\dis\eps_{\phi_2}(\Ll_2)=0.
\ec
For $1\leq s,k\leq 3$, we let $\Aa_{s,k}$ be the unique element of $\Ai(\phi_k)$ such that (compare (\ref{edge}))
\be
\eps_{\phi_k}(\Ll_s)=\inf_{i\in\Aa_{s,k}}\Ll_s(i).
\ee
Concretely, these are the sets
\be
\Aa_{1,1}=\big\{(1,0),(0,1)\big\},\quad \Aa_{1,2}=\Aa_{2,1}=\Aa_{2,2}=\big\{(0,0)\big\}.
\ee
We lift these sets to space-time by setting $A_{s,k}:=\big\{(i_1,i_2,-1):(i_1,i_2)\in\Aa_{s,k}\big\}$ $(1\leq s,k\leq 2)$ and we lift the spatial polar function $\Ll$ to space-time by setting $L(i_1,i_2,i_3):=\Ll(i_1,i_2)$ $(i=(i_1,i_2,i_3)\in\La)$, which in the language of (\ref{drift}) means that we choose the drift $v$ to be zero.

We choose for $\abf^\bullet_{1,1}$ and $\abf^\circ_{2,1}$ the uniform distribution on $A_{1,1}$ and we choose for all other $\abf^\star_{s,k}$'s the trivial distribution that gives probability one to $(0,0)$. The constants $B_\star$ from (\ref{BB}) are then given by
\be
B_\bullet=1\quand B_\circ=2.
\ee
Setting $\tet:=e^{-\la}$, one can then check that the constants $\al^\bullet_{s,k}$ and $\ti\al^\circ_{s,k}$ from (\ref{alpha}) and (\ref{varalpha}) are given by
\be
\al^\bullet_{1,1}=2\tet,\quad\ti\al^\circ_{2,1}=2\tet^{-1},\quand\al^\bullet_{s,k}=\ti\al^\circ_{s,k}=1\quad\mbox{otherwise}.
\ee
In line with Theorem~\ref{T:cycbd} we set $\hat\rbf_\circ=\hat\rbf_2=\rbf$. We define $\hat\rbf_1(1):=\hat r$ and $\hat\rbf_1(2):=1-\hat r$, which implies that the constants $\bet_{s,k}$ from (\ref{beta}) are given by
\be
\bet_{1,1}=\frac{r}{\hat r},\quad\bet_{1,2}=\frac{1-r}{1-\hat r},\quand\bet_{t,k}=1\quad\mbox{otherwise}.
\ee
It follows that the constants from (\ref{betcyc}) are given by
\be\ba{ll}
\dis\ga^\bullet_1=2\tet\frac{r}{\hat r}\vee\frac{1-r}{1-\hat r},\quad
&\dis\ga^\bullet_2=1,\\[5pt]
\dis\ga^\circ_1=1,\quad
&\dis\ga^\circ_2=2\tet^{-1},
\ec
and as a result
\be
C_\bullet=2\tet\frac{r}{\hat r}\vee\frac{1-r}{1-\hat r}\quand C_\circ=4\tet^{-1}.
\ee
We choose $\hat r$ with the aim of making $C_\bullet$ as small as possible. The optimal choice is the unique root of the equation $2\tet(r/\hat r)=(1-r)/(1-\hat r)$, which after a little algebra gives
\be
C_\bullet=1-r+2\tet r.
\ee
By choosing $\la$ large enough we can make $\tet=e^{-\la}$ as small as we wish. In particular, we can make sure that $\tet<1/2$. Then $C_\bullet<1$, which allows us to choose $\hat p>0$ such that $p_\circ:=1-(1-\hat p)^{-2}C_\bullet>0$. Theorem~\ref{T:cycbd} then allows us to conclude that the density of the upper invariant law satisfies
\be\label{rholow}
\ov\rho(p,r)\geq 1-\frac{p}{\hat p(1-\hat p)}\qquad(p\leq\eps),
\ee
with $\eps:=\hat p(1-\hat p)p_\circ C_\circ^{-1}$. In particular, this shows that $\lim_{p\to 0}\ov\rho(p,r)=1$ for all $r>0$, proving part~(i) of Theorem~\ref{T:coop}.

To also prove part~(ii) we assume that $\hat p=dr$ for some constant $d$ to be chosen later. A small calculation shows that
\be
\eps=\ha\tet d(\ha-\tet-d)r^2+O(r^3)\quad\mbox{as }r\to 0.
\ee
The optimal values of $d$ and $\tet$ are $d=\tet=\ffrac{1}{6}$, which gives
\be
\eps=\ha(\ffrac{1}{6})^3r^2+O(r^3)\quad\mbox{as }r\to 0.
\ee
Using this and $\hat p=dr=\ffrac{1}{6}r$, (\ref{rholow}) tells us that when $p\to 0$ and $r\to 0$ in such a way that $p\leq cr^2$, with $c<\ha(\ffrac{1}{6})^3$, then
\be
\ov\rho(p,r)\geq 1-\frac{cr^2}{\ffrac{1}{6}r(1-\ffrac{1}{6}r)}
\ee
for all $r$ small enough, so $\ov\rho(p,r)\to 1$.
\epro

\section{Shrinkers}\label{S:shrink}

\subsection{Compensated edge speeds}\label{S:compspeed}

In this section, we study shrinkers. In the present subsection, we prove Lemma~\ref{L:compspeed} which is closely related to an observation of Gray, who proved that each eroder is ``shift-equivalent'' to a shrinker \cite[Thm~18.2.1]{Gra99}. In Subsection~\ref{S:diverge} we prove Proposition~\ref{P:Peinf} which shows that Conjecture~\ref{C:shrink} about shrinkers cannot be proved with the methods of the present article.

We will use the following auxiliary result~\cite[Prop.1.7]{Zieg95}.

\bl[Farkas lemma] 
Let \label{L:Farkas} $S$ be a finite set and let $\Ll_s: \R^d\to\R$ be linear functions and $\eps_s\in\R$ for all $s\in S$. Then precisely one of the following two alternatives must hold:
\begin{itemize}
\item[{\rm I.}] there exists a $z\in\R^d$ such that $\Ll_s(z)\geq \eps_s$ for all $s\in S$,
\item[{\rm II.}] there exist constants $\la_s\geq 0$ $(s\in S)$ such that $\sum_{s\in S}\la_s\Ll_s(z)=0$ for all $z\in\R^d$ and $\sum_{s\in S}\la_s\eps_s>0$.
\end{itemize}
\el

\bpro[of Lemma~\ref{L:compspeed}]
Let $\Ll:\R^d\to \R^\sig$ be the linear polar function from Theorem~\ref{T:thm9} and $\eps_s:=\inf_{1\leq k\leq m}\eps_{\phi_k}(\Ll_s)$ $(\lis)$. We first choose the linear polar function $\Ll':\R^d\to \R^{\sig'}$.

Let $H_{s}:=\{z\in\R^d:\Ll_s(z)\geq \eps_s\}$. If there exists $z\in\cap_{s=1}^\sig H_s$, then $\sum_{s=1}^\sig \Ll_s(z)\geq \sum_{s=1}^\sig \eps_s>0$, contradicting the definition of a polar function. Therefore, $\bigcap_{s=1}^\sig H_s=\emptyset$. Let $\sigma'$ be the minimal number of these open halfspaces whose intersection is empty, we can assume without loss of generality that these halfspaces are $\{H_s, \; 1\leq s\leq \sig'\}$. If $\sig=\sig'$, we set $\Ll':=\Ll$. Assume that $\sigma'<\sigma$. We now use Lemma~\ref{L:Farkas} with $S=\{1, \dots, \sig'\}$. As $\bigcap_{s=1}^{\sig'} H_s=\emptyset$, condition II. must hold, that is, there exist constants $\lambda_1, \dots, \lambda_{\sig'}\geq 0$ such that
\be
\sum_{s=1}^{\sig'} \lambda_s \Ll_s(z)=0 \quad (z\in\R^d).
\ee
Assume that $\lambda_{\sig'}=0$. We can then apply Lemma~\ref{L:Farkas} again for the set $S'=\{1, \dots, \sig'-1\}$ and use the constants $\lambda_1, \dots, \la_{\sig'-1}$ to show that condition II. holds. Then we must have $\bigcap_{s=1}^{\sig'-1} H_s=\emptyset$, contradicting the minimality of the set of halfspaces whose intersection is empty. Therefore, we must have $\lambda_s>0$ for each $1\leq s \leq \sig'$.
Letting $\Ll'_s(z):=\lambda_s \Ll_s(z)$ $(1\leq s\leq\sig')$ it is then clear that $\Ll'$ is a polar function of dimension $\sig'\geq 2$. Furthermore, setting $\eps'_s:=\inf_{1\leq k\leq m}\eps_{\phi_k}(\Ll'_s)$ we claim that for each $1\leq s \leq \sig'$.
\be\label{eq:halfplaneequiv}\{z\in\R^d:\Ll'_s(z)\geq \eps'_s\}=\{z\in\R^d:\Ll_s(z)\geq \eps_s\}=H_s.\ee
This is indeed true, as the linearity of $\Ll$ and $\lambda_s>0$ implies that $\eps_{\phi_k}(\Ll'_s)=\lambda_s \eps_{\phi_k}(\Ll_s)$, and therefore $\eps'_s=\lambda_s\eps_s$.

Having fixed the polar function $\Ll'$ we now show that there exists $v\in\R^d$ such that $\eps_{\phi_k}^v(\Ll'_s)>0$ for all $1\leq s \leq\sig'$ and $1 \leq k \leq m$. As $\Ll'$ is linear, we have $\Ll'_s(i-v)=\Ll'_s(i)-\Ll'_s(v)$ for all $i\in\R^d$, therefore $\eps^v_{\phi_k}(\Ll'_s)=\eps_{\phi_k}(\Ll'_s)-\Ll'_s(v)$ for all $s$ and $k$. We then have
\be \eps^v_{\phi_k}(\Ll'_s)>0 \quad \forall 1\leq s \leq\sig', 1 \leq k \leq m \quad \Leftrightarrow \quad \inf_{1\leq k\leq m}\eps^v_{\phi_k}(\Ll'_s)=\eps'_s-\Ll'_s(v)>0 \quad \forall 1\leq s \leq\sig'.
\ee
That is, we need to prove that there exists a $v\in\R^d$ such that $\Ll'_s(v)<\eps'_s$ for all $s$, or equivalently by \eqref{eq:halfplaneequiv}, such that $v\in\bigcap_{s=1}^{\sig'} H_{s}^c$, where $H_s^c:=\R^d\setminus H_s$ denotes the complement of $H_s$. The rest of the proof goes along the same lines as the proof of \cite[Lemma~18.2.2.]{Gra99}.

We first show that for each $1\leq s'\leq \sig'$ there exists a $v_{s'}\in\R^d$ such that
\be\label{eq:vs}
\Ll'_{s'}(v_{s'})<\eps_{s'}, \quad \Ll'_{s}(v_{s'}) = \eps_{s} \text{ for any }s\neq s'.
\ee
Fix $s'$ and let $\partial H_{s}:=\{z\in\R^d: \Ll_s(z) = \eps_s\}$ denote the boundary of the closed halfspace $H_s$. As $\{H_s, \; 1\leq s\leq \sig'\}$ is a minimal set of halfspaces whose intersection is empty, we have that $I_{s'}:=\cap_{s\neq s'}H_s \neq \emptyset$. Furthermore, $I_{s'}$ is a closed convex set as it is the intersection of finitely many closed halfspaces. Let $\bar I_{s'}$ denote the set of points $z\in I_{s'}$ whose Eucledian distance to $\partial H_{s'}$ is minimal. As $\bigcup_{s=1}^\sig H_{s}^c=\R^d$, we have $I_{s'}\subset H^c_{s'}$. This implies that $\bar I_{s'}$ has a positive distance from $\partial H_{s'}$, is an affine linear subspace of $\R^d$ and is parallel to $\partial H_{s'}$. Therefore, there must exist a maximal set $\Sigma'$ of halfspaces such that
\be
\bar I_{s'} = \bigcap_{s\in\Sigma'}\partial H_s.
\ee
We now show that $\Sigma'=\{1, \dots, \sigma'\}\setminus\{s'\}$. For any $x\in \bar I_{s'}$ and $z\in \bigcap_{s\in\Sigma'} H_s$, every point of the line segment connecting $z$ to $x$ lies at least as far away from $\partial H_{s'}$ as $x$ itself. Since $\bar I_{s'}\subset H^c_{s'}$, the ray starting at $x$ and pointing in the direction of $z-x$ must also be entirely contained in $H^c_{s'}$. This implies that $\bigcap_{s\in\Sigma'} H_s$ is entirely contained in $H^c_{s'}$ as well, so that $\bigcup_{s\in \Sigma'\cup\{s'\}}H_s^c=\R^d$. As $\{H^c_s, \; 1\leq s\leq \sig'\}$ is a minimal set of halfspaces whose union covers $\R^d$ we then must have $\Sigma'=\{1, \dots, \sigma'\}\setminus\{s'\}$.

As $\bar I_{s'} = \bigcap_{s\neq s'}\partial H_s$ and $\bar I_{s'}\subset H^c_{s'}$, for any $v_{s'}\in\bar I_{s'}$ \eqref{eq:vs} holds. Having chosen $v_{s'}$ for all $1\leq s'\leq \sig'$, we set $v:=1/ \sig' \cdot (v_1+\dots+v_{\sig'})$.
Then indeed 
\be
\Ll'_s(v)=\frac 1 {\sig'} (\Ll'_s(v_1) + \dots + \Ll'_{\sig'}(v_{\sig'})) < \eps_s \quad \forall 1\leq s \leq \sig',
\ee
concluding the proof.
\epro

\subsection{Infinity of the Peierls sum}\label{S:diverge}

In this subsection, we prove Proposition~\ref{P:Peinf}. We have $\La=\Z^3$, $m:=2$ and we consider the monotone cellular automaton $\Phi^{p, r}$ that applies the maps $\phh^{\rm cc}$ and $\phh^{\rm id}$ (whose associated ``space-maps'', $\phi^{\rm cc}$ and $\phi^{\rm id}$, are defined in (\ref{phicc})) with probabilities $(1-p)r$ and $(1-p)(1-r)$, respectively.

We will use decorated Toom cycles of a special form. The height function is $h(i_1,i_2,i_3)=-i_3$. Given an integer $n\geq 1$ and a function $f:\{0,\ldots,n\}\to\Z$ with $f_0=0=f_n$ and $|f_k-f_{k-1}|\leq 1$ for all $1\leq k\leq n$, we can define a Toom cycle $(V',\psi)$ of length $2n+4$ by:
\begin{itemize}
\item $\psi_0=(0,0,0)$ and $\psi_k=(0,1,-k)$ for $1\leq k\leq n+1$,
\item $\psi_{n+2}=(0,0,-n)$, $\psi_{n+3}=(0,0,-n-1)$,
\item $\psi_{2n+4-k}=(f_k,-f_k,-k)$ for $0\leq k\leq n$.
\end{itemize}
Here $V'=\{0,\ldots,2n+4\}$ and, in line with Definition~\ref{D:cyc},
\be
V'_\circ=\{n+2\},\quad V'_\ast=\{n+1,n+3\},\quad V'_1=\{0,\ldots,n\},\quad V'_2=\{n+4,\ldots,2n+4\}.
\ee
Note that $\psi_{n+2}=\psi_{n+4}=(0,0,-n)$ is the only space-time point (apart from the origin) that is visited twice and that $\psi_{n+2}\in V'_\circ$ and $\psi_{n+4}\in V'_2$, in line with part~(ii) of Definition~\ref{D:cycle}.

To make $(V',\psi)$ into a decorated Toom cycle $(V',\psi,\kappa)$ in the sense of Definition~\ref{D:decor}, we let $g:\{0,\ldots,n\}\to\{1,2\}$ be a function such that $g_0=g_n=1$ and
\be
g_k=1\quad\mbox{if}\quad\nab{k}f\in\{-1,+1\}\quad\mbox{with}\quad\nab{k}f:=f_{k+1}-f_k\quad(0\leq k<n).
\ee
We then define a decoration function $\kappa:\psi(V')\to\{0,1,2\}$ by
\begin{itemize}
\item $\kappa(f_k,-f_k,-k)=g_k$ for all $0\leq k\leq n$,
\item $\kappa(0,1,-k)=2$ for all $1\leq k\leq n$,
\item $\kappa(0,0,-n-1)=\kappa(0,1,-n-1)=0$.
\end{itemize}
Taking into account the definition of the sets $\Aa_{s,k}$ in (\ref{Acount}), one can check that $\kappa$ satisfies conditions (i), (ii), and (iii)' of Definition~\ref{D:decor} so $(V',\psi,\kappa)$ is a decorated Toom cycle. We let $\hat\Di_0(n)$ denote the set of all decorated Toom cycles rooted at $0$ that are of this special form.

As a first step, we derive a stochastic expression for $\sum_{D\in\hat\Di_0(n)}\pi_{p,\rbf}(D)$, where $\pi_{p,\rbf}(D)$ is the presence-probability of $D$ as defined in Definition~\ref{D:decor}, with $\rbf(1):=r$ and $\rbf(2):=1-r$. Let $(\chi_k,\eta_k)_{k\geq 0}$ be i.i.d.\ such that
\be
\P\big[(\chi_k,\eta_k)=(l,1)\big]=\ffrac{1}{3}r\quad(l=-1,0,1)
\quand
\P\big[(\chi_k,\eta_k)=(0,2)\big]=1-r,
\ee
and let
\be\label{Mn}
S_n:=\sum_{k=0}^{n-1}\chi_k\quand M_n:=\sum_{k=0}^{n-1}1_{\{\eta_k=1\}}\qquad(n\geq 1).
\ee

\bl[Sum over special contours]
For\label{L:stochex} all $n\geq 1$, one has
\be\label{stochex}
\sum_{D\in\hat\Di_0(n)}\pi_{p,\rbf}(D)=p^2(1-p)^{2n+1}r^2(1-r)^n\E\big[3^{M_n}1_{\{S_n=0\}}\,\big|\,\chi_0=1\big].
\ee
\el

\bpro
There is a natural one-to-one correspondence between $\hat\Di_0(n)$ and the set of all functions $(\nab{k}f,g_k)_{0\leq k\leq n-1}$ that take values in the set $\{(-1,1),(0,1),(1,1),(0,2)\}$ and that satisfy $g_0=1$ and $\sum_{k=0}^{n-1}\nab{k}f=0$. Let $M_n$ be as defined in (\ref{Mn}) with $\eta_k$ replaced by $g_k$. Then the cardinality of $\{0\leq k\leq n:g_k=1\}$ is $M_n+1$ since also $g_n=1$. We have $|\psi(V')|=2n+3$ and there are two sinks, so the presence-probability $\pi_{p,\rbf}(D)$ as defined in (\ref{prespro}) of a contour $D$ that is defined in terms of $(\nab{k}f,g_k)_{0\leq k\leq n-1}$ is given by 
\be
\pi_{p,\rbf}(D) = p^2(1-p)^{2n+1}\big(\prod_{k=0}^n\rbf(g_k)\big)\big(\prod_{k=1}^n\rbf(2)\big)
=p^2(1-p)^{2n+1}r^{M_n+1}(1-r)^{2n-M_n}.
\ee
On the other hand,
\be
\P\big[(\chi_k,\eta_k)_{0\leq k\leq n-1}=(\nab{k}f,g_k)_{0\leq k\leq n-1}\big]=3^{-M_n}r^{M_n}(1-r)^{n-M_n},
\ee
so
\be
\pi_{p,\rbf}(D)=p^2(1-p)^{2n+1}r(1-r)^n3^{M_n}\P\big[(\chi_k,\eta_k)_{0\leq k\leq n-1}=(\nab{k}f,g_k)_{0\leq k\leq n-1}\big].
\ee
Summing over all functions $(\nab{k}f,g_k)_{0\leq k\leq n-1}$ with $g_0=1$ and $\sum_{k=0}^{n-1}\nab{k}f=0$, we find that
\be
\sum_{D\in\hat\Di_0(n)}\pi_{p,\rbf}(D)=p^2(1-p)^{2n+1}r(1-r)^n\E\big[3^{M_n}1_{\{S_n=0\}}1_{\{\chi_0=1\}}\big],
\ee
which by the fact that $\P[\eta_0=1]=r$ can be rewritten as (\ref{stochex}).
\epro

The next lemma shows that if $3^r(1-r)(1-p)^2>0$, then the expression in (\ref{stochex}) grows exponentially in~$n$.

\bl[Asymptotic lower bound]
For\label{L:aslow} all $p,r\in(0,1]$, one has
\be
\liminf_{n\to\infty}\Big(\sum_{D\in\hat\Di_0(n)}\pi_{p,\rbf}(D)\Big)^{1/n}\geq 3^r(1-r)(1-p)^2
\ee
\el

\bpro
We use Lemma~\ref{L:stochex}. Basic large deviation theory tells us that for each $\eps>0$, there exists a $\la_\eps>0$ such that $\P[M_n\leq(r-\eps)n]\leq e^{-\la_\eps n}$ for all $n$ large enough. Moreover, by the local central limit theorem, there exists a constant $c>0$ such that $\P[S_n=0]\geq cn^{-1/2}$ for all $n$ large enough. It follows that
\be
\sum_{D\in\hat\Di_0(n)}\pi_{p,\rbf}(D)\geq p^2(1-p)r^2\,\big(cn^{-1/2}-e^{-\la_\eps n}\big)\,\big[3^{(r-\eps)}(1-r)(1-p)^2\big]^n
\ee
for all $n$ large enough. Since $\eps>0$ is arbitrary, this implies the claim.
\epro

\bpro[of Proposition~\ref{P:Peinf}]
By Lemma~\ref{L:Pei} the Peierls sum in (\ref{Peinf}) is equal to $\sum_{D\in\bar\Di_0}\pi_{p,\rbf}(D)$ and can therefore be estimated from below by $\sum_{n=1}^\infty\sum_{D\in\hat\Di_0(n)}\pi_{p,\rbf}(D)$. As an immediate consequence of Lemma~\ref{L:aslow}, we see that for each $r\in(0,1)$ such that $3^r(1-r)>1$, there exists a $p'\in(0,1]$ such that $\sum_{n=1}^\infty\sum_{D\in\hat\Di_0(n)}\pi_{p,\rbf}(D)=\infty$ for all $0<p\leq p'$. Since $\log[3^r(1-r)]=r\log 3-r+O(r^2)$ as $r\to 0$ and $\log 3>1$, there exists an $r'\in(0,1]$ such that $3^r(1-r)>1$ for all $0<r\leq r'$, and Proposition~\ref{P:Peinf} follows.
\epro

\end{document}